\documentclass[10pt]{amsart}

\usepackage{amssymb}
\usepackage{pifont}
\usepackage{amscd}
\usepackage{amsmath}
\usepackage[dvips]{graphicx}
\usepackage{dsfont}
\usepackage[all]{xy}
\usepackage{color}
\usepackage{mathrsfs}
\usepackage{ulem}

\usepackage[colorlinks=true, linkcolor=blue]{hyperref}

\usepackage{tikz}
\usetikzlibrary{arrows.meta}

\theoremstyle{definition}

\newtheorem{Thm}{Theorem}[section]
\newtheorem{Lem}[Thm]{Lemma}
\newtheorem{Def}[Thm]{Definition}
\newtheorem{Cor}[Thm]{Corollary}
\newtheorem{Prop}[Thm]{Proposition}
\newtheorem{Prop-Def}[Thm]{Proposition-Definition}
\newtheorem{Ex}[Thm]{Example}
\newtheorem{Rem}[Thm]{Remark}

\newtheorem{Cri}[Thm]{Criterion}

\newcommand{\comp}{\mathop{\raisebox{+.3ex}{\hbox{$\scriptstyle\circ$}}}}

\newcommand{\xto}{\xrightarrow}

\numberwithin{equation}{section}

\title[BV structure]{The Hochschild cohomlogy ring  of a self-injective Nakayama algebra is a Batalin-Vilkovisky  algebra}
\author{Xiuli  Bian, Tomohiro Itagaki,   Wen Kou,  Weiguo Lyu     and Guodong Zhou }

\address{Xiuli Bian, Wen Kou,  Guodong Zhou
\newline School of Mathematical Sciences
\newline Key Laboratory of MEA (Ministry of Education)
\newline  Shanghai Key laboratory of PMMP
\newline  East China Normal University
\newline 200241, Shanghai
\newline P. R. China}
\email{52205500018@stu.ecnu.edu.cn}
\email{52285500004@stu.ecnu.edu.cn}
\email{gdzhou@math.ecnu.edu.cn}

\address{Weiguo Lyu
\newline Department of Mathematics  
\newline Shanghai Maritime University
\newline 201306, Shanghai
\newline P. R. China}
\email{wglv@shmtu.edu.cn}

\address{Tomohiro Itagaki 
\newline Department of Economics  
\newline Faculty of Economics 
 \newline Takasaki City University of Economics
 \newline Gunma 370-0801
\newline Japan}
\email{titagaki@tcue.ac.jp }

\newenvironment{SProof}[1][Sketch of Proof]{\begin{trivlist}
\item[\hskip \labelsep {\bfseries #1}]}{\flushright
$\Box$\end{trivlist}}

\newenvironment{Proof}[1][Proof]{\begin{trivlist}
\item[\hskip \labelsep {\bfseries #1}]}{\flushright
$\Box$\end{trivlist}}

\newcommand{\pl}{ \parallel }

\newcommand{\K}{{\mathbb K}}
\newcommand{\N}{{\mathbb N}}
\newcommand{\bbB}{{\mathbb B}}
\newcommand{\bbP}{{\mathbb P}}

\newcommand{\ot}{\otimes}
\newcommand{\di}{{\; | \;}}
\newcommand{\ul}{\underline}
\newcommand{\ol}{\overline}
\newcommand{\mode}{\;  (\mathrm{mod}\; e)}
\newcommand{\model}{\;  (\mathrm{mod}\; e_1)}
\newcommand{\lra}{\longrightarrow}

\newcommand{\rmH}{\mathrm{H}}
\newcommand{\rmHH}{\mathrm{HH}}
\newcommand{\Hom}{\mathrm{Hom}}
\newcommand{\rmIm}{\mathrm{Im}}
\newcommand{\Ker}{\mathrm{Ker}}
\newcommand{\rmchar}{\mathrm{char}}

\newcommand{\calB}{{\mathcal B}}
\newcommand{\calD}{{\mathcal D}}
\newcommand{\calL}{{\mathcal L}}

\newcommand{\End}{\mathrm{End}}

\newcommand{\Id}{\mathrm{Id}}
\newcommand{\ie}{\mathrm{i.e.}}

\newcommand{\rmBar}{{\mathrm{Bar}}}

\newcommand{\even}{\mathrm{even}}

\newcommand{\ra}{\rightarrow}
\newcommand{\sdp}{\times\kern-.2em\vrule height1.1ex depth-.05ex}
\newcommand{\epi}{\lra \kern-.8em\ra}

\newcommand{\Z}{{\mathbb Z}}

\newcommand{\calE}{{\mathcal E}}

 \normalbaselines
\begin{document}

\begin{abstract}

Lambre, Zhou and Zimmermann showed that the Hochschild cohomology ring of a Frobenius algebra with
semisimple Nakayama automorphism is a Batalin-Vilkovisky
algebra. They asked whether the semisimplicity condition is necessary. In this paper, we show that for a self-injective Nakayama algebra, the Hochschild cohomology ring  is always a Batalin-Vilkovisky algebra.

In course of proofs, we correct some inaccuracies in the literature, hoping not to introduce new errors.

\end{abstract}

\renewcommand{\thefootnote}{\alph{footnote}}
	\setcounter{footnote}{-1} 
	\renewcommand{\thefootnote}{\alph{footnote}}
	\setcounter{footnote}{-1} 

\maketitle

\tableofcontents

\section*{Introduction}

  Hochschild cohomology plays a central role in the study of associative algebras, it carries a rich algebraic structures: it is a graded commutative algebra structure via the cup product(also known as the Yoneda product), a graded Lie algebra structure via a Lie bracket of degree $-1$, and thus a Gerstenhaber algebra structure \cite{Ger63}, among others.

In recent decades, a new algebraic structure - known as a Batalin-Vilkovisky (BV) structure - on Hochschild cohomology has been studied. Roughly speaking a BV structure consists of a degree $-1$ operator that squares to zero and interacts with the cup product, to recover the Lie bracket. Such a structure is known to exist only for certain special classes of algebras.

Tradler first showed that the Hochschild cohomology ring of a finite dimensional symmetric algebra is a BV algebra \cite{Tra08}. Later, Lambre, Zhou and Zimmermann \cite{LZZ16}, and independently Volkov \cite{Vol16} discovered that Hochschild cohomology ring of a Frobenius algebra  with semisimple Nakayama automorphism has a BV algebra structure. 
Lambre, Zhou and Zimmermann further asked whether this result holds for Frobenius algebras with nonsemisimple  Nakayama automorphisms. Due to the computational complexity of Hochschild cohomology, this question  is rather hard.

When this paper was in course of  preparation, Herscovich and Li showed that for the Fomin-Kirillov algebra with three generators, which is a Frobenious algebra, its Hochschild cohomology ring is NOT a BV algebra \cite{HL22}.
Hence, the answer to the question of Lambre-Zhou-Zimmermann is negative in general.

In this  paper, we consider BV structure on the Hochschild cohomology ring of a  self-injective Nakayama algebra.  
In terms of quiver with relations,  basic self-injective Nakayama algebras over an algebraically closed field  are given as  truncated basic cycle algebras, which are truncated quotients of path algebras where the underlying graph is a single oriented cycle.      
 Let us introduce some notations:
Let $Z_e$ denote the quiver with $e$ vertices $\{ 1,\dots, e\}$ and $e$ arrows $\{\alpha_1,\dots,\alpha_e\}$ such that the origin $o(\alpha_i)=i$ of the arrow $\alpha_i$ equals the target $t(\alpha_{i-1})$ of the arrow $\alpha_{i-1}$ with $2\le i\le e$, and $o(\alpha_1)=1=t(\alpha_e)$.
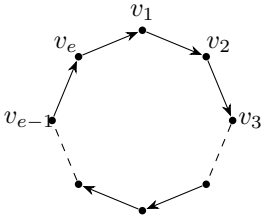
\begin{center}
\begin{tikzpicture}[>=Stealth, scale=1.2]
  \def \n {8} 
  \def \radius {1} 
  \def \angle {360/\n}
  \foreach \s in {1,...,\n} {
    \node[fill=black, circle, inner sep=1pt] (V\s) at ({\angle * (\s - 1)}:\radius) {};
  }
  \node at ({\angle * (0)}:\radius + 0.2) {$v_3$};
  \node at ({\angle * (1)}:\radius + 0.2) {$v_2$};
  \node at ({\angle * (2)}:\radius + 0.2) {$v_1$};
  \node at ({\angle * (4)}:\radius + 0.25) {$v_{e-1}$};
  \node at ({\angle * (3)}:\radius + 0.2) {$v_e$};
  \foreach \s [evaluate=\s as \next using {int(mod(\s, \n) + 1)}] in {1,...,\n} {
    \ifnum\s=\n
      \draw[-, dashed] (V\next) -- (V\s); 
    \else
      \ifnum\s=5
        \draw[-, dashed] (V\next) -- (V\s); 
      \else
        \draw[->] (V\next) -- (V\s); 
      \fi
    \fi
  }
\end{tikzpicture}
\end{center}
Moreover, let $\Lambda=\K Z_e/J^N$, where $N\ge 2$ and $J$ is the arrow ideal of $\K Z_e$. 
This class of algebras is one of the most studied classes of finite dimensional algebras. Zhang \cite{Zha97}, Locateli  \cite{Loc99} as well as Erdmann and Holm \cite{EH99} gave the Hochschild cohomology groups of self-injective Nakayama algebras. Erdmann and Holm \cite{EH99}, Bardzell, Locateli and Marcos \cite{BLM00} determined the cup product on their Hochschild cohomology groups, and
   Xu and Zhang  computed the Gerstenhaber Lie bracket \cite{XZ11}. For a self-injective Nakayama algebra $A$ over an algebraically closed field, Itagaki  \cite{Ita23}  determined the BV algebra structure on the Hochschild cohomology of $\rmHH^*(A)^{\nu \uparrow}$ in the sense of Volkov \cite{Vol16}, which coincides with the usual Hochschild cohomology ring when the   Nakayama automorphism is semisimple.  This  paper deals with self-injective Nakayama algebras regardless whether the   Nakayama automorphism is semisimple or not.


Our main result in this paper confirms the question of Lambre-Zhou-Zimmermann for self-injective Nakayama algebras. 
\begin{Thm} [see Theorem~\ref{thm:main-result}]
    Let $\Lambda=\K Z_e/ J^N$, $N\ge 2$, be a truncated basic cycle algebra over a field $\K$. The the Hochschild cohomology ring $\rmHH^*(\Lambda)$ of $\Lambda$ is a BV algebra.
\end{Thm}

Our method is computational,  one  byproduct  of which is that the BV   algebra structure is   determined  explicitly    in the case  that the Nakayama automorphism is not semisimple.  To this end, we need to review and revise all the results obtained in the literature about Hochschild cohomology of selfinjective Nakayama algebras,  including Hochschild cohomology spaces, cup product formulae, Hochschild cohomology ring in terms generators and relations, Gerstenhaber brackets.  During  this process, some inaccuracies in the literature are found and corrected. 
Based on the Gerstenhaber algebra structure computed  explicitly,  for the nonsemisimple case, a general criterion (Criterion~\ref{Criterion: nonsemisimple case}) is provided, which might be of independent interest, and which  gives an explicit construction of the BV operator.

\medskip

This paper is organized as follows. In Section~\ref{Section: preliminaries}, we review the definitions of Gerstanhaber algebras and BV algebras, as well as the constructions of Gerstenhaber and BV algebra structures on Hochschild cohomology of Frobenius algebras with semisimple Nakayama automorphism. Moreover, we racall the bilinear form and the Nakayama automorphism for self-injective Nakayama algebras. In Section~\ref{section: comparison morphism}, we recall the minimal projective bimodule resolutions for truncated quiver algebras. Moreover, we recall the comparison morphisms between the minimal resolution and the reduced bar resolution. In Section~\ref{Section:Hochschild cohomology groups}, we provide an explicit $\K$-basis for each Hochschild cohomology group of truncated basic cycle algebras, following Locateli \cite{Loc99}, and Bardzell, Locateli and Marcos \cite{BLM00}. In Section~\ref{Section: cup product}, a new product on $\bbP^*$ is described, still called cup product. This turns out to be a well-defined product on $\rmHH^*(\Lambda)$. We also explain why this product coincides with the Yoneda product. 
In Section~\ref{section: Hochschild cohomology ring}, we provide the ring structure for $\rmHH^*(\Lambda)$. We show that $\rmHH^*(\Lambda)$ is always finitely generated and provide a presentation of generators and relations.
In Section~\ref{Sect: Gerstenhaber}, we begin by describing the Gerstenhaber bracket on the Hochschild cohomology of truncated quiver algebras using the language of parallel paths. We then provide a precise characterization of the Gerstenhaber algebra structure of the Hochschild cohomology of truncated basic cycle algebras. In Section~\ref{section:BV-structure}, we explicitly compute the BV algebra structure on $\rmHH^*(\Lambda)$ using the result of \cite{LZZ16}, in the case where the Nakayama automorphism is semisimple.  A general criterion (Criterion~\ref{Criterion: nonsemisimple case}) is presented,   which enables dealing with the case of  nonsemisimple  Nakayama automorphisms.   
Combining the two cases, we obtain that $\rmHH^*(\Lambda)$ also carries a BV algebra structure regardless whether  the Nakayama automorphism is  semisimple or not.

\bigskip

\textbf{Conventions:} Throughout this paper, let $\K$ be a  field, we consider the path algebra $\Lambda=\K Q/I$ and denote the tensor product $\ot_{\K Q_0}$ over $\K Q_0$ by $\ot$ for simplicity. For a set $X$, the $\K$-vector space generated by $X$ is denoted by $\K X$. Paths will be written from left to right.

Whenever two paths in a path algebra are composable, we keep their concatenation as usual; if they are not composable, the expression is understood to be $0$.

\bigskip

\section{Preliminaries}\label{Section: preliminaries}\

In this section, we recall the definitions and the notations for Hochschild cohomology, Gerstenhaber algebras and Batalin-Vilkovisky (BV for short) algebras. Moreover, following \cite{LZZ16}, we recall the BV differential on Hochschild cohomology of a Frobenius algebra with semisimple Nakayama automorphism.

\subsection{Hochschild (co)homology and Gerstenhaber algebras}\

We begin by recalling the definition of a Gerstenhaber algebra.

\begin{Def}

A \textit{Gerstenhaber algebra} over $\K$ is a graded $\K$-vector space $V^\bullet=\bigoplus\limits_{n\in \N} V^n$ equipped with two linear maps: a cup product
\[
\cup\colon\, V^n \times V^m \rightarrow V^{n+m}, \quad (a, b) \mapsto
a\cup b
\]
and a Lie bracket of degree $-1$
\[
[- ,- ]\colon\, V^n\times V^m\rightarrow V^{n+m-1}, \quad (a, b)\mapsto
[a, b]
\]
such that

\begin{itemize}

\item [(1)] $(V^\bullet,  \cup)$ is a graded commutative algebra, that is, $a\cup b=(-1)^{|a|\cdot |b| }b\cup a$;

\item [(2)] $(V^\bullet, [- ,- ])$ is a graded Lie algebra of degree $-1$, i.e., we have that
\[
[a, b]=-(-1)^{(|a| -1)(|b| -1)}[b, a]
\]
and a graded Jacobi-identity
\[
(-1)^{(|a| -1)(|c| -1)}[[a, b], c]+(-1)^{(|b| -1)(|a| -1)}[[b, c], a]+(-1)^{(|c| -1)(|b| -1)}[[c, a], b]=0;
\]

\item [(3)] the cup product and the Lie bracket satisfy the compatible condition called the
Poisson rule, i.e.,
\[
[a\cup b,  c]= a\cup [b, c]+(-1)^{|b| (|c| -1)}[a, c]\cup b,
\]

\end{itemize}

where  $a, b, c$ are arbitrary homogeneous elements in $V^\bullet$, and $|a| $ denotes the degree of $a$.

\end{Def}

One classic example of a Gerstenhaber algebra is the Hochschild cohomology of a $\K$-algebra. The cohomology theory of associative algebras was first introduced by Hochschild in \cite{Hoc45}. Latter, Gerstenhaber showed in \cite{Ger63} that the Hochschild cohomology ring of a $\K$-algebra carries a Gerstenhaber algebra structure.

Given a $\K$-algebra $A$, its \textit{Hochschild cohomology groups} are defined as 
\[
\rmHH^n(A)\cong \mathrm{Ext}^n_{A^e}(A, A), \quad n \geq 0, 
\]
where $A^e=A \ot A^{op}$ is the enveloping algebra of $A$. There exists a projective resolution of $A$ as $A^e$-module, the so called \textit{standard bar resolution} $\mathrm{Bar}(A)$  which is given by $\mathrm{Bar}_r(A)=A\ot {A}^{\ot r}\ot A$, that is,
\begin{eqnarray*}
\mathrm{Bar}(A) \colon \cdots \to A\ot {A}^{\ot r}\ot A \xto{b_{r}} A\ot {A}^{\ot r-1}\ot A\to \cdots \to A\ot {A}\ot A\xto{b_1}A^{\ot 2}(\stackrel{\mu}{\to} A),
\end{eqnarray*}
where the map $\mu: A\ot A\to A$ is the multiplication of $A$, and
\begin{align*}
    b_r(a_0\ot a_1 \ot\cdots\ot a_r \ot a_{r+1})
    =& a_0a_1\ot {a_2}\ot\cdots\ot{a_r}\ot a_{r+1}\\
    &+ \sum_{i=1}^{r-1} (-1)^{i} a_0 \ot a_1 \ot \cdots \ot a_{i-1} \ot a_i a_{i+1} \ot a_{i+2} \ot \cdots \ot a_r \ot a_{r+1}\\
    &+ (-1)^r a_0 \ot a_1 \ot \cdots \ot a_{r-1} \ot a_r a_{r+1}.
\end{align*}

The \textit{Hochschild cohomology complex} is defined to be $C^*(A)=\Hom_{A^e}({\mathrm{Bar}}(A), A)$. Note that for each $r\ge 0$, 
\[
C^r(A)=\Hom_{A^e}(A\ot A^{\ot r}\ot A, A)\cong\Hom_\K(A^{\ot r}, A),
\]
in particular, $C^0(A)\cong\Hom_\K(\K, A)\cong A$. Thus $C^*(A)$ has the following form:
\[
C^*(A) \colon\, 0 \to A \xto{b^0} \Hom_\K(A, A) \to \cdots \to \Hom_\K(A^{\ot r}, A)
\xto{b^r} \Hom_\K(A^{\ot (r+1)}, A)\to \cdots,
\]
where the differential $b^*$ is given as follows: for $r\ge 0$, a linear map $f: A^{\ot r}\to A$, and $a_1, \dots, a_{r+1}\in A$,
\begin{align*}
    b^r(f) (a_1 \ot \cdots \ot a_{r+1})
=& (-1)^{r+1} a_1 f(a_2 \ot \cdots \ot a_{r+1})\\
&+ \sum_{i=1}^r(-1)^{r-i+1} f(a_1 \ot \cdots \ot a_{i-1} \ot a_i a_{i+1} \ot a_{i+2} \ot \cdots \ot a_{r+1})\\
&+ f(a_1 \ot \cdots \ot a_r) a_{r+1}.
\end{align*}
Note that the differential of $b^r$ here includes an additional sign $(-1)^{r+1}$ compared to the original definition. This modification is made to ensure that the definition of the morphism complex $C^*(A,A)=\Hom_{A^e}({\mathrm{Bar}}(A), A)$ is compatible with the differential graded (dg) setting. Specially, if we regard the algebra $A$ as a dg algebra with trivial differential, and the complex $\rmBar(A)$ as a dg $A$-module, then the differential $b^*$ on $C^*(A,A)$ is defined as follows:
\[
b^r(f)=d_A\comp f-(-1)^{|f| } f\comp b_r=(-1)^{r+1}f\comp b_r,
\]
where we write $| f| :=r$. Due to this modification, the definitions of the cup product and the Gerstenhaber bracket are also adjusted accordingly.


The \textit{cup product} $f\cup g \in C^{n+m}(A)=\Hom_\K(A^{\ot (n+m)}, A)$ for $f\in C^m(A)$ and $g\in C^n(A)$ is given by
\[
(f\cup g)(a_1 \ot \cdots \ot a_{n+m}): =(-1)^{m n} f(a_1 \ot \cdots \ot a_n) \cdot
g(a_{n+1} \ot \cdots \ot a_{n+m}).
\]
It is ready to see that the cup product is associative. By \cite[Page 278, (20)]{Ger63}, the differential is a derivation with respect to the cup product, that is,
\[
b(f\cup g)= b(f)\cup g +(-1)^{|f| } f\cup b(g),\; \forall f, g\in C^*(A),
\]
so  $(C^*(A),\, \cup,\, b)$ is a dg algebra. This also shows this cup product induces a well-defined product on Hochschild cohomology,
\[
\cup \colon\, \rmHH^m(A) \times \rmHH^n(A,A) \lra \rmHH^{n+m}(A),
\]
which is still called cup product. Moreover, by \cite[Theorem 3]{Ger63}, this cup product is graded commutative up to homotopy on the cochain level. As a consequence, the graded $\K$-vector space $\rmHH^*(A)=\bigoplus\limits_{n\ge 0}\rmHH^n(A)$ is a graded commutative algebra.

Furthermore, the \textit{Lie bracket} is defined as follows. Let $f \in C^m(A)$ and $g \in C^n(A)$.
If $m, n\ge 1$, then for $1\le i\le m$, define $f\circ_i g \in C^{n+m-1}(A)$ by
\[
(f\circ_i g)(a_1 \ot \cdots \ot a_{n+m-1}):=f(a_1 \ot \cdots \ot a_{i-1} \ot
g(a_i \ot \cdots \ot a_{i+n-1}) \ot a_{i+n} \ot \cdots \ot a_{n+m-1}),
\]
if $ m\ge 1$ and $n=0$, then $g\in A$ and define
\[
(f \circ_i g)(a_1 \ot \cdots \ot a_{m-1}):=f(a_1 \ot \cdots \ot a_{i-1} \ot g \ot a_i \ot \cdots \ot a_{m-1}),
\]
for any other case, set $f\circ_i g$ to be zero. Now we define
\[
f\circ g :=\sum_{i=1}^m (-1)^{(n-1)(i-1)}f \circ_i g,
\]
and
\[
[f, g] :=f \circ g-(-1)^{(n-1)(m-1)} g \circ f.
\]
It is a long routine verification that the shift $sC^*(A)$ (with $sC^*(A)^n=C^{n+1}(A)$) of the Hochschild cochain complex is a graded Lie algebra. By \cite[Page 280]{Ger63}, the differential is also a derivation with respect to the Gerstenhaber Lie bracket, that is,
\[
b([f, g])=[b(f), g]+(-1)^{|f| -1}[f, b(g)],\;  \forall f,\, g\in C^*(A).
\]
As a consequence,
$(sC^*(A),\, [-,-],\, b)$ is a dg Lie algebra.

There exists a compatibility condition between the cup product and the Lie bracket given by \cite[Corollary 1 to Theorem 5]{Ger63}. Altogether,  Gerstenhaber proved the following important result.

\begin{Thm}[{\cite[Page 267]{Ger63}}]  
The Hochschild cohomology ring $(\rmHH^*(A),\, \cup,\, [-,-])$  is a Gerstenhaber algebra.
\end{Thm}

\bigskip

\subsection{Hochschild cohomology and BV algebras.}\

With the development of Hochschild theory and homological algebra, it has been observed that in many situations a Gerstenhaber algebra often carries an additional operator structure, which leads to the notion of a Batalin-Vilkovisky algebra.

\begin{Def}\label{Def-BV-algebra} 
A \textit{Batalin--Vilkovisky  algebra} (BV algebra) is a Gerstenhaber algebra $(\rmH^\bullet,\, \cup,\, [-,-])$ together with an operator 
$\Delta\colon \rmH^\bullet \rightarrow \rmH^{\bullet-1}$
of degree $-1$ such that $\Delta\comp \Delta=0$ and
\[
[f,\, g]=(-1)^{|f| }(\Delta(f\cup g)- \Delta(f)\cup g-(-1)^{|f| }f\cup \Delta(g))
\]
for homogeneous elements $f, g\in \rmH^\bullet$.
\end{Def}

\begin{Rem} 
Our sign convention for Gerstenhaber algebras and BV algebras follows \cite[Proposition 1.2]{Get94} and \cite[Paragraphe 5.1]{Man99} and are different those in \cite{Ger63} and \cite{Tra08}. For the equivalence between these two sign convention, see \cite[Remark 2.5]{LZ14}.
\end{Rem}

Tradler showed that the Hochschild cohomology algebra of a finite dimensional symmetric algebra is a BV algebra \cite{Tra08}; a different proof was given by Menichi~\cite{Men04}. For a symmetric algebra $A$, the $\Delta$-operator on the Hochschild cohomology corresponds to the Connes' $\mathfrak{B}$-operator on the Hochschild homology via the duality between the Hochschild
cohomology and the Hochschild homology.

Lambre, Zhou and Zimmermann \cite{LZZ16} proved the following result generalizing  Tradler-Menichi's result.

\begin{Thm}[\cite{LZZ16}]  
The Hochschild cohomology ring of a Frobenius algebra with semisimple Nakayama automorphism is a BV algebra.
\end{Thm}

Let us recall some details about the above mentioned result of Lambre, Zhou and Zimmermann \cite{LZZ16}, as we shall need them in the sequel.

Recall (cf e.g. \cite[Section 1.10.1]{Zim14} or \cite{YS11}) that a finite dimensional $\K$-algebra $A$ is called a Frobenius algebra, if there is a non-degenerate associative bilinear form 
$\langle-, -\rangle:\, A\times  A\to \K$. 
Here the associativity means that 
$\langle ab, c\rangle=\langle a, bc\rangle$ 
for all $a,\, b$ and $c$ in $A$. Endow $D(A)=\Hom_\K(A, \K)$, the $\K$-dual of $A$,  with the canonical bimodule structure
\[
(afb)(c)=f(bca),\text{ for }  f\in D(A),\text{ and } a, b, c\in A.
\]
The property of being Frobenius is equivalent to saying  that $D(A)=\Hom_\K(A, \K)$ is isomorphic to $A$ as left or as right modules. It is readily seen that the map 
$a\mapsto \langle a, -\rangle$ 
for $a\in A$ gives an isomorphism of right modules between $A$ and $D(A)$, while the map  $a\mapsto \langle -, a\rangle$ gives the isomorphism of left modules.

For $a\in A$, there exists a unique $\sigma(a)\in A$ such that $\langle a,\  -\rangle=\langle -,\ \sigma(a)\rangle\in D(A)$. 
It is easy to see that $\sigma: A\to A$ is an algebra isomorphism and we call it the \textit{Nakayama automorphism of $A$} (associated to the bilinear form $\langle-, -\rangle$).
We denote by $A_{\sigma}$ the $A$-bimodule which is $A$ as $\K$-vector space with the $A$-bimodule defined by
\[
b\cdot c\cdot a :=bc\sigma(a),\text{ for } b,c\in A,\text{ and } a\in A^{op}.
\]
To define the operator $\Delta$ in \cite{LZZ16}, we need the definition of $\rmH_*(A,A_{\sigma})$. The complex used to compute the Hochschild homology $\rmH^*(A,A_{\sigma})$ is $C_*(A,A_{\sigma})=A_{\sigma}\ot_{A^e} \mathrm{Bar}(A)$. For $r\ge 0$, $C_r(A,A_{\sigma})=A_{\sigma}\ot_{A^e}(A\ot A^{\ot r}\ot A)\simeq A_{\sigma} \ot A^{\ot r}$ and the differential 
\[
d_r:\; C_r(A,A_{\sigma})=A_{\sigma}\ot A^{\ot r}\to C_{r-1}(A,A_{\sigma})=A_{\sigma}\ot A^{\ot r-1}
\]
sends $x\ot a_1\ot \cdots\ot a_r$ to
\[
x\sigma(a_1)\ot a_2\ot\cdots\ot a_r +\sum_{i=1}^{r-1}(-1)^i x\ot a_1\ot\cdots\ot a_{i-1}\ot a_i a_{i+1}\ot a_{i+2}\ot \cdots \ot a_r+(-1)^r a_r x\ot a_1\ot \cdots \ot a_{r-1}.
\]

In \cite{KK14}, Kowalzig and Kr\"{a}hmer defined a $\K$-linear map
\[
B_{\sigma}\colon\, C_r(A,A_\sigma) \lra C_{r+1}(A,A_\sigma) 
\]
by
\[
B_{\sigma}(a_0\ot a_1\ot \cdots\ot a_r)=\sum\limits_{i=1}^{r+1}(-1)^{ir} 1\ot a_{i}\ot \cdots\ot a_r\ot a_0\ot \sigma(a_1)\ot \cdots\ot \sigma(a_{i-1}).
\]
It is easy to check that $B_{\sigma}$ is a chain map satisfying $B_{\sigma}\comp B_{\sigma}=0$, which induces an operator $B_{\sigma}:\; \rmH_r(A,A_{\sigma})\to \rmH_{r+1}(A,A_{\sigma})$.

The isomorphism of $A$-bimodules $A_\sigma\rightarrow D(A),\;a\mapsto\langle  - , a\rangle$ induces an isomorphism of complexes 
\[
\partial:\ D(C_*(A,A_\sigma))\lra C^*(A),
\]
which is given by
\begin{equation*}
\begin{aligned}
D(C_r(A, A_\sigma))
& =\Hom_{\K}(A_{\sigma} \ot_{A^e}(A\ot A^{\ot r} \ot A) ,\K) \\
& \cong \Hom_{A^e}(A \ot A^{\ot r} \ot A,\Hom_\K(A_\sigma,\K)) \\
& \cong \Hom_{A^e} (A \ot A^{\ot r} \ot A, A)\\
& =C^r(A),
\end{aligned}
\end{equation*}
for  $r\ge 0$, which induces an isomorphism $D(\rmH_*(A,A_\sigma))\cong\rmHH^*(A)$ on homology group level. In fact, we can make $\partial$ and its inverse $\partial^{-1}$ explicit as follows:
\begin{equation*}
\begin{aligned}
\partial\colon\, \Hom_\K (A_\sigma\ot A^{\ot r},\K)
&\to\Hom_{A^e}(A \ot A^{\ot r} \ot A,A)\\
	f\quad\quad\quad\quad&\mapsto  \partial(f):\, 1\ot a_1\ot \cdots \ot a_r\ot 1\mapsto \sum\limits_{j=1}^N f(y_j\ot a_1\ot \cdots\ot a_r)x_j,\\
\partial^{-1} \colon\, \Hom_{A^e} (A\ot A^{\ot r}\ot A,A) &\to\Hom_\K(A_\sigma \ot A^{\ot r},\K)\\
	g\quad\quad\quad\quad\quad\quad&\mapsto  \partial^{-1}(g):\, a_0\ot a_1\ot \cdots\ot a_r \mapsto \langle  g(1\ot a_1\ot \cdots\ot a_r\ot 1), a_0 \rangle,
\end{aligned}
\end{equation*}
where $\{ x_i\di  1\le i\le N\}$ is a basis of $A$ and $\{ y_j\di  1\le j\le N \}$ the dual basis such that $\langle x_i,y_j\rangle=\delta_{i,j}$.

\medskip

Recall that for a finite dimensional vector space $V$,  a $\K$-linear transformation $\varphi\in \End_\K(V)$ is  called semisimple, if extended to the algebraic closure $\ol{\K}$ of $\K$, the $\ol{\K}$-linear transformation $\varphi\ot_\K\ol{\K}\in \End_{\ol{\K}}(V\ot_\K\ol{\K})$ is diagonalizable.  When talking about semisimplicity of the Nakayama automorphism of a Frobenius algebra, we always mean   this sense.

In the sequel, we could always assume that \textit{$\K$ is algebraically closed}. In fact, once the proofs and computations have been done over an algebraically closed field, then one can return to a general field  by using  the argument in Section 4 of \cite{LZZ16}. 

In the following, let $A$ be a Frobenius algebra with bilinear $\langle-, -\rangle:\, A\times  A\to \K$ and the semisimple Nakayama automorphism $\sigma$. Following \cite{LZZ16}, we recall the BV differential on $\rmHH^*(A)\simeq \rmHH_{(1)}^*(A)$.

Let $I_A$ be the set of eigenvalues of the automorphism $\sigma$. Then $I_A\subset \K$ because $\sigma$ is semisimple. Fix an eigenvalue $\lambda\in I_A$ and let $A_{\lambda}$ denote the eigenspace associated with $\lambda$. It is straightforward to verify that for $\lambda,\mu\in I_A$, one has $A_{\lambda} \cdot A_{\mu} \subseteq A_{\lambda\mu}$. When $\lambda\mu\not\in I_A$, we adopt the convention that $A_{\lambda\mu}=0$. Denote by $\widehat{I_A}:=<I_A>$ the submonoid of $\K^{\times}$ generated by $I_A$.

For each $\mu\in\widehat{I_A}$ put
\[
C_r^{(\mu)}(A,A_{\sigma}):=\bigoplus_{\mu_i\in I_A,\prod \mu_i=\mu} A_{\mu_0}\ot A_{\mu_1}\ot\cdots\ot A_{\mu_r}.
\]
The differential $d_r:\; C_r(A,A_{\sigma}) \to C_{r-1}(A,A_{\sigma})$ restricts to this subspace and denote its restriction by $d_r^{(\mu)}$, then $(C_*^{(\mu)}(A,A_{\sigma}),d_*^{(\mu)})$ is a subcomplex of $(C_*(A,A_{\sigma}),d_*)$. Denote 
\[
\rmH_r^{(\mu)}(A,A_{\sigma}):= \rmH_r(C_*^{(\mu)}(A,A_{\sigma}),d_*^{(\mu)}).
\]
We hence obtain a vector space homomorphism $\rmH_r^{(\mu)}(A,A_{\sigma})\to \rmH_r(A,A_{\sigma})$.

An analogous definition exists for cohomology. For $\mu\in\widehat{I_A}$, let $C^r_{(\mu)}(A)$ be those Hochschild cochains $\varphi\in C^r(A)$ such that we have $\varphi(A_{\mu_1}\ot\cdots\ot A_{\mu_r})\subset A_{\mu\mu_1\cdots\mu_r}$ for all eigenvalues $\mu_1,\ldots,\mu_r$ of $\sigma$. The restriction is $b^r_{(\mu)}:=b^r\mid _{C^r_{(\mu)}(A)}:\; C^r_{(\mu)}(A)\to C^{r+1}_{(\mu)}(A)$. Put
\[
\rmH^r_{(\mu)}(A,A):=\rmH^r(C^*_{(\mu)}(A),b^*_{(\mu)}).
\]
The subcomplex $(C^*_{(\mu)}(A),b^*_{(\mu)})$ of $(C^*(A),b^*)$ defines a morphism of graded vector spaces $\rmH_{(\mu)}^*(A,A)\to \rmHH^*(A)$.

By \cite[Theorem 4.1]{LZZ16}, the BV differential $\Delta$ on $\rmHH^*(A)$ is defined by the following commutative diagram:
\begin{eqnarray*}
\xymatrix{\rmHH^r(A) \ar[d]_{\simeq} \ar[r]^{\Delta_r}&\rmHH^{r-1}(A)\\
H^r_{(1)}(A,A) \ar[d]_{\partial^{-1}} &H^{r-1}_{(1)}(A,A) \ar[u]_{\simeq}\\
D(H_r^{(1)}(A,A_{\sigma})) \ar[r]^{D(B_{\sigma})} &D(H_{r-1} ^{(1)} (A,A_{\sigma})) \ar[u]_{\partial}}
\end{eqnarray*}
In particular, we have, for $f\in \rmH^r(A,A)$ whose corresponding element in $\rmH^r_{(1)}(A,A)$ is $f'$,
\[
 \Delta(f)(\alpha\ot a_1\ot\cdots\ot a_{r-1}\ot \beta) =\sum_{j=1}^N \sum_{i=1}^r  (-1)^{i(r-1)} \langle f'(a_i \ot \cdots\ot a_{r-1}\ot y_j\ot \sigma(a_1)\ot \cdots\ot \sigma(a_{i-1})), \; 1\rangle \alpha x_j\beta.
\]

\medskip

\bigskip

\section{Minimal resolution for truncated quiver algebras and comparison morphisms}\label{section: comparison morphism}\

In this section, we recall the minimal projective bimodule resolutions for truncated quiver algebras following Bardzell \cite{Bar97} which in fact dealt with general monomial algebras.
Moreover, we recall the comparison morphisms between the minimal resolution and the reduced bar resolution as constructed in \cite{ACT09}.

\medskip

Let $Q=(Q_0, Q_1, o, t)$ be a \textit{quiver}. That is, $Q_0$ is the set of vertices, $Q_1$ the set of arrows, the map $o: Q_1\to Q_0$ gives the origin of an arrow and the map $t: Q_1\to Q_0$ is the target of an arrow. We always assume that $Q$ is a finite quiver, i.e., $Q_0$ and $Q_1$ are finite sets. The path $p=\alpha_1\alpha_2\cdots \alpha_r$ of \textit{length} $l( p) =r$ is defined by the condition that $t(\alpha_i)=o(\alpha_{i+1})$ for all $1\leq i\leq r-1$, and vertices are viewed as paths of length $0$. For $n\geq 0$, $Q_n$ denotes the set of all paths of length $n$ and $Q_{\geq n}$ is the set of all paths with length at least $n$. Denote  by $\K Q$ {the path algebra} of $Q$, that is the space generated by all paths of finite length and endowed with the multiplication given by concatenation of paths.

For an algebra given by a bound quiver $\Lambda=\K Q/I$ with $I\subseteq Q_{\geq 1}$, one usually uses a variant of the normalized bar resolution, which is the \textit{reduced bar resolution} introduced by Gerstenhaber and Schack \cite{GS86} and popularized by Cibils \cite{Cib90}. Let ${\Lambda_+}=\Lambda/(\mathbb{K}Q_0)$. Since $I$ is asked to be contained in $Q_{\geq 1}$, the spaces $KQ_n\; (n\geq 0),\; \Lambda,\;  \Lambda_+$ etc are all $\K Q_0$-bimodules. The reduced bar resolution of $\Lambda=\K Q/I$  has the form:
\[
\bbB_*\colon\; \cdots \to \Lambda \ot {\Lambda_+}^{\ot r} \ot \Lambda \xto{b_r} \Lambda \ot
{\Lambda_+}^{\ot r-1} \ot \Lambda \to \cdots \to \Lambda \ot
{\Lambda_+}\ot \Lambda\xto{b_1} \Lambda \ot \Lambda (\xto{b_0} \Lambda),
\]
where $b_0: \Lambda \ot \Lambda \to \Lambda$ is (induced by) the multiplication of $\Lambda$, and the differential $b_r$ with $r\ge 1$ is given by
\begin{align*}
b_r(a_0[ a_1 \mid  \cdots \mid  a_r ] a_{r+1})
& = a_0 a_1 [ a_2 \mid \cdots\mid  a_r] a_{r+1} \\
& +\sum_{i=1}^{r-1} (-1)^i a_0 [a_1 \mid \cdots\mid a_{i-1} \mid a_i a_{i+1}\mid a_{i+2}\mid \cdots\mid a_r ] a_{r+1} \\
& +(-1)^r a_0[ a_1 \mid \cdots \mid a_{r-1}] a_r a_{r+1},
\end{align*}
where we use the bar notation $a_0[ a_1 \mid \cdots \mid a_r ] a_{r+1} =a_0\ot a_1\ot \cdots\ot a_r\ot a_{r+1}$ for $a_0,\; a_{r+1}\in \Lambda$ and $a_1,\cdots,a_r\in\Lambda_+$.

It is not difficult to verify that the constructions in Section~\ref{Section: preliminaries} carry over verbatim from the bar resolution to the reduced bar resolution; see, for example, \cite{San08}.

\medskip

For an integer $N\ge 2$ and a finite quiver $Q$, let $\Lambda={\K Q}/J^N$ be a truncated quiver algebra. The set
\[
\calB=\bigcup\limits_{n=0}^{N-1}Q_n
\]
is a $\K$-basis of $\Lambda$. We shall make no distinction between an element $p\in \bigoplus\limits_{n=0}^{N-1}\K Q_n$ and its quotient in $\Lambda$. Let us introduce some useful notations. For $n\ge 0$, define
\[
\chi(n)=
\begin{cases}
 kN,           &\text{if } n=2k,\\
 kN+1,         &\text{if } n=2k+1.
\end{cases} 
\]
For composable arrows $\alpha_1, \dots, \alpha_n\in Q_1$, we write
\[
\alpha_{1,\; n}:=\alpha_1\cdots \alpha_n.
\]

Bardzell constructed a minimal projective bimodule resolution $(\bbP_*,d_*)$ for monomial algebras in \cite{Bar97}. More explicitly, for truncated quiver algebra $\Lambda$,
\[
\bbP_*\colon\;
\cdots \lra P_r\xto{d_r} P_{r-1} \lra \cdots \lra P_1 \xto{d_1}
P_0(\xto{\varepsilon} \Lambda),
\]
where $P_r= \Lambda \ot \K Q_{\chi(r)} \ot \Lambda$ and the differential $d_*$ is given by
\[
\begin{aligned}
d_{2k}(1\ot \alpha_{1,\; kN}\ot 1)
    & =\sum\limits_{i=0}^{N-1}\alpha_{1,\; i}\ot \alpha_{i+1,\; (k-1)N+i+1}\ot \alpha_{(k-1)N+i+2,\; kN},\\
d_{2k+1}(1\ot \alpha_{1,\; kN+1}\ot1)
    & =\alpha_1\ot \alpha_{2,\; kN+1}\ot1-1\ot \alpha_{1,\; kN}\ot \alpha_{kN+1},
\end{aligned}
\]
where all $\alpha_i$ are arrows.

\medskip

Ames, Cagliero and Tirao constructed the comparisons morphisms between $\bbP_*$ and $\bbB_*$ \cite{ACT09}. Let us recall their construction.

The chain map
$\mu_*\colon  \bbP_*\lra\bbB_*$ is defined as follows. Let $\mu_0\colon \bbP_0=\Lambda\ot \K Q_0\ot \Lambda \lra \bbB_0=\Lambda \ot \Lambda$ be the identity map.
If $n=2k$ with $k\ge 1$ and $p=1 \ot \alpha_{1, kN}\ot 1$  with all $\alpha_i\in Q_1$, define 
\[
\begin{aligned}
   \mu_{2k}(p)=\sum_{\substack{(x_1,\dots,x_k)\in \mathbb{N}^k\\ 1\le x_i\le N-1}} 1\bigl[ 
   & \underbrace{\alpha_{1,\; x_1}}_{x_1}\mid \underbrace{\alpha_{x_1+1}}_1 \mid \underbrace{\alpha_{x_1+2,\; x_1+x_2+1}}_{x_2} \mid \underbrace{\alpha_{x_1+x_2+2}}_1 \mid \cdots \mid\\
   & \underbrace{\alpha_{k+\sum_{j=1}^{k-1} x_j,\;  k-1+\sum_{j=1}^k x_j}}_{x_k} \mid \underbrace{\alpha_{k+\sum_{j=1}^k x_j}}_1 \bigr]\, \underbrace{\alpha_{k+1+\sum_{j=1}^k x_j,\; kN}}_{k(N-1)-\sum_{j=1}^k x_j}.
\end{aligned}
\]
If $n=2k+1$ with $k\geq 0$ and $p=1\ot \alpha_{1,\; kN+1}\ot 1$ with all $\alpha_i\in Q_1$, define
\[
\begin{aligned}
    \mu_{2k+1}(p) = \sum_{\substack{(x_1,\dots,x_k)\in \mathbb{N}^k\\ 1\le x_i\le N-1}} 1\bigl[ 
    & \underbrace{\alpha_1}_1 \mid \underbrace{\alpha_{2,\;  x_1+1}}_{x_1} \mid \underbrace{\alpha_{x_1+2}}_1 \mid \underbrace{\alpha_{x_1+3,\; x_1+x_2+2}}_{x_2} \mid \cdots \mid \\
    & \underbrace{\alpha_{k+\sum_{j=1}^{k-1} x_j}}_1 \mid \underbrace{\alpha_{k+1+\sum_{j=1}^{k-1} x_j,\; k+\sum_{j=1}^k x_j}}_{x_k} \mid \underbrace{\alpha_{k+1+\sum_{j=1}^k x_j}}_{1} \bigr] \underbrace{\alpha_{k+2+\sum_{j=1}^k x_j,\; kN+1}}_{k(N-1)-\sum_{j=1}^k x_j}.
\end{aligned}
\]

The chain map $\omega_*\colon \bbB_*\lra\bbP_*$ is defined by $\omega_0:=\Id$. For $n\geq 1$ and any 
\[
p=1[p _1\mid \cdots\mid p _n]1\in\Lambda \ot \Lambda_+^n \ot \Lambda,
\]
write $l( p) =\sum\limits_{i=1}^{n}l( p _i) $, and let $p_1\cdots p_n=\alpha_{1, l(p) } \in \K Q_{l(p)}$. When $n=2k$ with $k\geq 0$, define
\[
\omega_{2k}(p)  =  
\begin{cases}
1\ot\underbrace{\alpha_{1,\;kN}}_{kN}\ot\underbrace{\alpha_{kN+1,\; l(p)}}_{l(p) -kN}, &\text{if } p_{2i-1}p_{2i}=0\in \Lambda \text{ for all }  i=1,\dots,k, \\
0,& \text{otherwise}.
\end{cases}
\]
when $n=2k+1$ with $k\geq 0$, define    
\[
\omega_{2k+1}(p)  =  
\begin{cases}
\sum_{j=0}^{l(p_1) -1}\underbrace{\alpha_{1,\;j}}_j\ot\underbrace{\alpha_{j+1,\; kN+j+1}}_{kN+1}\ot\underbrace{\alpha_{kN+j+2,\;l(p) }}_{l(p) -kN-j-1},&\text{if } p_{2i} p_{2i+1}=0\in \Lambda \text{ for all } i=1,\dots,k, \\
0, &\text{otherwise}.
\end{cases}
\]

\begin{Rem} 
The authors of \cite{ACT09} proved the comparison morphisms $\mu_*:\; \bbP_*\lra\bbB_*$ and $\omega_*:\; \bbB_*\lra\bbP_*$ by direct inspection. In 2018 Redondo and Rom\'{a}n constructed comparison morphisms for general monomial algebras and based on these comparison morphisms \cite{RR18}, Redondo and Rossi~Bertone investigated $L_\infty$-structures \cite{RR22}.  In a forthcoming work,  we will give another construction of comparison morphisms for monomial algebras by using algebraic Morse theory and also investigate higher structures.

These comparison morphisms can also be constructed by using the method of contracting homotopies (see, for instance, \cite[Section 2]{IIVZ15}).
Notice that a contracting homotopy on the reduced bar resolution $\bbB_*$ can be found in Section~3.1 in \cite{ACT09}, and a contracting homotopy on $\bbP_*$, was constructed by Sk\"{o}ldberg \cite{Sko08}

\end{Rem}

\medskip

Applying the functor $\Hom_{\Lambda^e}(-,\Lambda)$ to the reduced bar resolution $\bbB_*$,
we obtain the \textit{reduced Hochschild  cochain complex} (see Cibils \cite[Proposition 2.2]{Cib90}), denoted by $\bbB^*$:
\[
\bbB^*\colon  0 \lra\Hom_{\Lambda^e}(\Lambda\ot\Lambda, \Lambda)\xto{b^0}  \cdots \lra \Hom_{\Lambda^e}(\Lambda\ot{\Lambda_+}^{\ot r}\ot\Lambda, \Lambda)
 \xto{b^{r}} \Hom_{\Lambda^e}(\Lambda\ot{\Lambda_+}^{\ot r+1}\ot\Lambda, \Lambda)\lra \cdots ,
\]
and for any $f\in \Hom_{\Lambda^e}(\Lambda\ot{\Lambda_+}^{\ot r}\ot\Lambda, \Lambda)\cong \Hom_{(\K Q_0)^e}( {\Lambda_+}^{\ot r}, \Lambda)$, the differential $b^r$ is given by
\begin{align*}
b^{r}(f)(1[ a_1 \mid \cdots \mid a_{r+1}] 1)
 =& (-1)^{r+1} f(a_1[ a_2 \mid \cdots \mid a_{r+1}] 1) \\
& +\sum_{i=1}^{r}(-1)^{r-i+1} f(1 [a_1 \mid \cdots \mid a_{i-1} \mid a_i a_{i+1} \mid a_{i+2} \mid \cdots \mid  a_{r+1}] 1) \\
&+ f( 1[ a_1 \mid \cdots\mid a_{r} ] a_{r+1}).
\end{align*}

Applying the functor $\Hom_{\Lambda^e}(-,\Lambda)$ to the minimal projective bimodule resolution $\bbP_*$, we  obtain the cochain complex $\bbP^*$:
\[
\bbP^*\colon 0 \lra\Hom_{\Lambda^e}(\Lambda\ot\Lambda, \Lambda)\xto{d^0}  \cdots \lra \Hom_{\Lambda^e}(\Lambda\ot\K Q_{\chi(r)}\ot\Lambda, \Lambda)
\xto{d^{r}} \Hom_{\Lambda^e}(\Lambda\ot \K Q_{\chi(r+1)}\ot\Lambda, \Lambda)\lra \cdots,
\]
with $d^{2k}(f)=-f\circ d_{2k+1}$, and $d^{2k+1}(f)=f\circ d_{2k+2}$ for $k\geq 0$.

The comparison morphisms $\mu_*$ and $\omega_*$, induce the corresponding morphisms $\mu^*$ and $\omega^*$ between $\bbB^*$ and $\bbP^* $, defined by
\[
\mu^r:\Hom_{\Lambda^e}(\Lambda\ot{\Lambda_+}^{\ot r}\ot \Lambda, \Lambda)\lra\Hom_{\Lambda^e}(\Lambda\ot{\K  Q_{\chi(r)}}\ot\Lambda, \Lambda),\ f\mapsto f\circ \mu_r,
\]
and
\[
\omega^r:\Hom_{\Lambda^e}(\Lambda\ot{\K  Q_{\chi(r)}}\ot\Lambda, \Lambda)\lra\Hom_{\Lambda^e}(\Lambda\ot{\Lambda_+}^{\ot r}\ot\Lambda, \Lambda),\ g\mapsto g\circ \omega_r.
\]
In the following sections, we shall use these morphisms to compute the cup product, the Lie bracket and the BV structure on the complex $\bbP^*$.

\bigskip

\section{The Hochschild cohomology spaces}\label{Section:Hochschild cohomology groups}\

Let $\Lambda={\K Q}/J^N$ be a truncated quiver algebra with $N\geq 2$. Truncated basic cycle algebra is a special case of truncated quiver algebra. In this section, we give a new description of the complex $\bbP^*$ and provide an explicit $\K$-bases of equivalence classes for Hochschild cohomology groups of truncated  basic cycle algebras, following Locateli \cite{Loc99}, and Bardzell, Locateli and Marcos \cite{BLM00}. To keep the notation simple, a representation of each equivalence class will be used to denote the class itself.

\subsection{The Hochschild cohomology spaces for truncated quiver algebras}\ 

The notion of \textit{parallel paths}, introduced by Cibils \cite{Cib90}, provide a convenient description of the spaces $\bbP^*=\Hom_{\Lambda^e}(P_*,\Lambda)$. Recall that two paths $p$ and $q$ are called \textit{parallel} if they have the same origin and target; we write $p\pl  q$. For subsets $X$ and $Y$ of paths in $Q$, set 
\[
(X\pl  Y):=\{(p, q)\in X\times Y\di \ p \pl  q\},
\]
and denote by $\K(X\pl  Y)$ the $\K$-vector space with basis $(X\pl  Y)$.

\begin{Lem}[{\cite{Cib90}, \cite[Lemma 1]{Loc99}}]\label{Lem: Parallel paths}

For every $r\ge 0$, there is an isomorphism of $\K$-vector spaces
\[
\bbP^r=\Hom_{\Lambda^e}(P_r,\Lambda)\cong \K( Q_{\chi(r)}\pl \calB),
\]
where for each $(\gamma, b)\in (Q_{\chi(r)}\pl \mathcal B)$, the map 
\[
f_{(\gamma, b)}\colon P_r=\Lambda\ot \K Q_{\chi(r)}\ot\Lambda\to \Lambda
\]
is given on basis elements by sending $p\ot \gamma'\ot q\in P_r$ to $\delta_{\gamma, \gamma'} pbq$, where $p,\; q\in \calB$ and $\gamma'\in Q_{\chi(r)}$.

\end{Lem}

Locateli described the differential of the  complex $\bbP^*$ of a truncated quiver algebra $\Lambda$ in terms of parallel paths.
 
\begin{Prop}[{\cite[Section 3]{Loc99}}]\label{Prop:new description of P}

The differential in the  complex
\[
\bbP^*: \cdots \to \bigoplus_{i=0}^{N-1} \K(Q_{kN}\pl  Q_i) \xto{d^{2k}} 
\bigoplus_{i=0}^{N-1} \K(Q_{kN+1}\pl  Q_i)\xto{d^{2k+1}} 
\bigoplus_{i=0}^{N-1} \K(Q_{(k+1)N}\pl  Q_i)\to \cdots 
\]
has the form (for $k\geq 0$)
\[
d^{2k}=
\begin{bmatrix}
0 & 0 & 0 & \cdots & 0 & 0 \\
d^{2k}_0 & 0 & 0 & \cdots & 0 & 0 \\
0 & d^{2k}_1 & 0 & \cdots & 0 & 0 \\
0 & 0 & d^{2k}_2 & \cdots & 0 & 0 \\
\vdots & \vdots & \vdots & \ddots & \vdots & \vdots \\
0 & 0 & 0 & \cdots & d^{2k}_{N-2} & 0
\end{bmatrix},\qquad
d^{2k+1}=
\begin{bmatrix}
0 & 0 & \cdots &  0 \\
\vdots & \vdots & \ddots & \vdots \\
0 & 0 & \cdots &  0 \\
d^{2k+1}_{0} & 0 & \cdots & 0
\end{bmatrix}.
\]
Here, for $i=0,\dots,N-2$, 
\[
d^{2k}_i\colon \K(Q_{kN} \pl  Q_i) \to \K(Q_{kN+1} \pl  Q_{i+1})
\]
is given by
\[
d^{2k}_i(p,q)=-\sum_{\alpha\in Q_1}(\alpha p,\alpha q) +\sum_{\beta\in Q_1}(p \beta, q \beta)
\]
for  $(p,q)\in (Q_{kN}\pl  Q_i)$, and 
\[
d^{2k+1}_0\colon \K(Q_{kN+1} \pl  Q_0) \to \K(Q_{(k+1)N} \pl  Q_{N-1})
\]
is given by
\[
d^{2k+1}_0(p,q)=\sum_{uv\in Q_{N-1}}(u p v,  u v).
\]
for $(p,q)\in (Q_{kN+1} \pl Q_0)$.
\end{Prop}

\medskip

\subsection{The case $e=1$ for truncated basic cycle algebras $\Lambda=\K Z_e/J^N$}\ 

In this subsection, we briefly discuss here the  well known case $e=1$  for truncated basic cycle algebras $\Lambda=\K Z_e/J^N$   without proof for completeness.
In the main text of this paper,  we will consider only  the case  $e\geq 2$,

The   truncated basic cycle algebra   $\Lambda=\K Z_1/J^N$ with $N\ge 2$  is  given by the bounded quiver $Z_1$ with exactly one vertex $1$ and one loop $\alpha$, and we have the relation $\alpha^N=0$.

The Hochschild cohomology groups of $\Lambda$ are given as follows:
\begin{itemize}
\item[(1)] 
\[
\rmHH^0(\Lambda)=\K \{ (e_1,\alpha^i)\di 0\le i\le N-1 \},
\]

  \item[(2)]  for $k\ge 1$,   \begin{equation*}
	 \rmHH^{2k}(\Lambda)=
	 \begin{cases}
	 \K  \left\{(\alpha^{kN},\alpha^i) \di 0\le i\le N-1 \right\} , & \text{if }    \rmchar(\K) \di N,\\
	 \K\left\{(\alpha^{kN},\alpha^i) \di 0\le i\le N-2 \right\} , &\text{otherwise},
	 \end{cases}
      \end{equation*}
  
 \item[(3)]  for $k\ge 1$,   
\begin{equation*}
	 \rmHH^{2k-1}(\Lambda)=
	 \begin{cases}
	 \K  \left\{(\alpha^{(k-1)N+1},\alpha^i) \di 0\le i\le N-1 \right\} , & \text{if }    \rmchar(\K) \di N,\\
	 \K\left\{(\alpha^{(k-1)N+1},\alpha^i) \di 1\le i\le N-1 \right\} , &\text{otherwise}.
	 \end{cases}
\end{equation*}
\end{itemize}
Then $\dim (\rmHH^0(\Lambda))=N$ and for $n\ge 1$,
\begin{equation*}
    \dim(\rmHH^n(\Lambda))=
    \begin{cases}
        N, & \text{if } \rmchar (\K)\di N, \\
        N-1, & \text{otherwise}.
    \end{cases}
\end{equation*}

If $\rmchar(\K)\nmid N$, the Hochschild cohomology ring $\rmHH^*(\Lambda)$ is isomorphic to a graded commutative algebra generated by $x_0,y,z$ with degree $0,1,2$ respectively, subject  to the relations $x_0^N=0$, $y x_0^{N-1}=0$, $x_0^{N-1}z=0$ and $y^2=0$. Moreover, $\rmHH^*(\Lambda)$ is a Gerstenhaber algebra with all nonzero Gerstenhaber bracket between the basis elements are 
 \[
 [x_0,y]=-x_0,\; \text{and}\; [y,z]=-N z.
 \]
The BV differential $\Delta$ is given by
\[
\Delta(yx_0^az^b)=(bN+N-a-1)x_0^az^b,\quad \text{for }\; 0\le a\le N-2,\; b\ge 0.
\]

If $\rmchar (\K)\di N$, then the Hochschild cohomology ring $\rmHH^*(\Lambda)$ is isomorphic to a graded commutative algebra generated by  generated by $ x_0,y, w, z'$ with degree $ 0,1,1,2$ respectively, subject  to the relations $x_0^N=0$, $y= w x_0$ and $w^2=-\frac{N(N-1)}{2}x_0^{N-2}z'$. All nonzero Gerstenhaber bracket between the basis elements are 
 \[
 [x_0,y]=-x_0,\; [x_0,w]=-1,\; \text{and}\; [y,w]=-w.
 \]
The BV differential $\Delta$ is given by
\[
\Delta(w z'^b)=0\ \text{and }	\Delta(w x_0^az'^b)=-a \, x_0^{a-1}z'^b,\quad \text{for }\; 1\le a\le N-1,\; b\ge 0.
\]

\medskip

\subsection{The case of $e\geq  2$ for  truncated basic cycle algebras}\ 

\textit{From now on, when  considering  the truncated basic cycle algebra $\Lambda={\K Z_e}/J^N$ with $N\ge 2$, we always assume that  $e\ge 2$. 
}

For a vertex $l$, we denote by $\alpha_l^j$ the path of length $j$ starting at $l$ with the convention that when $j=0$, $\alpha_l^0=e_l$. For any integer $x\in \Z$, we denote by $\ul x$ (resp. $\ol{x}$) the unique representative of $x$ modulo $e$ satisfying $1\le \ul x \le e$ (resp. $0\le \ol x\le e-1$). Then $Z_e$ is a basic cycle with $e$ vertices $\{1, \dots,e\}$ and $e$ arrows $\{\alpha_1, \dots, \alpha_e\}$ such that $o(\alpha_i)=i$ and $t(\alpha_i)=\ul{i+1}$.

Notice that for $k\ge 0$ and $0\le i \le N-1$, the set $(Q_{kN}\pl  Q_i)$ is nonempty only when $i\equiv kN\mode$, then
\[
(Q_{kN}\pl  Q_i)=\{(\alpha_l^{kN},\alpha_l^i)\di l=1, \dots, e\}.
\]
Similarly, $(Q_{kN+1}\pl  Q_i)$ is nonempty only when $i\equiv kN+1\mode$, and then 
\[
(Q_{kN+1}\pl  Q_i)=\{ (\alpha_l^{kN+1}, \alpha_l^i)\di l=1, \dots, e\}.
\]

\begin{Cor}[{\cite[Section 2]{BLM00}}]\label{Cor:new description of P for basic cycle}
For the truncated basic cycle algebra $\Lambda$, the differentials $d^{2 k}$ and $d^{2 k+1}$ in the following complex
\[
 \bbP^*: \cdots \lra \bigoplus_{i=0}^{N-1} \K(Q_{kN}\pl  Q_i) \xto {d^{2k}} 
 \bigoplus_{i=0}^{N-1} \K(Q_{kN+1}\pl  Q_i) \xto{d^{2k+1}} 
 \bigoplus_{i=0}^{N-1} \K(Q_{(k+1)N}\pl  Q_i)\lra \cdots
\]
have the same form as in Proposition~\ref{Prop:new description of P}.
For $i=0,\dots,N-2$ with $kN\equiv i\mode$, the map
\[
d^{2k}_i\colon \K(Q_{kN}\pl  Q_i)  \to  \K(Q_{kN+1}\pl  Q_{i+1})
\]
is given on basis elements by
\[
d^{2k}_i(\alpha_l^{kN},\alpha_l^i)=-(\alpha_{\ul{l-1}}^{kN+1},\alpha_{\ul{l-1}}^{i+1})+(\alpha_l^{kN+1},\alpha_l^{i+1}),\qquad 1\leq l\leq e.
\]
If $kN+1\equiv 0\mode$, then
\[
d^{2k+1}_0\colon \K(Q_{kN+1}\pl  Q_0)  \to  \K(Q_{(k+1)N}\pl  Q_{N-1})
\]
is given by
\[
d^{2k+1}_0(\alpha_l^{kN+1},e_l)=\sum_{t=0}^{N-1}(\alpha_{\ul{l-t}}^{(k+1)N},  \alpha_{\ul{l-t}}^{N-1}),\qquad 1\le l\le e.
\]
\end{Cor}

Bardzell, Locateli and Marcos \cite[Section 5]{BLM00} determined explicit basis for Hochschild cohomology groups of a truncated basic cycle algebra.

To simplify the presentation, we introduce the following notations. 

For any $k\ge 0$ and $0\le i,j\le N-1$ with $i\equiv kN\mode$ and $j\equiv kN+1\mode$, set
\[
x_{kN, i}:= \sum_{l=1}^e( \alpha_l^{kN} ,\alpha_l^i),  \qquad y_{kN+1, j}:= \sum_{l=1}^e(\alpha_l^{kN+1} ,\alpha_l^j).
\]
For any $k\ge 0$ and $1\le i\le N-1$ with $i\equiv kN\mode$, we write
\[
v_{1;\; kN+1,i}:=( \alpha_1^{  kN+1} ,\alpha_1^{i}).
\]
When $N\equiv 1\mode$, for $1\le l\le N-1$, we put
\[
u_{l;\; 0, N-1}:=( e_l ,\alpha_l^{N-1}).
\]

\begin{Thm}[{Compare with \cite[Section 5]{BLM00}}]\label{Thm: basis} 
Let $\Lambda={\K Z_e}/J^N$ with $N \ge 2$ and $e\ge 2$. Then
\begin{itemize}

\item [(1)] $\rmHH^0(\Lambda)=
      \begin{cases}
      \K \left[\left\{x_{0, i}\di  0\le i\le N-2,\; i \equiv 0 \mode\right\}\cup 
      \left\{u_{l;\;  0,N-1}\di  1\le l\le e\right\} \right],  &\text{if }  N\equiv 1 \mode,\\
       \K  \left\{x_{0, i}\di 0\le i\le N-2,\; i \equiv 0 \mode\right\},  &\text{if } N\not\equiv 1 \mode;
      \end{cases}$

\item [(2)] for $k\ge 1$,
\[
   \begin{aligned}
	 \rmHH^{2k}(\Lambda)
     &=
	 \begin{cases}
	 \K  \left\{x_{kN, i}\di 0\le i\le N-1,\; i \equiv N-1 \mode\right\} , & \text{if } \rmchar(\K) \di N  \text{ and }  k N \equiv N-1 \mode,\\
	 \K\left\{x_{kN, i}\di 0\le i\le N-2,\; i \equiv kN \mode\right\} , &\text{otherwise};
	 \end{cases}\\
     &= \begin{cases}
	 \K  \left\{x_{kN, i}\di 0\le i\le N-1,\; i \equiv kN \mode\right\}, & \text{if } \rmchar(\K) \di N \text{ and } \gcd(N,e)=1,\\
	 \K\left\{x_{kN, i}\di 0\le i\le N-2,\; i \equiv kN \mode\right\} , &\text{otherwise};
	 \end{cases}
     \end{aligned}
 \]

\item [(3)] for $k\ge 0$,
      \[
      \begin{aligned}
	 \rmHH^{2 k+1}(\Lambda) &=
	 \begin{cases}
	 \K \left\{y_{kN+1, i }\di 0\le i\le N-1,\; i \equiv 0  \mode\right\}, &\text{if }  \rmchar (\K) \di N  \text{ and } kN  \equiv -1 \mode,\\
      \K \left\{v_{1;\; kN+1, i}\di 1\le i\le N-1,\; i\equiv kN+1 \mode\right\}, & \text{if } \rmchar (\K) \di e, \\
	 \K\left\{y_{kN+1, i}\di 1\le i\le N-1, i \equiv kN+1  \mode\right\}, &\text{otherwise};
	 \end{cases}\\
     &= \begin{cases}
	 \K \left\{y_{kN+1, i }\di 0\le i\le N-1,\; i \equiv kN+1  \mode\right\}, &\text{if } \rmchar (\K) \di N \text{ and } \gcd(N,e)=1,\\
      \K \left\{v_{1;\; kN+1, i}\di 1\le i\le N-1,\; i\equiv kN+1 \mode\right\}, & \text{if } \rmchar (\K) \di e, \\
	 \K\left\{y_{kN+1, i}\di 1\le i\le N-1, i \equiv kN+1  \mode\right\}, &\text{otherwise}.
	 \end{cases}
     \end{aligned}
    \]

\end{itemize}

\end{Thm}

\begin{SProof}

The case $\rmHH^0(\Lambda)$ is easy. For the rest of the proof, let $k\ge 0$.

For $0\le i\le N-2$ and $i\equiv kN \mode$, under the ordered  basis $\left\{( \alpha_l^{kN} ,\alpha_l^{i}), 1\le l\le e\right\}$, the matrix of
\[
d_i^{2k}\colon \K(Q_{kN}\pl  Q_i) \lra \K(Q_{kN+1}\pl  Q_{i+1})
\]
has the form:
\[
d_i^{2k}=
\begin{bmatrix}
 		1 & 0 & 0 & \cdots & 0 & -1 \\
 		-1 & 1 & 0 & \cdots & 0 & 0 \\
 		0 & -1 & 1 & \cdots & 0 & 0 \\
 		\vdots & \vdots & \vdots & \ddots & \vdots & \vdots \\
 		0 & 0 & 0 & \cdots & 1 & 0 \\
 		0 & 0 & 0 & \cdots & -1 & 1
\end{bmatrix}_{e}\xto[\text{}]{\text{Elementary row operations}}
\begin{bmatrix}
 	1 & 0 & 0 & \cdots & 0 & -1 \\
 	0 & 1 & 0 & \cdots & 0 & -1 \\
 	0 & 0 & 1 & \cdots & 0 & -1 \\
 	\vdots & \vdots & \vdots & \ddots & \vdots & \vdots \\
 	0 & 0 & 0 & \cdots & 1 & -1 \\
 	0 & 0 & 0 & \cdots & 0 & 0
\end{bmatrix}_{e}.
\]
So the rank of $d_i^{2k}$ is $e-1$. We see that for $0\leq i\leq N-2$,
\[
\Ker (d_i^{2k})=
\begin{cases} 
\K\{\sum\limits_{l=1}^e( \alpha_l^{kN} ,\alpha_l^i)\}=\K\{x_{kN, i}\},&\text{if }  i\equiv kN \mode,\\
0, & \text{otherwise};
\end{cases}
\]
and 
\[
\rmIm (d_i^{2k})=
\begin{cases} 
\K\{-(\alpha_{ l-1}^{kN+1},\alpha_{l-1}^{i+1})+(\alpha_l^{kN+1},\alpha_l^{i+1})\di 2\leq l\leq e \},& \text{if } i\equiv kN \mode,\\ 
0, & \text{otherwise}.
\end{cases}
\]

Let $N=me+t$ with $m\ge 0$, $0\le t \le e-1$. When $(k-1)N+1\equiv 0\mode$, the matrix of
\[
d_0^{2k-1}\colon \K(Q_{(k-1)N+1}\pl  Q_0)  \lra  \K(Q_{kN}\pl  Q_{N-1})
\]
has  the form:
\[
\left[
\begin{array}{ccccccccccc}
		m+1 & m & m & m & \cdots & m & m & m+1 & m+1 & \cdots & m+1\\
		m+1 & m+1 & m & m & \cdots & m & m & m & m+1 &\cdots & m+1\\
		m+1 & m+1 & m+1 & m & \cdots & m & m & m & m & \cdots & m+1\\
		\vdots & \vdots & \vdots & \vdots & \ddots & \vdots & \vdots &\vdots &\vdots&\ddots&\vdots\\
		m & m & m & m & \cdots & m & m+1&m+1 & m+1 &\cdots &m+1
\end{array}
\right]_{e}
\]
where in each row we have $t$ entries $m+1$ and $e-t$ entries $m$, which, by elementary matrix operations, can be transformed into 
\[
\begin{bmatrix}
	1 & 0 & 0 & \cdots & 0 & -1 \\
	0 & 1 & 0 & \cdots & 0 & -1 \\
	0 & 0 & 1 & \cdots & 0 & -1 \\
	\vdots & \vdots & \vdots & \ddots & \vdots \\
	0 & 0 & 0 & \cdots & 1 & -1 \\
	0 & 0 & 0 & \cdots & 0 & N
\end{bmatrix}_{e}
\]
Then  
\[
\mathrm{rank}(d_0^{2k-1})=
\begin{cases} 
e, & \text{if }\,  (k-1)N+1\equiv 0\mode \text{ and } \rmchar (\K) \nmid   N, \\
e-1, & \text{if }\,  (k-1)N+1\equiv 0\mode \text{ and } \rmchar (\K)   \di  N, \\
0, &  \text{if }\,  (k-1)N+1\not\equiv 0\mode.
\end{cases}
\]
When $(k-1)N+1\equiv 0\mode$ and $ \rmchar (\K)  \di  N$, 
\[
y_{(k-1)N+1, 0}=\sum_{l=1}^e( \alpha_l^{(k-1)N+1} ,e_l) \in \Ker(d_0^{2k-1})
\]
and generates $\Ker(d_0^{2k-1})$, while the image of $d_0^{2k-1}$ is spaned by 
$\left\{\sum\limits_{t=0}^{N-1}(\alpha_{\ul{l-t}}^{kN},  \alpha_{\ul{l-t}}^{N-1}), 1\leq l\leq e-1\right\}$; in all other cases, $d_0^{2k-1}$ is an isomorphism.

Notice that the condition that 
$(k-1)N+1\equiv 0\mode \text{ and } \rmchar (\K) \di  N$ 
is contradictory to $\rmchar (\K) \di e$, so when $\rmchar (\K) \di e$, $d_0^{2 k-1}$ is an isomorphism.

For $k\geq 1$, 
\[
\rmHH^{2k}(\Lambda)=\left(\bigoplus_{i=0}^{N-2}\Ker(d_i^{2k})\right)\bigoplus \left(\K(Q_{kN}\pl  Q_{N-1})/\rmIm(d^{2k-1}_0)\right).
\]
It is not difficult to see  that 
\begin{equation*}
	\mathrm{HH}^{2k}(\Lambda)=
	\begin{cases}
		 \K  \left\{x_{kN, i}\di 0\le i\le N-1,\; i \equiv N-1 \mode\right\} , & \text{if } \rmchar(\K) \di N  \text{ and }  k N  \equiv N-1 \mode,\\
		 \K\left\{x_{kN, i}\di 0\le i\le N-2,\; i \equiv kN \mode\right\} , &\text{otherwise};
	\end{cases}
\end{equation*}

For $k\ge 0$,  
\[
\rmHH^{2k+1}(\Lambda)=\Ker(d^{2k+1}_0) \bigoplus 
\left(\bigoplus_{i=1}^{N-1}\K(Q_{kN+1} \pl  Q_i)/{\rmIm (d_{i-1}^{2k})}\right).
\]
The $\K$-vector space $\K(Q_{kN+1}\pl  Q_{i})/{\rmIm d_{i-1}^{2k}}$ is one-dimensional for $1\le i\le N-1$ with $i\equiv k N+1 \mode$. When $\rmchar(\K) \nmid e$, its basis can be chosen as
$y_{kN+1, i}=\sum\limits_{l=1}^e( \alpha_l^{kN+1} ,\alpha_l^i)$. In fact, the class of $\sum\limits_{l=1}^e( \alpha_l^{kN+1} ,\alpha_l^i)$ in $\rmHH^{2k+1}(\Lambda)$ is equal to that of $e( \alpha_1^{kN+1} ,\alpha_1^i)$ which does not vanish as $\rmchar (\K) \nmid e$. However, if $\rmchar (\K)\di  e$, one can only choose one representative among $( \alpha_l^{kN+1} ,\alpha_l^{i})$ with  $1\le l\le e$, say, for instance, $v_{1;\; kN+1,i}=( \alpha_1^{kN+1} ,\alpha_1^i)$.

We have shown that
\begin{equation*}
	\rmHH^{2 k+1}(\Lambda)=
	\begin{cases}
	\K \left\{y_{kN+1, i }\di 0\le i\le N-1,\; i \equiv 0  \mode\right\}, &\text{if }  \rmchar (\K) \di N  \text{ and } kN  \equiv -1 \mode,\\
    \K \left\{v_{1;\; kN+1, i}\di 1\le i\le N-1,\; i\equiv kN+1 \mode\right\}, & \text{if } \rmchar (\K) \di e, \\
	\K\left\{y_{kN+1, i}\di 1\le i\le N-1, i \equiv kN+1  \mode\right\}, &\text{otherwise}.
	\end{cases}
\end{equation*}

\end{SProof}

\begin{Rem} 

As is noticed in the proof of the above result, when $\rmchar(\K) \di e$, the  choice of the representative $\sum\limits_{l=1}^e( \alpha_l^{kN+1} ,\alpha_l^i)$ of $\rmHH^{2k+1}(\Lambda)$ in \cite[Section 5]{BLM00} is not correct because in this case,   
\[
\sum\limits_{l=1}^e( \alpha_l^{kN+1} ,\alpha_l^i)= e( \alpha_1^{kN+1} ,\alpha_l^i)=0\in   \rmHH^{2k+1}(\Lambda).
\]
However, our choice of the representative $v_{1;\; kN+1,i}= (\alpha_1^{kN+1} ,\alpha_1^i)$ of $\rmHH^{2k+1}(\Lambda)$ is NOT canonical, one can take $(\alpha_j^{kN+1} ,\alpha_j^i)$ for any $1\leq j\leq e$.

\end{Rem}

\begin{Cor} \label{Cor: dimension formulae}

Let $N=me+t$ with $m\ge 0$, $0\le t \le e-1$. For $k\geq 1$, define 
\[
\calD_k=\calD_k(N, e)= \left\{x_{kN,i} \di 0\le i\le N-2,\; i \equiv kN \mode\right\} 
\]
and write the number of the set $D_k(N, e)$ as $\# (\calD_k(N, e))=[\frac{N-2-\ol{kN}}{e}]+1$.
Then the dimension of the Hochschild cohomology groups of $\Lambda=\K Z_e/J^N$ are given as follows:
\[
\dim(\rmHH^0(\Lambda))=
\begin{cases} 
m, &\text{if } t=0,\\
m+e, & \text{if } t=1,\\
m+1, & \text{if } 2\leq t\leq e-1;
\end{cases}
\]
for $n\ge 1$,		
\begin{equation*}
		\dim(\mathrm{HH}^{n}(\Lambda))=
		\begin{cases}
			\# (D_{[\frac{n}{2}]}(N, e))+1, & \text{if } \rmchar (\K) \di N  \text{ and }{[\frac{n+1}{2}]}N \equiv N-1 \mode,\\
			\# (D_{[\frac{n}{2}]}(N, e)), &\text{otherwise},
		\end{cases}
\end{equation*} 		
where $[x]$ denotes the maximal integer not greater than $x$.

\end{Cor}

\begin{SProof} 

The dimension formulae for $\rmHH^0(\Lambda)$ is easy to obtain. For other cases, notice that the map sending $x_{kN,i}=\sum\limits_{l=1}^e( \alpha_l^{kN} ,\alpha_l^i)$ to $y_{kN+1,i+1}=\sum\limits_{l=1}^e( \alpha_l^{kN+1} ,\alpha_l^{i+1})$ establishes a bijection between
$\calD_k$ and $\calE_{k} = \calE_k(N, e)=\left\{y_{kN+1, i+1}\di 0\le i\le N-2, i \equiv kN  \mode\right\} $ for any $k\geq 0$.

\end{SProof}

\begin{Rem} 

Note that our dimension formulae  in Corollary~\ref{Cor: dimension formulae} differ  from those in \cite[Theorem 2]{Loc99} and \cite[Section 2]{BLM00}. The reason is that, as shown in Corollary~\ref{Cor: dimension formulae}, the equality 
\[
\dim (\rmHH^{2k}(\Lambda))=\dim(\rmHH^{2k+1}(\Lambda))
\]
does not always hold for $k\ge 1$. Moreover, as observed in \cite[Section 2]{BLM00}, the dimension $\dim(\rmHH^1(\Lambda))$ can also be given by 
\[
\dim(\rmHH^1(\Lambda))=
\begin{cases} 
m, &\text{if } t=0,1,\\
m+1, & \text{if } 2\leq t\leq e-1.
\end{cases}
\]

Our dimension formulae coincide  with those contained in    \cite[4.4 Proposition and 5.3 Proposition]{EH99}.  

The dimension formulae for $\rmHH^{\geq 1}(\Lambda)$ in \cite{Zha97} is only valid in characteristic zero. 
\end{Rem}

\bigskip

\section{Cup product}\label{Section: cup product}\

In this section, we present the cup product formulas between selected representatives in the Hochschild cohomology groups using comparison morphisms, whereas Bardzell, Locateli and Marcos \cite{BLM00} used the Yoneda product.

\medskip

The definition of cup product that defined in Section \ref{Section: preliminaries} carries over for the reduced bar resolution  of $\Lambda$. More precisely, given two cochains $f\in \bbB^m$ and $ g\in \bbB^n$, the cup product of $f$ and $g$ is the cochain 
\[
f\cup_{\bbB} g\in \bbB^{m+n}=\Hom_{\Lambda^e}(\Lambda \ot \Lambda_+^{\ot m+n} \ot \Lambda, \Lambda)
\]
defined by
\[
f \cup_{\bbB} g(1[ a_1 \mid \cdots\mid a_{m+n} ]1) =(-1)^{m n}f(1[ a_1\mid \cdots\mid a_{m} ]1 )g(1[ a_{m+1} \mid \cdots\mid a_{m+n}]1 ).
\]

Using the comparison morphisms $\mu^*$ and $\omega^*$ defined in Section~\ref{section: comparison morphism}, Ames, Cagliero, and Tirao \cite{ACT09} described the cup product on the complex $\bbP^*$ via the following method: For $f\in \bbP^m$ and  $g\in \bbP^n$,
\[
f \cup_{\bbP} g =\mu^{m+n}(\omega^m f \cup_{\bbB} \omega^n g )=(f\circ \omega_m\cup_{\bbB} g \circ  \omega_n)\circ \mu_{m+n}.
\]

From now on, for simplicity, we shall write $\vee:=\cup_\bbP$ and $a_{1,n}=a_1\cdots a_n$.

\begin{Prop}[{\cite[Proposition 7.5]{ACT09}}]\label{Prop:cup product for truncated quiver algebras}

Let $f\in \bbP^m$ and $g\in \bbP^n$.  Let $k,\, h\ge 0$, we have

\begin{enumerate}

\item   if $m=2k$ and $n=2h$, then
        \[
        f\vee g(1\ot a_{1, kN} b_{1, hN}\ot 1)=f(1 \ot a_{1, kN}\ot 1)g (1 \ot b_{1, hN}\ot 1),
        \]
			
		\item if $m=2k$ and $n=2h+1$, then
        \[
        f\vee g(1 \ot a_{1, kN} b_{1,  hN+1}\ot 1)=f(1 \ot a_{1, kN}\ot 1)g(1 \ot b_{1, hN+1}\ot 1),
        \]
			
		\item if $m=2k+1$ and $n=2h$, then
        \[
        f\vee g(1 \ot a_{1, kN+1}b_{1, hN}\ot 1)=f(1 \ot a_{1, kN+1}\ot 1)g(1 \ot b_{1, hN}\ot 1),
        \]
			
		\item if $m=2k+1$ and $n=2h+1$, then
        \[
        \begin{aligned}
            & f\vee g(1\ot a_{1, (k+h+1)N}\ot 1)  \\
            &=-\sum\limits_{1\le i<j\le N} a_{1, i-1} f(1 \ot a_{i, i+kN}\ot 1)a_{i+kN+1, j+kN-1} g(1 \ot a_{j+kN, j+(k+h)N}\ot 1) a_{j+(k+h)N+1, (k+h+1)N}.  
        \end{aligned}
        \]
            
\end{enumerate}
        
\end{Prop}

\begin{Rem}

For a truncated quiver algebra, the cochain complex $\bbP^*$ is NOT a differential graded algebra, but an  $A_\infty$-algebra. In fact, $\bbP^*$ is homotopy equivalent to $\bbB^*$, which itself is a differential graded algebra, so by homotopy transfer technique, $\bbP^*$ is an $A_\infty$-algebra. However, when $f, g, h \in \bbP^*$ have odd degrees, one sees easily that $(f\vee g)\vee h\neq f\vee (g\vee h)$, so $\bbP^*$ is a genuine $A_\infty$-algebra. It would be very interesting to determine this higher structure for $\bbP^*$. This will be done in a forthcoming paper using homotopy transfer technique.
		
\end{Rem}

By direct inspection, we obtain the following consequence of Proposition~\ref{Prop:cup product for truncated quiver algebras}.

\begin{Cor}\label{Cor: cup product for truncated basic cycles}

Let $\Lambda={\K Z_e}/J^N$ be a truncated basic cycle algebra with $N\ge 2$ and $e\ge 2$. Let $m, n\ge 0$,  $1\le l, p\le e$ and $0\le r, s\le N-1$. Assume that $\chi(m)\equiv r \mode$  and $\chi(n)\equiv s \mode$. Then the cup product of two elements $(\alpha_l^{\chi(m)},\alpha_l^r)\in(Q_{\chi(m)}\pl  \calB)$ and  $(\alpha_p^{\chi(n)},\alpha_p^s)\in(Q_{\chi(n)}\pl  \calB)$ are given as follows:
\begin{enumerate}

	\item   If $m$ or $n$ is even, then 
    \[
    (\alpha_l^{\chi(m)},\alpha_l^r)\vee (\alpha_p^{\chi(n)},\alpha_p^s)=
    \begin{cases}
    \delta_{p,\ul{l+r}}(\alpha_l^{\chi( m+n)},\alpha_l^{r+s})\in (Q_{\chi(m+n)}\pl  \calB), & \text{if } r+s\le N-1,\\
    0, & \text{otherwise}.
    \end{cases}
    \]

	\item   If $n$ and $m$ are odd, then
	\begin{equation*}
	\begin{aligned}
	(\alpha_l^{\chi(m)},\alpha_l^r)\vee (\alpha_p^{\chi(n)},\alpha_p^s)
    &=-\sum\limits_{1\le i<j\le N}\delta_{p,\ul{l+r+j-i-1}}(\alpha_{\ul{l-i+1}}^{\chi(m+n)}, \alpha_{\ul{l-i+1}}^{N-2+r+s})\\
	&=-\sum\limits_{i=0}^{N-2}a_i(\alpha_{\ul{l-i}}^{\chi(m+n)},\alpha_{\ul{l-i}}^{N-2+r+s})\in (Q_{\chi(m+n)}\pl  \calB),
	\end{aligned}
	\end{equation*}
	where for each $i\in [0, N-2]$, $a_i$ is the number of $j$ such that $j\equiv p-l-r+i+1 \mode$ and $i+1\le j\le N-1$. In particular, whenever $r+s>1$,  $(\alpha_l^{\chi(m)},\alpha_l^r)\vee (\alpha_p^{\chi(n)},\alpha_p^s)=0$ .
            
\end{enumerate}
		
\end{Cor}

From Corollary~\ref{Cor: cup product for truncated basic cycles}, we obtain the cup product of two chosen basis elements in $\rmHH^*(\Lambda)$.

\begin{Prop}\label{Prop: cup product of basis elements}

Let $\Lambda={\K Z_e}/J^N$ be a truncated quiver algebra with $N \ge 2$. Let $m,n\ge 0$. Then the cup product of two basis elements in $\rmHH^{m}(\Lambda)$ and $\rmHH^n(\Lambda)$ are given as follows:

\begin{enumerate}
			
\item ($\rmHH^0\vee \rmHH^0$) Assume that $m=n=0$.
				
	\begin{itemize}
    
	\item[(i)]  Suppose $N\not \equiv 1 \mode$. Then for any $0\le i, j\le N-2$ with $i, j\equiv 0\mode$,
    \begin{equation*}
	x_{0,i} \vee x_{0,j}=
	\begin{cases}
	x_{0,i+j}, & \text{if }  i+j\le N-2 ,\\
	0, &\text{otherwise}.
	\end{cases}
	\end{equation*}

	\item[(ii)]  Suppose $N\equiv 1 \mode$. Then
    
		\begin{itemize}
        
		\item[(a)]  for any $1\le l, p\le e$,  $u_{l; \; 0,N-1} \vee u_{p;\; 0,N-1}=0$,
        
		\item[(b)]  for any $1\le l\le e$ and $0\le i\le N-2$ with $i\equiv 0\mode$,  \begin{equation*}
		u_{l; \;  0,N-1} \vee x_{0,i}=
		\begin{cases}
		u_{l; \;  0,N-1}, & \text{if }  i=0 ,\\
		0, &\text{otherwise},
		\end{cases}
		\end{equation*} 	
        
		\item[(c)]  for any $0\le i, j\le N-2$ with $i, j\equiv 0\mode$,
        \begin{equation*}
		x_{0,i} \vee x_{0,j}=
		\begin{cases}
		x_{0,i+j}, & \text{if } i+j\le N-2 ,\\
		\sum\limits_{l=1}^{e} u_{l;\; 0,N-1}, & \text{if } i+j=N-1,\\
		0, & \text{otherwise}.
		\end{cases}
		\end{equation*} 
        
	  \end{itemize} 	
				
	\end{itemize}

\item ($\rmHH^{0}\vee \rmHH^{2k}$) Assume that $m=0$ and $n=2k$ with $k\ge 1$.
			
	\begin{itemize}
    
	\item[(i)] Suppose that  $N\not \equiv 1 \mode$ and that  $\rmchar (\K)\nmid N$ or $kN\not \equiv N-1\mode$. Then for any $0\le i,j\le N-2$ with $i\equiv 0\mode$ and $j\equiv kN\mode$,
	\begin{equation*}
	x_{0,i} \vee x_{kN,j}=
	\begin{cases}
	x_{kN,i+j}, & \text{if } i+j\le N-2 ,\\
	0, & \text{otherwise}.
    \end{cases}
	\end{equation*}
    In particular, if $\rmchar(\K)\nmid N$, $kN\equiv N-1\mode$ and $i+j=N-1$, $x_{0,i}\vee x_{kN,j}=x_{kN,N-1}\equiv 0$, as $x_{kN,N-1}$ is the coboundary of $\frac{1}{N}y_{(k-1)N+1,0}$.
                
	\item[(ii)] Suppose that $N\not \equiv 1 \mode$ and that $\rmchar (\K)\mid  N$ and $kN\equiv N-1\mode$. Then for any $0\le i\le N-2$ with $i\equiv 0\mode$ and  $0\le j\le N-1$ with   $j\equiv N-1\mode$,
	\begin{equation*}
	x_{0,i} \vee x_{kN,j}=
	\begin{cases}
	x_{kN,i+j}, & \text{if } i+j\le N-1 ,\\
	0, &\text{otherwise}.
	\end{cases}
	\end{equation*} 	
     
	\item[(iii)] Suppose that  $N \equiv 1 \mode$ and that  $\rmchar (\K)\nmid N$ or $kN\not\equiv N-1 \mode$, which means that ($N \equiv 1 \mode$ and    $\rmchar (\K)\nmid N$) or ($N \equiv 1 \mode$ and $e\nmid k$). Then
    
		\begin{itemize}
        
		\item[(a)]  for any $1\le l\le e$ and $0\le j\le N-2$ with $j\equiv kN\mode$,
        \[
        u_{l;\; 0,N-1} \vee x_{k N,j}=0,
        \]
        in particular, if $N=me+1$ for $m\ge 1$, $\rmchar(\K)\nmid N$, $e\di k$ and $j=0$, $u_{l;\; 0,N-1}\vee x_{kN,0}=(\alpha_l^{kN},\alpha_l^{N-1}) \equiv 0$, as $(\alpha_l^{kN},\alpha_l^{N-1})$ is the coboundary of $(\alpha_l^{(k-1)N+1},e_l)-\frac{m}{N}y_{(k-1)N+1,0}$,
                    
		\item[(b)]  for any $0\le i,j\le N-2$ with $i\equiv 0\mode$ and $j\equiv kN\mode$,
		\begin{equation*}
		x_{0,i} \vee x_{kN,j}=
		\begin{cases}
		x_{kN,i+j}, & \text{if } i+j\le N-2 ,\\
		0, &\text{otherwise}.
		\end{cases}
		\end{equation*}
        In particular, if $\rmchar(\K)\nmid N$, $e\di k$ and $i+j=N-1$, $x_{0,i}\vee x_{kN,j}=x_{kN,N-1}\equiv 0$, as $x_{kN,N-1}$ is the coboundary of $\frac{1}{N} y_{(k-1)N+1,0}$.
		\end{itemize}

	\item[(iv)] Suppose that $N \equiv 1 \mode$ and that $\rmchar (\K)\mid  N$ and $kN \equiv N-1 \mode$, which means that $N \equiv 1 \mode$, $\rmchar (\K)\mid  N$ and $e\mid k$. Then
    
		\begin{itemize}
        
		\item[(a)] for any $1\le l\le e$ and $0\le j\le N-1$ with  $j\equiv 0\mode$,
		\begin{equation*}
		u_{l;\; 0,N-1} \vee x_{k N,j}=
		\begin{cases}
		(\alpha_l^{kN},\alpha_l^{N-1})\equiv \frac{1}{e} x_{kN,N-1}, & \text{if } j=0,\\
		0, &\text{otherwise},
		\end{cases}
		\end{equation*}
        as $(\alpha_l^{kN},\alpha_l^{N-1})-(\alpha_1^{kN},\alpha_1^{N-1})$ is the coboundary of $(\alpha_l^{(k-1)N+1},e_l)-(\alpha_1^{(k-1)N+1},e_1)$,
                    
		\item[(b)]  for any $0\le i\le N-2$ and $0\le j\le N-1$ with  $i, j\equiv 0\mode$,
		\begin{equation*}
		x_{0,i} \vee x_{kN,j}=
		\begin{cases}
		x_{kN,i+j}, & \text{if } i+j \le N-1 ,\\
		0, & \text{otherwise}.
		\end{cases}
		\end{equation*}
        
		\end{itemize}
	
	\end{itemize}

\item ($\rmHH^{2h}\vee \rmHH^{2k}$) Assume that $m=2h$ and $n=2k$ with $h,k\ge 1$.

	\begin{itemize}

    \item[(i)] Suppose that $\rmchar(\K)\nmid N$ or $\gcd(N,e)\neq 1 $. Then for any  $0\le i,j\le N-2$ with $i\equiv hN \mode$ and $j\equiv kN \mode$,
	\begin{equation*}
	x_{hN,i} \vee x_{kN,j}=
	\begin{cases}
	x_{(h+k)N,i+j}, & \text{if } i+j\le N-2 ,\\
	0, &\text{otherwise}.
	\end{cases}
	\end{equation*} 	
    In particular, if $\rmchar(\K)\nmid N$, $\gcd(N,e)=1$ and $i+j=N-1$, then $x_{hN,i} \vee x_{kN,j}=x_{(h+k)N,N-1}\equiv 0$.
 
     \item[(ii)] Suppose that $\rmchar(\K)\di N$ and $\gcd(N,e)=1$. Then for any  $0\le i,j\le N-1$ with $i\equiv hN\mode$ and $j\equiv kN\mode$,
	\begin{equation*}
	x_{hN,i} \vee x_{kN,j}=
	\begin{cases}
	x_{(h+k)N,i+j}, & \text{if } i+j\le N-1 ,\\
	0, &\text{otherwise}.
    \end{cases}
	\end{equation*}
		
	\end{itemize}

\item ($\rmHH^{0}\vee \rmHH^{2k+1}$) Assume that $m=0$ and $n=2k+1$, $k\ge 0$.

	\begin{itemize}
    
	\item[(i)] Suppose that $N\not \equiv 1 \mode$ and $\rmchar(\K)\di e$. Then for any $0\le i\le N-2$, $1\le j\le N-1$ with $i\equiv 0\mode$ and $j\equiv kN+1\mode$,
	\begin{equation*}
	x_{0,i} \vee v_{1;\; kN+1,j}=
	\begin{cases}
	v_{1;\; kN+1,i+j}, & \text{if }  i+j\le N-1 ,\\
	0, &\text{otherwise}.
	\end{cases}
	\end{equation*}
				
	\item[(ii)] Suppose that $N\not \equiv 1 \mode$, $\rmchar(\K)\nmid e$  and ($\rmchar(\K)\nmid N$ or $kN+1 \not \equiv 0 \mode$). Then for any $0\le i\le N-2$, $1\le j\le N-1$ with $i\equiv 0\mode$ and $j\equiv kN+1\mode$,
	\begin{equation*}
	x_{0,i} \vee y_{kN+1,j}=
	\begin{cases}
	y_{kN+1,i+j}, & \text{if } i+j\le N-1 ,\\
	0, &\text{otherwise}.
	\end{cases}
	\end{equation*}

	\item[(iii)] Suppose that $N\not \equiv 1 \mode$,  $\rmchar(\K)\di N$ and $kN+1 \equiv 0 \mode$. Then for any $0\le i\le N-2$, $0\le j\le N-1$ with $i\equiv 0\mode$ and $j\equiv 0 \mode$,
	\begin{equation*}
	x_{0,i} \vee y_{kN+1,j}=
	\begin{cases}
	y_{kN+1,i+j}, & \text{if } i+j\le N-1 ,\\
	0, &\text{otherwise}.
	\end{cases}
	\end{equation*}
				
	\item[(iv)] Suppose that $N \equiv 1 \mode$ and $\rmchar(\K)\di e$. Then
    
	  \begin{itemize}
      
	  \item[(a)] for any $1\le l\le e$ and $1\le j\le N-1$ with  $j\equiv kN+1\mode$, 
        \[
        u_{l;\; 0,N-1} \vee v_{1;\; k N+1,j}=0,
        \]
    
		\item[(b)] for any $0\le i\le N-2$, $1\le j\le N-1$ with $i\equiv 0\mode$ and $j\equiv kN+1\mode$,
		\begin{equation*}
		x_{0,i} \vee v_{1;\; kN+1,j}=
		\begin{cases}
		v_{1;\; kN+1,i+j}, & \text{if } i+j\le N-1 ,\\
		0, &\text{otherwise}.
		\end{cases}
		\end{equation*}
        
		\end{itemize}
				
    \item[(v)] Suppose that $N \equiv 1 \mode$, $\rmchar(\K)\nmid e$ and ($\rmchar(\K)\nmid N$ or $kN +1\not \equiv 0 \mode$). Then
    
	  \begin{itemize}
      
		\item[(a)] for any $1\le l\le e$ and $1\le j\le N-1$ with  $j\equiv kN+1\mode$, 
        \[
        u_{l;\; 0,N-1} \vee y_{k N+1,j}=0,
        \]
					
		\item[(b)] for any $0\le i\le N-2$, $1\le j\le N-1$ with $i\equiv 0\mode$ and $j\equiv kN+1\mode$,
		\begin{equation*}
		x_{0,i} \vee y_{kN+1,j}=
		\begin{cases}
		y_{kN+1,i+j}, & \text{if } i+j\le N-1 ,\\
		0, &\text{otherwise}.
		\end{cases}
		\end{equation*}
        
		\end{itemize}
				
	\item[(vi)] Suppose that $N \equiv 1 \mode$, $\rmchar(\K)\di N$ and $kN+1  \equiv 0 \mode$,  which implies, in particular, that $\rmchar(\K)\nmid e$. Then
    
		\begin{itemize}
        
		\item[(a)] for any $1\le l\le e$ and $0\le j\le N-1$ with  $j\equiv 0\mode$,
		\begin{equation*}
		u_{l;\; 0,N-1} \vee y_{k N+1,j} =
		\begin{cases}
		(\alpha_l^{kN+1},\alpha_l^{N-1})\equiv \frac{1}{e} y_{kN+1,N-1}, & \text{if } j= 0,\\
		0, & \text{otherwise},
		\end{cases}
		\end{equation*}
        as $(\alpha_l^{kN+1},\alpha_l^{N-1})-(\alpha_{\ul{l-1}}^{kN+1},\alpha_{\ul{l-1}}^{N-1})$ is the coboundary of $(\alpha_l^{kN},\alpha_l^{N-2})$,
                    
		\item[(b)] for any $0\le i\le N-2$, $0\le j\le N-1$ with $i\equiv 0\mode$ and $j\equiv 0\mode$,
		\begin{equation*}
		x_{0,i} \vee y_{kN+1,j}=
		\begin{cases}
		y_{kN+1,i+j}, & \text{if } i+j\le N-1 ,\\
		0, &\text{otherwise}.
		\end{cases}
		\end{equation*}
		\end{itemize}
		\end{itemize}

\item ($\rmHH^{2h}\vee \rmHH^{2k+1}$) Assume that $m=2h$ and $n=2k+1$, $h\ge 1$ and $k\ge 0$.

	\begin{itemize}
    
	\item[(i)] Suppose  that $\rmchar(\K)\di e$. Then for any $0\le i\le N-2$, $1\le j \le N-1$, with $i\equiv hN\mode$ and $j\equiv kN+1\mode$,
	\begin{equation*}
	x_{hN,i} \vee v_{1;\; kN+1,j}=
	\begin{cases}
	(\alpha_{\ul{1-i}}^{(h+k)N+1},\alpha_{\ul{1-i}}^{i+j})\equiv v_{1;\; (h+k)N+1,i+j}, & \text{if } i+j\le N-1,\\
	0, &\text{otherwise}.
	\end{cases}
	\end{equation*}
				 
	\item[(ii)] Suppose that  $\rmchar(\K)\nmid e$ and ( $\rmchar(\K)\nmid N$ or $\gcd(N,e)\not=1$). Then for any $0\le i\le N-2$, $1\le j \le N-1$, with $i\equiv hN\mode$ and $j\equiv kN+1\mode$,
	\begin{equation*}
	x_{hN,i} \vee  y_{kN+1,j}=
	\begin{cases}
	y_{(h+k)N+1,i+j}, & \text{if } i+j\le N-1 ,\\
	0, &\text{otherwise}.
	\end{cases}
	\end{equation*}
				
	\item[(iii)] Suppose that  $\rmchar(\K)\di N$ and $\gcd(N,e)=1$, which implies, in particular, that  $\rmchar(\K)\nmid e$. Then for any $0\le i,j\le N-1$ with $i\equiv hN\mode$ and $j\equiv kN+1\mode$,
	\begin{equation*}
	x_{hN,i} \vee y_{kN+1,j}=
	\begin{cases}
	y_{(h+k)N+1,i+j}, & \text{if } i+j\le N-1 ,\\
	0, &\text{otherwise}.
	\end{cases}
	\end{equation*}
	\end{itemize}

\item ($\rmHH^{2h+1}\vee \rmHH^{2k+1}$) Assume that $m=2h+1$ and $n=2k+1$ with  $h, k\ge 0$.

	\begin{itemize}
    
	\item[(i)]  Suppose that $\rmchar(\K)\di e$. Then for any $1\le i,j\le N-1$ with $i\equiv hN+1\mode$ and $j\equiv kN+1\mode$,
    \[
    v_{1;\; hN+1,i} \vee  v_{1;\; kN+1,j}=0.
    \]
				
	\item[(ii)] Suppose that $\rmchar(\K)\nmid e$ and ($\rmchar(\K)\nmid N$ or $\gcd(N,e)\not=1$). Then for any $1\le i,j\le N-1$ with $i\equiv hN+1\mode$ and $j\equiv kN+1\mode$,
    \[
    y_{hN+1,i} \vee  y_{kN+1,j}=0.
    \]
				
	\item[(iii)] Suppose that $\rmchar(\K)\di N$ and $\gcd(N,e)=1$. Then for any $0\le i,j\le N-1$ with $i\equiv hN+1\mode$ and $j\equiv kN+1\mode$,
	\begin{equation*}
	y_{hN+1,i} \vee y_{kN+1,j}=
	\begin{cases}
	-\frac{N(N-1)}{2} x_{(h+k+1)N,N-2+i+j}, & \text{if } i+j\le 1,\\
	0, & \text{otherwise}.
	\end{cases}
	\end{equation*}
    
	\end{itemize}
			
\end{enumerate}	

\end{Prop}

\begin{SProof}

We provide proofs for the case (6)(iii), since the others can be obtained by a similar manner or directly from Corollary~\ref{Cor: cup product for truncated basic cycles}.

\textbf{(6)(iii)} Suppose that $\rmchar(\K)\di N$ and $\gcd(N,e)=1$. Following Corollary~\ref{Cor: cup product for truncated basic cycles} (2), we have
\begin{equation*}
	(\alpha_l^{hN+1},\alpha_l^i)\vee (\alpha_{l'}^{kN+1},\alpha_{l'}^j)=
	\begin{cases}
	-\sum\limits_{r=0}^{N-2}a_r(\alpha_{\ul{l-r}}^{(h+k+1)N},\alpha_{\ul{l-r}}^{N-2+i+j}), & \text{if } i+j\le 1,\\
	0, & \text{otherwise},
	\end{cases}
\end{equation*}
where for $i+j\le 1$ and each $r\in [0,N-2]$, $a_r$ is the number of solutions satisfying $s\equiv l'-l-i+r+1(mod\ e)$ and $r+1\leq s\leq N-1$.

It is not difficult to know that $a_r=\left[\frac{N-r-2-\ol{l'-l-i}+e}{e}\right]$, then we have
\begin{equation*}
	\begin{aligned}
	y_{hN+1,i}\vee y_{kN+1,j}&=\sum\limits_{l=1}^e(\alpha_l^{hN+1},\alpha_l^i)\vee\sum\limits_{{l'}=1}^e(\alpha_{l'}^{kN+1},\alpha_{l'}^j)\\
	&=-\sum\limits_{l=1}^e\sum\limits_{{l'}=1}^e\sum\limits_{r=0}^{N-2}\left[\frac{N-r-2-\ol{l'-l-i}+e}{e}\right](\alpha_{\ul{l-r}}^{(h+k+1)N},\alpha_{\ul{l-r}}^{N-2+i+j})\\
	&=-\sum\limits_{l=1}^e\sum\limits_{r=0}^{N-2}(N-r-1)(\alpha_{\ul{l-r}}^{(h+k+1)N},\alpha_{\ul{l-r}}^{N-2+i+j})\\
	&=-\sum\limits_{r=0}^{N-2}(N-r-1)\Big(\sum\limits_{l=1}^e(\alpha_{l}^{(h+k+1)N},\alpha_{l}^{N-2+i+j})\Big)\\
    &=-\frac{N(N-1)}{2}x_{(h+k+1)N,N-2+i+j}.
	\end{aligned}
\end{equation*}
		
\end{SProof}

\begin{Rem}
Without using the above formulae in Corollary~\ref{Cor: cup product for truncated basic cycles} derived from comparison morphisms, Bardzell, Locateli and Marcos \cite{BLM00}  introduced a product $\vee$ on $\bbP^*$ as follows:  for $(\alpha,\pi)\in \bbP^m$ and $(\beta,\tau)\in \bbP^n$,
\begin{equation*}
(\alpha,\pi) \vee (\beta,\tau)=\left\{
\begin{array}{ll}
(\alpha\beta,\pi\tau), &\text{if }  m  \text{ or } n \text{ is even};\\
0,  &\text{otherwise}.
\end{array}
\right.
\end{equation*}
They proved that $\vee$ is an associative product on $\bbP^*$ and it induces an associative product on $\mathrm{HH}^{*}(\Lambda)$. They also claimed that this product coincides with the Yoneda product and hence with the cup product on the Hochschild cohomology groups.

However, by Proposition~\ref{Prop: cup product of basis elements}(6), the cup product of two elements of odd degrees may be nonzero. In fact,    in the setup of Proposition~\ref{Prop: cup product of basis elements} (6)(iii), if $i+j\le 1$,
\[
y_{hN+1,i} \vee y_{kN+1,j}= -\frac{N(N-1)}{2} x_{(h+k+1)N,N-2+i+j} \neq 0, 
\]
whenever $\rmchar(\K)\nmid \frac{N(N-1)}{2}$, i.e., $\rmchar(\K)=2$ and $N\equiv 2\ (\mathrm{mod}\ 4)$.

The error in \cite[Section 4]{BLM00} lies in the fact when $n\geq 1$, $\{\rmHH^{n+i}(A, -), i\geq 0\}$ do NOT form a $\delta$-functor, as $\rmHH^n(A, -)$ with $n\geq 1$ is not left exact. 

\end{Rem}

\begin{Rem} A careful examination shows that the cup product formulae coincident with those in \cite{EH99}.

\end{Rem}

\begin{Rem}
It is well known that the cup product coincides with the Yoneda product, denoted by $\star$, on the Hochschild cohomology groups. In fact, by carefully choosing some liftings, one can show that these two products coincide on the cochain level.

We only consider a special case, the other cases are similar. Assume that $\rmchar \K\di N$ and  $\gcd (N,e)=1$ and let $g=\sum\limits_{l=1}^e(\alpha_l,\alpha_l^i)\in\rmHH^1(\Lambda)$ with $0\le i\le N-1$ and $i\equiv 1\mode$. We further have $i\neq 0$ because $e\neq 1$. We will show that $f\vee  g=f\star g$ for any $f\in \rmHH^*(\Lambda)$.

Let $f\in \rmHH^n(\Lambda)$ represented by a map $f: P_n\to \Lambda$. The Yoneda product $f\star g$ is given by the following commutative diagram:
\begin{eqnarray*}
\xymatrix{
P_{n+1}\ar[r]^{d_{n+1}}\ar[d]^{g_{n+1}} & P_{n}\ar[d]^{g_{n}}\ar[r]^{d_{n}} & \cdots\ar[r]^{d_3} & P_{2}\ar[d]^{g_2}\ar[r]^{d_2} & P_{1}\ar[d]^{g_1}\ar[rd]^{g} \\
P_{n}\ar[r]^{-d_{n}}\ar[d]^{f} & P_{n-1}\ar[r]^{-d_{n-1}} & \cdots\ar[r]^{-d_2} & P_{1}\ar[r]^{-d_1} & P_{0}\ar[r]^{\varepsilon} & \Lambda \\
\Lambda}
\end{eqnarray*}
where the $g_s$, $s\ge 1$ are the lifting of $g$. Then $f\star g=f\circ g_{n+1}$. We will construct explicitly the lifting $g_i$, $1\le i\le n+1$ such that $f\vee g = f \star g$.

For $k\ge 1$, define the maps $g_{2k}:\Lambda \ot \K Q_{kN} \ot \Lambda \to \Lambda \ot \K Q_{(k-1)N+1} \ot \Lambda$ by
\[
g_{2k}(1\ot \alpha_{l'}^{kN}\ot 1)=-\sum\limits_{p=i}^{N-1}\sum\limits_{t=i-1}^{p-1} \alpha_{l'}^t\ot \alpha_{\ul{l'+t}}^{(k-1)N+1}\ot \alpha_{\ul{l'+t+(k-1)N+1}}^{N-2-t+i},
\]
and $g_{2h+1}:\Lambda \ot \K Q_{hN+1} \ot \Lambda \to \Lambda \ot \K Q_{hN} \ot \Lambda$, $h\ge 0$ by
\[
g_{2h+1}(1\ot \alpha_{l'}^{hN+1}\ot 1)=1\ot \alpha_{l'}^{hN}\ot \alpha_{\ul{l'+hN}}^{i}.
\]
It is not hard to check that this is  well-defined.

If $n=2h$ for $h\ge 0$, let $f=\sum\limits_{l=1}^e(\alpha_l^{hN},\alpha_l^j)$ with $0\le j\le N-1$ and $j\equiv hN\mode$. Then we have
\begin{equation*}
		f\circ g_{2h+1}(1\ot \alpha_{l'}^{hN+1}\ot 1)
		=
        \begin{cases}
        \alpha_{l'}^{i+j},     & \text{if }  i+j\le N-1, \\
        0,    & \text{otherwise},
        \end{cases}
\end{equation*}
with $l'=1,\ldots,e$. Then 
\[
f\star g=f\circ g_{2h+1}=
\begin{cases}
\sum\limits_{l=1}^e(\alpha_{l}^{hN+1}, \alpha_{l}^{i+j}),  & \text{if }  i+j\le N-1 \\
0,    & \text{otherwise}
\end{cases}
=f\vee g.
\]

Assume that $n=2h+1$ for $h\ge 0$ and $f=\sum\limits_{l=1}^e(\alpha_l^{hN+1},\alpha_l^j)$ with $0\le j\le N-1$ and $j\equiv hN+1\mode$. Then 
\begin{equation*}
	\begin{aligned}
		f\circ g_{2(h+1)}(1\ot \alpha_{l'}^{(h+1)N}\ot 1)
		&=-\sum\limits_{p=i}^{N-1}\sum\limits_{t=i-1}^{p-1}\alpha_{l'}^{N-2+i+j}\\
		&=
        \begin{cases}
        -\frac{N(N-1)}{2}\alpha_{l'}^{N-1}, & \text{ if } i=1 \text{ and }j=0,\\
        0, & \text{otherwise},
        \end{cases}
	\end{aligned}
\end{equation*}
with $l'=1,\ldots,e$. Then 
\begin{equation*}
f\star g=f\circ g_{2(h+1)}
    =
    \begin{cases}
    -\frac{N(N-1)}{2}\sum\limits_{l=1}^e(\alpha_{l}^{(h+1)N}, \alpha_{l}^{N-1}), & \text{ if } i=1 \text{ and }j=0\\
    0, & \text{otherwise}
    \end{cases}
=f\vee g.
\end{equation*}
\end{Rem}
    
\bigskip

\section{Hochschild cohomology ring}\label{section: Hochschild cohomology ring}\

In this section, we will provide the ring structure of $\rmHH^*(\Lambda)$. Bardzell, Locateli and Marcos determined the ring structure of $\rmHH^*(\Lambda)$ by providing generators and relations {\cite[Section 5]{BLM00}}. However, certain computational inaccuracies occur in \cite{BLM00}, and we present a revised version here.

\medskip

We always let $\Lambda=\K Z_e/{J^N} $ be a truncated basic cycle algebra with $N,\, e\ge 2$. Denote $\gcd(N,e)=d\ge 1$ and $N=dN_1$, $e=de_1$. 
For $0\le j\le e_1-1$, $k_j$ (resp. $I$) is the smallest integer greater than zero such that $k_jN_1\equiv j\ (\mathrm{mod}\  e_1)$ (resp. $(I-1)N+1\equiv 0\mode$). In particular, we have $k_0=e_1$ and $\{k_1,\ldots,k_{e_1-1}\}=\{1,\ldots,e_1-1\}$. Consider the following elements (whenever they exist):
\begin{equation*}
	\begin{aligned}
	 x_{0,0} & =\sum\limits_{l=1}^e (e_l,e_l) =1 \in \rmHH^0(\Lambda), \\
	 x_{0,e} & =\sum\limits_{l=1}^e (e_l,\alpha_l^e) \in \rmHH^0(\Lambda), \\
	 u_{r;\; 0,N-1} & = (e_r,\alpha_r^{N-1}) \in \rmHH^0(\Lambda),\; r=1,\dots,e,\\
	 y_{1,1} & =\sum\limits_{l=1}^e (\alpha_l,\alpha_l) \in \rmHH^{1}(\Lambda), \\
	 v_{1;\; 1,1} &= (\alpha_1,\alpha_1) \in \rmHH^1(\Lambda), \\
	 x_{k_j N,dj} & =\sum\limits_{l=1}^e (\alpha_l^{k_jN},\alpha_l^{dj}) \in \rmHH^{2k_j}(\Lambda),\; j=0,1,\dots,e_1-1,\\
	 y_{(I-1)N+1,0} &=\sum\limits_{l=1}^e ( \alpha_l^{(I-1)N+1},e_l) \in \rmHH^{2I-1}(\Lambda).
	\end{aligned}
\end{equation*}

\begin{Prop}\label{Prop:ring structure1} (Compare with \cite[Theorem 5.4]{BLM00})
Assume that $N\le e$. If $\rmchar(\K) \nmid N$ or $\gcd(N,e)\not= 1$, then $\rmHH^*(\Lambda)$ is isomorphic to the graded commutative algebra generated by $y$ and $z_j$ for $j=0,\dots, [\frac{N-2}{d}]$ with $|y| =1$ and $|z_j| =2 k_j$, subject to the following relations: 
    
\begin{itemize}

\item $y^2=0$, and

\item for any  $1\le i\le j\le [\frac{N-2}{d}]$,
	\begin{equation*}
		z_iz_j=
		\begin{cases}
			z_{i+j}, & \text{if }  i+j\le  [\frac{N-2}{d}] \text{ and } k_i+k_j<e_1,\\
			z_{i+j}z_0  , & \text{if } i+j\le  [\frac{N-2}{d}] \text{ and } k_i+k_j> e_1,\\
			0, & \text{if } i+j> [\frac{N-2}{d}].
		\end{cases}
	\end{equation*}
    
\end{itemize}

In particular, if $N=e$ and $\rmchar(\K)\nmid N$, then $\rmHH^*(\Lambda)$ is isomorphic to the graded commutative polynomial ring $\K [y,z]/{y^2}$ with $|y| =1$ and $|z| =2$; see \cite[Corollary 5.5]{BLM00}.

\end{Prop}

\begin{Proof}
If $\rmchar(\K) \nmid e$. 
Since $N\le e$, it is clear that $\rmHH^0(\Lambda)=\K \{x_{0,0}\}$ and $\rmHH^1(\Lambda)=\K \{y_{1,1}\}$.

For $k\ge 1$, by Theorem~\ref{Thm: basis} (2), $kN\equiv dj\mode$ for some $0\le dj\le N-2$ if $\rmHH^{2k}(\Lambda)\not=0$, that is, $kN_1\equiv j \model$, for some $0\le j\le [\frac{N-2}{d}]$. Following from {\cite[Theorem 5.4]{BLM00}}, we can write $k=k_j+ak_0$ for some $a\ge 0$. Hence, $\rmHH^{2k}(\Lambda)\neq 0$ only when $k\ge 1$ satisfies $kN\equiv dj\mode$ for some $0\le dj\le N-2$, namely, 
\[
k\in \bigcup_{j=0}^{[\frac{N-2}{d}]}(k_j+\N k_0),
\]
and in this case, it is one dimensional generated by 
\[
x_{kN,dj} =x_{(k_j+ak_0)N,dj}={x_{k_0N,0}} ^a \vee x_{k_jN,dj},
\]
where we use ${x_{k_0N,0}} ^a$ to denote $\underbrace{x_{k_0N,0}\vee\cdots\vee x_{k_0N,0}}_{a}$. A similar statement holds for odd degree cohomology groups: $\rmHH^{2k+1}(\Lambda)\neq 0$ only when $k\in \bigcup\limits_{j=0}^{[\frac{N-2}{d}]}(k_j+\mathbb{N} k_0)$, and in this case, it is one dimensional generated by 
\[
y_{kN+1,dj+1}=  y_{1,1} \vee x_{(k_j+ak_0)N,dj}=y_{1,1} \vee {x_{k_0N,0}}^{a}\vee x_{k_jN,dj}.
\]
As an easy consequence,  $\dim \rmHH^n(\Lambda)= 1$ only when $n=0,1$ or
\begin{equation*}
	\begin{array}{ll}
		n=2k \text{ or } 2k+1 \text{ with } 1\le k=k_j+ak_0\text{ for some } 0\le j\le [\frac{N-2}{d}] \text{ and } a\ge 0;
	\end{array}
\end{equation*}
otherwise, it vanishes.

We obtain that $\rmHH^{2k_j}(\Lambda)=\K\{x_{k_jN,dj}\}$ for $j=0, \dots,  [\frac{N-2}{d}]$ together with $\rmHH^1(\Lambda)=\K \{y_{1,1}\}$ generate all positive degree cohomology groups. One could define a homomorphism of graded commutative algebras  from the graded polynomial algebra $\K[y,z_0, \dots, z_{[\frac{N-2}{d}]}]$ with $|y| =1$ and $|z_j| =2k_j, 0\le j\le [\frac{N-2}{d}]$ to $\rmHH^*(\Lambda)$ by  sending  $1$ to $x_{0,0}$, $y$ to $y_{1,1}$ and $z_j$ to $x_{k_jN,dj}$.

We know that ${y_{1,1}}^2=0$ by Proposition~\ref{Prop: cup product of basis elements} (6)(ii). The remain relations are obtained by Proposition~\ref{Prop: cup product of basis elements} (3)(i). In fact, for any  $1
\le i\le j\le [\frac{N-2}{d}]$,
\begin{eqnarray*}
& & x_{k_i N,di} \vee  x_{k_jN, dj}=
x_{(k_i+k_j)N,d(i+j)}=
    \begin{cases}
    x_{k_{i+j}N,d(i+j)},  &\text{if } i+j\le  [\frac{N-2}{d}]\text{ and } k_i+k_j<e_1,\\
    x_{k_{i+j}N,d(i+j)}\vee x_{k_0N, 0},  &\text{if } i+j\le  [\frac{N-2}{d}]\text{ and } k_i+k_j> e_1,\\
    0,&\text{if } i+j>  [\frac{N-2}{d}].
    \end{cases} 
\end{eqnarray*}
Denote by $A$ the graded commutative algebra introduced in the statement of this result. So  it induces a graded  surjective homomorphism  $\varphi$ from $A$ to  $\rmHH^*(\Lambda)$.

Now we show that whenever $\rmHH^n(\Lambda)=0$, $A_n=0$ and when  $\dim(\rmHH^n(\Lambda))=1$, $\dim(A_n)\leq 1$. Together with the surjectivity of $\varphi$, we see that $\varphi$ is an isomorphism. In fact, by direct inspection, it is easy to see that each monomial in $ y, z_i,\;  i=0,\dots, [\frac{N-2}{d}]$ can be turned into the form  $z_0^a z_i,\, y,$ or $yz_0^a z_i$ for $0\le i\le [\frac{N-2}{d}],\; a\ge 0$. Considering the degrees of these elements, we obtain the result.

If $\rmchar(\K)\di e$, the difference with the proof of $\rmchar(\K)\nmid e$ is that since now $\rmHH^1(\Lambda)=\K \{v_{1;\; 1,1}\}$, replace everywhere $y_{kN+1,j}$ by $v_{1;\; kN+1,j}$ will finish the proof. We omit the details.

\end{Proof}

\begin{Prop}\label{Prop:ring structure2} (Compare with \cite[Theorem 5.6]{BLM00})
Assume that $N\le e$. If $\rmchar(\K)\di N$ and $\gcd(N,e)=1$, then $\rmHH^*(\Lambda)$ is isomorphic to the graded commutative algebra generated by $ y,z_j (0\le j\le N-1), w $ with $|y| =1$, $|z_j| =2k_j$ and $|w| =2I-1$, subject to the relations:

\begin{itemize}

\item   $y^2=0$,

\item   $y z_{N-1}=0$,

\item   $w y=-\frac{N(N-1)}{2}z_{N-1}$,

\item   for $1\le i\le j\le N-1$, 	
        \begin{equation*}
	      z_i z_j=
		\begin{cases}
		z_{i+j}, & \text{if } i+j\le N-1 \text{ and } k_i+k_j < e,\\
		z_{i+j}z_0 , & \text{if } i+j\le N-1 \text{ and } k_i+k_j > e,\\
		0, & \text{if } i+j> N-1,
		\end{cases}
	  \end{equation*}
    
\item   for $1\le i \le N-1$,  
    \begin{equation*}
	  w z_i=
	  \begin{cases}
		y z_{i-1}, & \text{if } k_i+I-1\le e,\\
		y z_{i-1} z_0, & \text{if } k_i+I-1> e,
	  \end{cases}
    \end{equation*}
    
\item 	\begin{equation*}
		w^2=
		\begin{cases}
		-\frac{N(N-1)}{2}z_{N-2}, & \text{if } 2I-1\le e,\\
		-\frac{N(N-1)}{2} z_{N-2} z_0, & \text{if } 2I-1> e.
		\end{cases}
	  \end{equation*}
    
\end{itemize}

\end{Prop}

\begin{Proof} 

Note that the conditions that $\rmchar(\K)\di N$ and $\gcd(N,e)=1$ implies $\rmchar(\K) \nmid e$. Since $N\le e$, by Theorem~\ref{Thm: basis} (1) and (3), $\rmHH^0(\Lambda)=\K \{x_{0,0}\}$ and $\rmHH^1(\Lambda)=\K \{y_{1,1}\}$ (by $\rmchar(\K) \nmid e$).

For $k\ge 1$, by Theorem~\ref{Thm: basis} (2), if $\rmHH^{2k}(\Lambda)\neq 0$, then $kN\equiv i\mode$ for some $0\le i\le N-1$. However, since $N\le e$, the integer $i$ is unique (it exists only if $\rmHH^{2k}(\Lambda)\neq 0$). Then we can write $k=k_i+ak_0$ for some $a\ge 0$. Then by Proposition~\ref{Prop: cup product of basis elements} (3), $\rmHH^{2k}(\Lambda)$ is generated by 
\[
x_{kN,i}=x_{(k_i+ak_0)N,i}={x_{k_0N,0}}^{a}\vee x_{k_iN,i}.
\]
Consequently, $\rmHH^{2k_i}(\Lambda)=k\{x_{k_iN,i}\},\;  i=0, \dots, N-1$ generate all  cohomology groups of even degree at least two. Only when $k\ge 1$ satisfies $kN\equiv i\mode$ for some $0\le i\le N-1$, that is, $k\in \bigcup\limits_{i=0}^{N-1}(k_i+\N k_0)$, $\rmHH^{2k}(\Lambda)\neq 0$, and in this case, it is one dimensional generated by $x_{kN,i}={x_{k_0N,0}}^a \vee x_{k_i N,i}$.

The basis of $\rmHH^{2k+1}(\Lambda)$ for $k\ge 1 $ is more complicated. We distinguish it into two cases.

If $kN+1\equiv 0\mode$, we have $\rmHH^{2k+1}(\Lambda)\neq 0$, and $\rmHH^{2 k+1}(\Lambda)$ is one dimensional with basis $\{y_{k N+1, 0}\}$. In this case, $k=I-1+ak_0$ for some $a\ge 0$ and 
\[
y_{kN+1,0}=y_{(I-1+ak_0)N+1,0}={x_{k_0N,0}}^{a}\vee y_{(I-1)N+1,0}.
\]
Note that in this case, $I=k_{N-1}$.

If $kN+1\not \equiv 0\mode$, by Theorem~\ref{Thm: basis} (3), if $\rmHH^{2k+1}(\Lambda)\neq 0$, then $kN+1\equiv i+1\mode$ for unique $i\in\{0,\ldots,N-2\}$. Then $k=k_i+ak_0$ for some $a\ge 0$. By Proposition~\ref{Prop: cup product of basis elements} (3) and (5), $\rmHH^{2k+1}(\Lambda)$ is generated by 
\[
y_{kN+1,i+1}= y_{(k_i+ak_0)N+1,i+1}=y_{1,1} \vee {x_{k_0N,0}}^{a}\vee x_{k_iN,i}.
\]
We obtain that only when $k\in \bigcup\limits_{i=0}^{N-2}(k_i+\N k_0)$, $\rmHH^{2k+1}(\Lambda)\neq 0$ and in this case, it is one dimensional generated by $y_{kN+1,i+1}= y_{1,1}  \vee {x_{k_0N,0}}^{a}\vee x_{k_iN,i}$.

It follows that $\{y_{1,1};\; x_{k_jN,j},\; 0\le j\le N-1;\; y_{(I-1)N+1,0}\}$ generates $\rmHH^*(\Lambda)$ as an algebra. As a consequence,
\[
\dim \rmHH^n(\Lambda)=
\begin{cases} 
1, & n=0 \text{ or } 1,\\
1, & n=2k  \text{ with } 1\le k=k_i+ak_0 \text{ for some } 0\le i \le N-1 \text{ and } a\ge 0, \\
1, & n=2k+1 \text{ with } 1\le k=k_i+ak_0 \text{ for some } 0\le i \le N-2 \text{ and } a\ge 0, \\
1, & n= 2k+1 \text{ with } k=I-1+ak_0 \text{ for some } a\ge 0, \\
0, & \text{otherwise}.
\end{cases}
\]

By Proposition~\ref{Prop: cup product of basis elements} (6)(iii), we obtain the relations ${y_{1,1}}^{2}=0$,
\[
 y_{(I-1)N+1,0} \vee y_{1,1}=-\frac{N(N-1)}{2} x_{IN,N-1}=-\frac{N(N-1)}{2} x_{k_{N-1}N,N-1}
\]
and
\begin{equation*}
 {y_{(I-1)N+1,0}}^{2}=-\frac{N(N-1)}{2} x_{(2 I-1)N,N-2}=
 \begin{cases}
 -\frac{N(N-1)}{2} x_{k_{N-2}N,N-2}, & \text{if } 2I-1\le e,\\
 -\frac{N(N-1)}{2} x_{k_{N-2}N,N-2}\vee x_{k_0N,0} , & \text{if } 2I-1>e.
 \end{cases}
 \end{equation*}
The relations $y_{1,1} \vee x_{k_{N-1},N-1}=0$ and
 \begin{equation*}
 y_{(I-1)N+1,0} \vee x_{k_iN,i}=	
 \begin{cases}
 y_{1,1} \vee x_{k_{i-1},i-1}, & \text{if } k_i+I-1\le e,\\
 y_{1,1} \vee x_{k_{i-1},i-1}\vee x_{k_0N,0} , & \text{if } k_i+I-1>e,
 \end{cases}
 \end{equation*}
for $1\le i\le N-1$, can be obtained from Proposition~\ref{Prop: cup product of basis elements} (5)(iii). Similar to the proof of Proposition~\ref{Prop:ring structure1}, the remain relations can be obtained from Proposition~\ref{Prop: cup product of basis elements} (3)(ii).

So there exists a surjective graded algebra homomorphism from the graded commutative algebra $A$ in the statement to $\rmHH^*(\Lambda)$ sending $y$ to $y_{1,1}$, $z_j$ to $x_{k_jN,j}$ and $w$ to $y_{(I-1)N+1,0}$.

To show that this is an isomorphism, as in the proof of Proposition~\ref{Prop:ring structure1},  it is easy to see that each monomial in $ y, z_i,\;  i=0,\dots, N-1, w $ can be turned into the form $z_0^a z_i,\; a\ge 0, 0\le i\le N-1,$ or  $ y$  or  $ yz_0^a z_i, \; a\ge 0,  0\le i\le N-2$ or $wz_0^a,\; a\ge 0$. This shows that whenever $\rmHH^n(\Lambda)=0$, $A_n=0$ and when  $\dim(\rmHH^n(\Lambda))=1$, $\dim(A_n)\leq 1.$  Together with the surjectivity of $\varphi$, we see that $\varphi$ is an isomorphism. 

\end{Proof}

\begin{Prop}\label{Prop:ring structure3} (Compare with \cite[Theorem 5.7]{BLM00})
Assume that $N> e$ and $N\not\equiv 1\mode$. If $\rmchar(\K) \nmid N$ or $\gcd(N,e)\not= 1$, then $\rmHH^*(\Lambda)$ is isomorphic to the graded commutative algebra  generated by $x_0$, $y$ and $z_j$, $j=0,1,\ldots, e_1-1$ with $| x_0| =0$, $|y| =1$ and $|z_j| =2k_j$. Moreover, these generators satisfy the following relations:

\begin{itemize}

\item $x_0^{[ \frac{N-2}{e}]+1 }=0$,
		
\item $y^2=0$,
		
\item for $[\frac{\,\ol{N-2}\,}{d}]+1\le j\le e_1-1$, $x_0^{[ \frac{N-2}{e}]} z_j=0$,

\item  for $1\le i\le j\le e_1-1$, let $0\le s=\frac{1}{d}\ol{di+dj}\le e_1-1$, $z_i z_j=x_0^{[\frac{di+dj}{e}]}  z_s z_0^{\frac{k_i+k_j-k_s}{k_0}} $.
		
\end{itemize}

\end{Prop}

\begin{Proof}

If $\rmchar(\K) \nmid e$. Since $N\not\equiv 1\mode$, by Theorem~\ref{Thm: basis} (1), 
\[
\rmHH^0(\Lambda)
=\K\{x_{0,ae}\di 0\le a \le [\frac{N-2}{e}]\}.
\]
From Proposition~\ref{Prop: cup product of basis elements} (1)(i), we know that $x_{0,ae}={x_{0,e}}^{a}$, for $1\le a\le [\frac{N-2}{e}]$. Then we have that $x_{0,0}$ and $x_{0,e}$ generate $\rmHH^0(\Lambda)$. Also, by Theorem~\ref{Thm: basis} (3),  we have
\[
\rmHH^1(\Lambda)=\K\{y_{1,ae+1}\di 0\le a \le [\frac{N-2}{e}]\}.
\]
Since by Proposition~\ref{Prop: cup product of basis elements} (4)(ii), $y_{1,ae+1}=y_{1,1}\vee x_{0, ae}=y_{1,1}\vee {x_{0,e}}^{a}$, we know that $y_{1,1}$ and $x_{0,e}$ generate $\rmHH^1(\Lambda)$.

For $k\ge 1$, by Theorem~\ref{Thm: basis} (2), we have 
\[
\rmHH^{2k}(\Lambda)=\K\{x_{kN,i}\di 0\le i\le N-2,\; i\equiv kN\mode\}.
\]
Let $0\le i\le N-2$ such that $i\equiv kN\mode$, we have that
\[
i=\ol{i}+a e=dj+a e
\]	
for some $0\le j\le e_1-1$ and $a\ge 0$, then $kN\equiv dj\mode$, so we have $k=k_j+b k_0$ for some $b\ge 0$. It follows from Proposition~\ref{Prop: cup product of basis elements} (2)(i) and (3)(i) that
\[
x_{kN,i}=x_{(k_j+b k_0)N ,dj+a e} ={x_{0,e}}^{a}\vee {x_{k_0N,0}}^{b}\vee x_{k_jN,dj}.
\]
Consequently, $x_{0,e}$ and $x_{k_jN,dj}$, $j=0,1,\ldots,e_1-1$ generate all even degree cohomology groups.

By Theorem~\ref{Thm: basis} (3),
\[
\rmHH^{2k+1}(\Lambda)=\K\{y_{kN+1,i+1}\di 0\le i\le N-2,\ i\equiv kN\mode\}.
\]
By 	Proposition~\ref{Prop: cup product of basis elements} (5)(ii), for $0\le i\le N-2 $ with $ i\equiv kN\mode$,
\[
y_{kN+1,i+1}=y_{(k_j+bk_0)N+1,dj+ae+1}=y_{1,1}\vee {x_{0,e}}^{a}\vee {x_{k_0 N,0}}^{b}\vee x_{k_j N,d j},
\]
where $i=\ol{i}+a e=dj+a e$	for some $0\le j\le e_1-1$ and $a\ge 0$, and  $k=k_j+b k_0$ for some $b\ge 0$. It follows that $y_{1,1}$, $x_{0,e}$ and $x_{k_jN,dj}$, $j=0,1,\ldots,e_1-1$ generate all odd degree cohomology groups.

We conclude that the set  $\{x_{0,e},y_{1,1}, x_{k_jN,dj}\di 0\le j\le e_1-1\}$ generates $\rmHH^*(\Lambda)$ as an algebra. As a consequence, we have $\dim \rmHH^{n}(\Lambda) =[\frac{N-2-\ol{kN}}{e}]+1$
if $n=2k$ or $n=2k+1$ for $k\ge 0$.

By Proposition~\ref{Prop: cup product of basis elements} (6)(ii), we know that ${y_{1,1}}^{2}=0$. It follows from Proposition~\ref{Prop: cup product of basis elements} (1)(i) that
\[
{x_{0,e}}^{[\frac{N-2}{e}]}\not=0\text{ and } {x_{0,e}}^{[\frac{N-2}{e}]+1}=0,
\]
because $e[\frac{N-2}{e}]\le N-2$ and $e[\frac{N-2}{e}]+e\ge N-1$.
For $[\frac{\,\ol{N-2}\,}{d}]+1\le j\le e_1-1$, by Proposition~\ref{Prop: cup product of basis elements} (2)(i), 
\[
{x_{0,e}}^{[\frac{N-2}{e}]}\vee x_{k_jN,dj}=0,
\]
since $e[\frac{N-2}{e}]+d([\frac{\,\ol{N-2}\,}{d}]+1)\ge N-1$ and $e[\frac{N-2}{e}]+d[\frac{\,\ol{N-2}\,}{d}]\le N-2$.

For $1\le i\le j\le e_1-1$,  by Proposition~\ref{Prop: cup product of basis elements} (3)(i),
\[
x_{k_iN,di} \vee x_{k_jN, dj}=x_{(k_i+k_j)N,di+dj}.
\]
Now let $s=\frac{1}{d}\ol{di+dj}$, then $k_i+k_j=k_s+ \frac{k_i+k_j-k_s}{k_0} k_0$ and  $di+dj=ds+[\frac{di+dj}{e}]e$, so
\[
x_{k_iN,di} \vee x_{k_jN, dj} =x_{k_sN,ds} \vee {x_{0,e}}^{[\frac{di+dj}{e}]} \vee {x_{k_0N,0}}^{\frac{k_i+k_j-k_s}{k_0}}.
\]
It is easy to see that the set of relations above consititutes a Gr\"{o}bner basis with respect to the lexicographic order $y_{1,1}>x_{0,e}>x_{k_1 N,1}>\cdots >x_{k_{e_1-1}N,e_1-1}>x_{k_0 N,0} $.
Denote by $A$ the graded commutative  algebra introduced in the statement of this result. So there exists an isomorphism graded algebra homomorphism $\varphi$  from $A$ to $\rmHH^*(\Lambda)$  sending  $x_0$ to $x_{0,e}$, $y$ to $y_{1,1}$, and $z_j$ to $x_{k_jN, dj}$, $j=0,1,\ldots,e_1-1$.

If $\rmchar(\K) \di e$, the difference with the proof of $\rmchar(\K) \nmid e$ is that since now
\[
\rmHH^{2k+1}(\Lambda)=\K \{v_{1;\;kN+1,i+1}\di 0\le i\le N-2,\ i\equiv kN\mode\},
\]
replace everywhere $y_{kN+1,i+1}$ by $v_{1;\;kN+1,i+1}$ will finish the proof. We omit the details.

\end{Proof}

\begin{Ex}\label{Example: generators and relations}

Let $\K$ be a field of characteristic two. Consider the algebra $\Lambda$ defined by the quiver with relations
\[
\xymatrix{1\ar@<1ex>[r]^\alpha & 2\ar@<1ex>[l]^{\beta},  & \alpha\beta\alpha\beta=0=\beta\alpha\beta\alpha}.
\]
Then we have $\rmHH^*(\Lambda)\cong \K[x,y,z]/\langle x^2,y^2\rangle$ with $|x| =0$, $|y| =1$ and $|z| =2$. This isomorphism maps $x$ to $(e_1,\alpha\beta)+(e_2,\beta\alpha)$; $y$ to $(\alpha,\alpha)$; $z$ to $((\alpha\beta)^2,e_1)+((\beta\alpha)^2,e_2)$.
    
\end{Ex}

\begin{Prop}\label{Prop:ring structure4} (Compare with \cite[Theorem 5.8]{BLM00})
Assume that $N>e$, $N\not\equiv 1\mode$, $\rmchar(\K)\di N$ and $\gcd(N,e)=1$. Then $\rmHH^*(\Lambda)$ is isomorphic to the graded commutative algebra generated by $\{x_0,y,z_j,w\di 0\le j\le e-1\}$ with  $| x_0| =0$, $|y| =1$,  $| z_j| =2k_j$ and  $|w| =2I-1$, subject to the following relations:

\begin{itemize}

\item $x_0^{[ \frac{N-1}{e}]+1 }=0$,

\item $y^2=0$,

\item for $\ol{N-1}+1 \le i\le e-1$, $x_0^{[\frac{N-1}{e}]} z_i=0$,

\item $w x_0=y z_{e-1}$,

\item $w y=-\frac{N(N-1)}{2}x_0^{[\frac{N-1}{e}]} z_{\ol{N-1}}$,

\item $y x_0^{[\frac{N-1}{e}]} z_{\ol{N-1}}=0$,
	
\item for $1\le i\le j\le e-1$, $z_i z_j=x_0^{[\frac{i+j}{e}]} z_s z_0^{\frac{k_i+k_j-k_s}{k_0}} $ with $0\le s=\ol{i+j}\le e-1$,
		
\item for $1\le i\le e-1$, $w z_i=y z_{i-1} z_0^{\frac{k_i+I-1-k_{i-1}}{k_0}} $,

\item $w^2=-\frac{N(N-1)}{2}x_0^{[\frac{N-1}{e}]} z_{\ol{N-2}} z_0^{\frac{2I-1-k_{\ol{N-2}}}{k_0}} $.

\end{itemize}
	
\end{Prop}

\begin{Proof}

The spaces $\rmHH^0(\Lambda)$ and $\rmHH^1(\Lambda)$ are described in the same way as in the proof of Proposition~\ref{Prop:ring structure3} for the case $\rmchar(\K)\nmid e$.

For $k\ge 1$, by Theorem~\ref{Thm: basis} (2), 
\[
\rmHH^{2k}(\Lambda)=\K\{x_{kN,i}\di 0\le i\le N-1,\; i\equiv kN\mode\},
\]
it has dimension $[\frac{N-1-\ol{kN}}{e}]+1$. Let $0\le i\le N-1$ satisfying $i\equiv kN\mode$, we have
\[
	x_{kN,i}={x_{0,e}}^{a}\vee {x_{k_0N,0}}^{b}\vee x_{k_j N,j},
\]
where $i=j+a e$	for some $0\le j\le e-1$ and $a\ge 0$, and  $k=k_j+b k_0$ for some $b\ge 0$. Consequently, (products of) $x_{0,e}$ and $x_{k_jN,j}$, $j=0,1,\dots,e-1$ generate all even degree cohomology groups.

The case of $\rmHH^{2k+1}(\Lambda), k\ge 1 $ is more complicated. We distinguish it into two cases.

When $kN+1\not\equiv 0\mode$, by Theorem~\ref{Thm: basis} (3),
\[
\rmHH^{2k+1}(\Lambda)=\K\{y_{kN+1,i+1}\di 0\le i\le N-2,\ i\equiv kN\mode\}.
\]
In this case, the dimension of $\rmHH^{2k+1}(\Lambda)$ is $[\frac{N-2-\ol{kN}}{e}]+1$. If $0\le i\le N-2$ and $i\equiv kN\mode$, then 
\[
y_{kN+1,i+1}=y_{1,1}\vee {x_{0,e}}^{a}\vee {x_{k_0 N,0}}^{b}\vee x_{k_j N,j},
\]
where $i=j+a e$	for some $0\le j\le e-2$ and $a\ge 0$, and  $k=k_j+b k_0$ for some $b\ge 0$. So in this case, $\rmHH^{2k+1}(\Lambda)$ is generated by (products of) $y_{1,1}$, $x_{0,e}$ and $x_{k_jN,j}$, $j=0,1,\dots, e-2$.

Consider the case where $kN+1\equiv 0\mode$. Then $k\equiv I-1\equiv k_{e-1}\mode $. We have
\begin{equation*}
	\begin{array}{ll}
	\rmHH^{2k+1}(\Lambda)
    &=\K\{y_{kN+1,i}\di 0\le i\le N-1,\; i\equiv 0\mode\}\\
	&=\K\{y_{kN+1,ae}\di 0\le a \le [\frac{N-1}{e}]\}.
	\end{array}
\end{equation*}
The dimension of $\rmHH^{2k+1}(\Lambda)$ is $[\frac{N-1}{e}]+1$. By definition, we have $k=I-1
+c k_0$, for some $c\ge 0$. By Proposition~\ref{Prop: cup product of basis elements} (5)(iii) and (4)(iii),
\[
y_{kN+1,ae}=y_{(I-1)N+1+c k_0N, ae}={x_{0,e}}^a\vee {x_{k_0N,0}}^{c}\vee y_{(I-1)N+1,0}.
\]
In this case, $\rmHH^{2k+1}(\Lambda)$ is generated by (products of) $x_{0,e}$, $x_{k_0N,0}$ and $y_{(I-1)N+1,0}$. Note that in this case, since $I-1 =k_{e-1}$, for $1\le a \le [\frac{N-1}{e}]$, we also have 
\[
y_{kN+1,ae} =y_{1, 1} \vee {x_{0,e}}^{a-1} \vee {x_{k_0N,0}}^{c}\vee  x_{k_{e-1}N,e-1}.
\]
In particular, we obtain the relation
\begin{equation} \label{eq:5.7 f}  
y_{(I-1)N+1,0} \vee x_{0,e}=y_{1, 1} \vee x_{k_{e-1}N,e-1}.
\end{equation}

Then we obtain that $x_{0,e}$, $y_{1,1}$, $x_{k_jN,j}$, $j=0,1,\dots,e-1$ and $y_{(I-1)N+1,0}$ generate $\rmHH^*(\Lambda)$ as a graded algebra. In particular,
\begin{equation*}
	\begin{array}{rl}\dim(\rmHH^{n}(\Lambda))
    &=
	\begin{cases}
    [\frac{N-2}{e}]+1,&\text{if } n=0,\, 1,\\
	{[\frac{N-1-\ol{kN}}{e}]+1},&\text{if } n=2k,\; k\ge 0,\\
	{[\frac{N-2-\ol{kN}}{e}]+1},&\text{if } n=2k+1,\; k\ge 1 \text{ and } kN+1\not\equiv 0\mode,\\
    {[\frac{N-1}{e}]+1},&\text{if } n=2k+1,\; k\ge 1 \text{ and } kN+1 \equiv 0\mode.
	\end{cases}
    \end{array}
\end{equation*}

\medskip

Now we consider relations among the generators.

By Proposition~\ref{Prop: cup product of basis elements} (1)(i), we have
\[
{x_{0,e}}^{[\frac{N-2}{e}]}\not=0\text{ and } {x_{0,e}}^{[\frac{N-2}{e}]+1}=0.
\]
In here, we have $[\frac{N-2}{e}]=[\frac{N-1}{e}]$ because of $N\not \equiv 1\mode$. We have obtained the relation Equation~\eqref{eq:5.7 f}:
\[
y_{(I-1)N+1,0} \vee x_{0,e} = y_{1, 1} \vee x_{k_{e-1}N,e-1}.
\]

By Proposition~\ref{Prop: cup product of basis elements} (6)(iii), we have relations: ${y_{1,1}}^2=0$;  since $I=k_{\ol{N-1}}$, by Proposition~\ref{Prop: cup product of basis elements} (6)(iii) and (2)(ii),
\begin{equation*}
	\begin{array}{ll}
	y_{(I-1)N+1,0}\vee y_{1,1}=-\frac{N(N-1)}{2}x_{IN,N-1}=-\frac{N(N-1)}{2}{x_{0,e}}^{[\frac{N-1}{e}]}\vee x_{k_{\ol{N-1}}N,\ol{N-1}};
	\end{array}
\end{equation*}
since $N-2=\ol{N-2}+[\frac{N-2}{e}]e$ and $2I-1=k_{\ol{N-2}}+\frac{2I-1-k_{\ol{N-2}}}{k_0}k_0$, by Proposition~\ref{Prop: cup product of basis elements} (6)(iii) and (2)(ii),
\begin{equation*}
    \begin{array}{ll}
    {y_{(I-1)N+1,0}}^2=-\frac{N(N-1)}{2} x_{(2I-1)N,N-2}=-\frac{N(N-1)}{2} {x_{0,e}}^{[\frac{N-2}{e}]}\vee x_{k_{\ol{N-2}}N,\ol{N-2}}\vee {x_{k_0N,0}}^{\frac{2I-1-k_{\ol{N-2}}}{k_0}}.
    \end{array}
\end{equation*}

Following from Proposition~\ref{Prop: cup product of basis elements} (3)(ii), for $0\le i\le \ol{N-1}$ and $\ol{N-1}+1\le j\le e-1$, since $e[\frac{N-1}{e}]+i\le N-1$ and $e[\frac{N-1}{e}]+j\ge N$, by Proposition~\ref{Prop: cup product of basis elements} (3)(ii),
\[
{x_{0,e}}^{[\frac{N-1}{e}]}\vee x_{k_iN,i}\not=0 \text{ and } {x_{0,e}}^{[\frac{N-1}{e}]}\vee x_{k_jN,j}=0.
\]
For $1\le i\le j\le e-1$, let $0\le s=\ol{i+j}\le e-1$, by Proposition~\ref{Prop: cup product of basis elements} (3)(ii) and (2)(ii),
\[
x_{k_iN,i} \vee x_{k_jN,j}={x_{0,e}}^{[\frac{i+j}{e}]} \vee x_{k_sN,s}\vee {x_{k_0N,0}}^{\frac{k_i+k_j-k_s}{k_0}}.
\]
By Proposition~\ref{Prop: cup product of basis elements} (5)(iii), as $1+e[\frac{N-1}{e}]+\ol{N-1}= N$,
\[
y_{1,1} \vee {x_{0, e}}^{[\frac{N-1}{e}]} \vee x_{k_{\ol{N-1}}N, \ol{N-1}}=0,
\]
and for $1\le i\le e-1$,
\[
y_{(I-1)N+1,0}\vee x_{k_iN,i}=y_{(k_i+I-1)N+1,i}=y_{1,1}\vee x_{(k_i+I-1)N,i-1}=y_{1,1}\vee x_{k_{i-1}N,i-1} \vee {x_{k_0N,0}}^{\frac{k_i+I-1-k_{i-1}}{k_0}}.
\]

It is easy to see that the set of relations above consititutes a Gr\"{o}bner basis with respect to the lexicographic order $y_{(I-1)N+1,0}>y_{1,1}>x_{0,e}>x_{k_1 N,1}>\cdots >x_{k_{e-1}N,e-1}> x_{k_0 N,0} $.
Therefore, there exists an isomorphism of graded algebras $\phi$ from the graded commutative algebra in the statement, which we denote by $A$, to $\rmHH^*(\Lambda)$ by sending $x_0$ to $x_{0,e}$, $y$ to $y_{1,1}$, $z_j$ to $x_{k_jN,j}$, $j=0,1,\dots,e-1$ and $w$ to $y_{(I-1)N+1,0}$.

\medskip

\end{Proof}

In the previous $4$ propositions we considered the case $N\not\equiv 1\mode$. In the following $4$ propositions (Compare to \cite[Theorem 5.10]{BLM00}) we will instead focus on the case $N\equiv 1\mode$, and describe the generators and relations of $\rmHH^*(\Lambda)$ as a graded algebra.

We notice that in this case $k_j=j$ for any $1\le j\le e-1$ and  $I=k_0=e$.

\begin{Prop}\label{Prop: ring structure5}

	Assume that $N=e+1$ and $\rmchar (\K)\nmid N$. Then the algebra $\rmHH^*(\Lambda)$ is isomorphic to the graded commutative algebra generated by $\widetilde{x_l}$, $1\le l\le e$, $y$, $z$ and $z'$ with $| \widetilde{x_l}| =0$, $|y| =1$,  $|z| =2$, and  $| z'| =2e$, subject to the relations:
	\begin{itemize}
		\item $y^2=0$,
		
		\item for $1\le l\le p\le e$, $\widetilde{x_l}\widetilde{x_p}=0$,
		
		\item for $1\le l\le e$,  $y\widetilde{x_l}=0$,   $\widetilde{x_l}z=0$,  and $\widetilde{x_l}z'=0$,
		
		\item $z^{e}=0$.
	\end{itemize}
    
\end{Prop}

\begin{Proof}

The proof closely follows that the Proposition~\ref{Prop: ring structure 7} and is simpler; hence, we omit it.

\end{Proof}

\begin{Prop}\label{Prop: ring structure 6}

Assume that $N=e+1$ and $\rmchar (\K)\di N$. Then the algebra $\rmHH^*(\Lambda)$ is generated by 	$\widetilde{x_l},\, 1\le l\le e$, $y$, $z$, $w$, and $z'$, with $| \widetilde{x_l}| =0$, $|y| =1$,  $|z| =2$, $|w| =2e-1$, and $| z'| =2e$. Furthermore, these generators satisfy the relations:

\begin{itemize}

	\item  $y^2=0$,
		
	\item $w y=-\frac{N(N-1)}{2}z^e $,
		
	\item $w^2=-\frac{N(N-1)}{2} z^{e-1} z'$,
		
	\item for $1\le l\le p\le e$, $\widetilde{x_l}\widetilde{x_p}=0$,

    \item for $1\le l\le e-1$, $\widetilde{x_l}z'=\widetilde{x_e}z'$,
		
	\item for $1\le l\le e$, $z\widetilde{x_l}=0$, $y \widetilde{x_l}=0$ and $w \widetilde{x_l}=\frac{1}{e} y z^{e-1}$,

    \item $z^e=\sum_{l=1}^e \widetilde{x_l}z'$,
		
	\item $w z=yz'$.

\end{itemize}

\end{Prop}

\begin{Proof}

The proof closely follows that the Proposition~\ref{Prop: ring structure 8} and is simpler; hence, we omit it.

\end{Proof}

\begin{Prop}\label{Prop: ring structure 7}

Assume that $N=me+1$, $m> 1$ and $\rmchar (\K)\nmid N$. Then the algebra $\rmHH^*(\Lambda)$ is isomorphic to the  graded commutative algebra generated by $x_0$,   $\widetilde{x_l}$ ($1\le l\le e$), $y$, $z$, and $z'$ with $| x_0| =0$, $| \widetilde{x_l}| =0$, $|y| =1$, $|z| =2$, and $| z'| =2e$, which satisfy the relations:

\begin{itemize}
			
	\item $x_0^m=\sum\limits_{l=1}^e\widetilde{x_l}$,
			
    \item $z^e=x_0 z' $,

	\item for $1\le l\le p\le e$, $\widetilde{x_l}\widetilde{x_p}=0$,
			
	\item for $1\le l\le e$,  $x_0 \widetilde{x_l}=0$,  $y \widetilde{x_l}=0$,  $z \widetilde{x_l}=0$, and $\widetilde{x_l}z'=0$,
			
	\item $y^2=0$.
    
\end{itemize}
	
\end{Prop}

\begin{Proof}

Assume first $\rmchar (\K)\nmid e$.
Since $N\equiv 1\mode$, by Theorem~\ref{Thm: basis} (1),
\[
\rmHH^0(\Lambda)=\K\{x_{0,ae}=x_{0,e}^a,\, u_{l;\; 0,N-1} \di 0\le a\le m-1 ,\, 1\le l\le e\}.
\]
Then $x_{0,e}$ and $u_{l;\; 0,N-1}$, $1\le l\le e$ generate $\rmHH^0(\Lambda)$.
Also, by Theorem~\ref{Thm: basis} (3), 
\[
\rmHH^1(\Lambda)=\K\{y_{1,ae+1}=y_{1,1}\vee x_{0,e}^a \di 0\le a\le m-1\}.
\]
So (products of) $x_{0,e}$ and $y_{1,1}$ generate $\rmHH^1(\Lambda)$.

For $k\ge 1$, by Theorem~\ref{Thm: basis} (2), we have
\[
\rmHH^{2k}(\Lambda)=\K\{x_{kN,i}\di 0\le i \le me-1,\; i\equiv k\mode\}=\K\{x_{kN,ae+\ol{k}}\di 0\le a\le m-1\},
\]
and $x_{kN,ae+\ol{k}}=x_{0,e}^a \vee x_{N,1}^{\ol{k}}\vee x_{eN,0}^{[\frac{k}{e}]}$. Consequently, $\rmHH^{2k}(\Lambda)$ is $m$-dimensional and (products of)  $x_{0,e}\in\rmHH^0(\Lambda)$, $x_{N,1}\in\rmHH^2(\Lambda)$ and $x_{eN,0}\in\rmHH^{2e}(\Lambda)$ generate all cohomology groups of even degree at least two.

For $k\ge 1$, by Theorem~\ref{Thm: basis} (3), we have
\[
\rmHH^{2k+1}=\K\{y_{kN+1,i+1}\di 0\le i\le me-1,\; i\equiv k\mode\}=\K\{y_{kN+1,ae+\ol{k}+1}\di 0\le a\le m-1\},
\]
and $y_{kN+1,ae+\ol{k}+1}=y_{1,1}\vee  x_{0,e}^a \vee x_{N,1}^{\ol{k}}\vee x_{eN,0}^{[\frac{k}{e}]}$. Then $\rmHH^{2k+1}(\Lambda)$ is $m$-dimensional, and (products of) $y_{1,1}\in \rmHH^1(\Lambda)$, $x_{0,e}\in\rmHH^0(\Lambda)$, $x_{N,1}\in\rmHH^2(\Lambda)$ and $x_{eN,0}\in\rmHH^{2e}(\Lambda)$ generate all cohomology groups of odd degree.

It follows that $\{x_{0,e}, u_{l;\; 0,N-1}, y_{1,1}, x_{N,1}, x_{eN,0}\di 1\le l\le e\}$ generates $\rmHH^*(\Lambda)$ as an algebra. Moreover, we obtain that
\begin{equation*}
	\dim(\rmHH^n(\Lambda))=
	\begin{cases}
	m+e, & \text{if } n=0,\\
	m, & \text{if } n\neq 0.
	\end{cases}
\end{equation*}

By Proposition~\ref{Prop: cup product of basis elements} we obtain the following relations:
\begin{itemize}
    \item $x_{0,e}^m= \sum\limits_{l=1}^e u_{l;\; 0,N-1}$, 

    \item $u_{l;\; 0,N-1} \vee u_{p;\; 0,N-1}=0$ for $1\le l\le p\le e$,

    \item $x_{0,e}\vee u_{l;\; 0,N-1} =0$, $y_{1,1} \vee u_{l;\; 0,N-1}=0$, $x_{N,1} \vee u_{l;\; 0,N-1}=0$, and $u_{l;\; 0,N-1} \vee x_{eN,0}=0$,  for $1\le l\le e$,

    \item $y_{1,1}^2=0$,

    \item $x_{N,1}^e=x_{0,e} \vee x_{eN,0}$.
\end{itemize}
It is easy to see that the set of relations above constitutes a Gr\"{o}bner basis with respect to the lexicographic order $y_{1,1}> x_{N,1}>x_{0,e}>u_{1;\; 0,N-1}>\cdots >u_{e;\; 0,N-1}>x_{eN,0}$. Denote by $A$ the graded commutative algebra introduced in the statement. Then there exists a graded algebra isomorphism $\phi$ from $A$ to $\rmHH^*(\Lambda)$ sending $x_0$ to $x_{0,e}$, $\widetilde{x_l}$ to $u_{l;\; 0,N-1}$, $l=1,\dots,e$, $y$ to $y_{1,1}$, $z$ to $x_{N,1}$, and $z'$ to $x_{eN,0}$.

If $\rmchar (\K)\di e$, the difference with the proof above is that since now $\rmHH^{1}(\Lambda)=\K\{v_{1;\; 1,ae+1}\di 0\le a\le m-1\}$, replace everywhere $y_{kN+1,i}$ by $v_{1;\; kN+1,i}$ will finish the proof. We omit the details.

\end{Proof}

\begin{Prop}\label{Prop: ring structure 8}

Assume that $N=me+1$, $m> 1$ and $\rmchar (\K)\di N$. Then the algebra $\rmHH^*(\Lambda)$ is generated by $x_0$, $\widetilde{x_l}$, $1\le l\le e$, $y$, $z$, $w$ and $z'$ with $| x_0| =| \widetilde{x_l}| =0$, $| y|=1 $, $| z'| =2e$, $|z| =2$ and $|w| =2e-1$. Furthermore, these generators satisfy the relations:

\begin{itemize}

    \item $x_0^m=\sum\limits_{l=1}^e\widetilde{x_l}$,	

    \item $w x_0=y z^{e-1}$,

	\item $y^2=0$,
		
	\item $w y=-\frac{N(N-1)}{2}e \widetilde{x_e}z'$,
		
	\item $w^2=-\frac{N(N-1)}{2} z^{e-1} x_0^{m-1}z'$,
		
	\item for $1\le l\le p\le e$, $\widetilde{x_l}\widetilde{x_p}=0$,
		
	\item for $1\le l\le e$, $x_0\widetilde{x_l}=0$, $z\widetilde{x_l}=0$, $y \widetilde{x_{l}} =0$ and $w\widetilde{x_l}=\frac{1}{e} y z^{e-1}x_0^{m-1} $,

    \item for $1\le l\le e-1$, $\widetilde{x_l}z'=\widetilde{x_e}z'$,
		
	\item $z^e=x_0  z' $,
		
	\item $w z=yz'$.
    
	\end{itemize}
    
\end{Prop}

\begin{Proof}

As in the proof of Proposition~\ref{Prop: ring structure 7} in the case $\rmchar(\K)\nmid e$, we obtain that $x_{0,e}$ together with $u_{l;\; 0,N-1}$ for $1\le l\le e$ generate $\rmHH^0(\Lambda)$. Moreover, $\rmHH^1(\Lambda)$ is generated by (products of) $x_{0,e}$ and $y_{1,1}$.

For $k\ge 1$, by Theorem~\ref{Thm: basis} (2), if $kN\equiv N-1\mode$, $\ie\ k\equiv 0\mode$, we have
\[
\rmHH^{2k}(\Lambda)=\K\{x_{kN,ae}=x_{0,e}^a \vee x_{e N,0}^{\frac{k}{e}}\di 0\le a\le m\}.
\]
In this case, $\rmHH^{2k}(\Lambda)$ is $(m+1)$-dimensional.
If $k\not\equiv 0\mode$,
\[
\rmHH^{2k}(\Lambda)=\K\{x_{kN,ae+\ol{k}}=x_{0,e}^a \vee x_{N,1}^{\ol{k}} \vee x_{eN,0}^{[\frac{k}{e}]}\di 0\le a\le m-1\}.
\]
In this case, $\rmHH^{2k}(\Lambda)$ is $m$-dimensional. Moreover, (products of) $x_{0,e}\in\rmHH^0(\Lambda)$, $x_{eN,0}\in\rmHH^{2e}(\Lambda)$ and $x_{N,1}\in\rmHH^{2}(\Lambda)$ generate all cohomology groups of even degree at least two.

For $k\ge 1$, by Theorem~\ref{Thm: basis} (3), if $kN+1\equiv 0\mode$, $\ie,\, k\equiv e-1\mode$, we have
\[
\rmHH^{2k+1}=\K\{y_{kN+1,ae}=y_{(e-1)N+1,0}\vee x_{0,e}^a \vee x_{e N,0}^{[\frac{k}{e}]}\di 0\le a\le m\}.
\]
Then in this case $\rmHH^{2k+1}(\Lambda)$ is $m+1$-dimensional. If $kN+1\not \equiv 0\mode$,
\[
\rmHH^{2k+1}=\K\{y_{kN+1,ae+\ol{k}+1}=y_{1,1} \vee x_{0,e}^a \vee x_{N,1}^{\ol{k}}\vee x_{eN,0}^{[\frac{k}{e}]}\di 0\le a\le m-1\}.
\]
In this case, $\rmHH^{2k+1}(\Lambda)$ is $m$-dimensional. Consequently, $x_{0,e}\in\rmHH^0(\Lambda)$, $y_{1,1}\in \rmHH^1(\Lambda)$, $x_{eN,0}\in\rmHH^{2e}(\Lambda)$, $ x_{N,1}\in\rmHH^2(\Lambda)$ and $y_{(e-1)N+1,0}\in\rmHH^{2e-1}(\Lambda)$ generate all cohomology groups of odd degree.

It follows that $\{x_{0,e},u_{l;\; 0,N-1},y_{1,1}, x_{N,1}, x_{eN,0}, y_{(e-1)N+1,0}\di 1\le l\le e\}$ generates $\rmHH^*(\Lambda)$ as an algebra. As a consequence, we obtain that
\begin{equation*}
	\dim(\rmHH^n(\Lambda))=
    \left\{
	\begin{array}{ll}
	m+e, & \text{if } n=0,\\
	m+1, & \text{if } n=2k \text{ and } k\equiv0\mode  \text{ or } n=2k+1 \text{ and } k\equiv e-1\mode,\\
	m, & \text{if } n=2k \text{ and } k\not\equiv0\mode \text{ or } n=2k+1 \text{ and } k\not \equiv e-1\mode.
	\end{array}
	\right.
\end{equation*}

Now we consider relations among generators. By Proposition~\ref{Prop: cup product of basis elements}, we obtain the following relations:
\begin{itemize}
    \item $x_{0,e}^m =\sum\limits_{l=1}^e u_{l;\; 0,N-1}$,

    \item ${x_{N,1}}^e= x_{0,e} \vee x_{eN,0}$,

    \item ${y_{1,1}}^2=0$,

    \item $y_{(e-1)N+1,0}^2 =-\frac{N(N-1)}{2} x_{N,1}^{e-1}\vee x_{0,e}^{m-1} \vee x_{e N,0}$,

    \item $y_{(e-1)N+1,0} \vee y_{1,1} =-\frac{N(N-1)}{2} x_{0,e}^m \vee x_{e N,0}$,

    \item $y_{(e-1)N+1,0}\vee x_{N,1}=y_{1,1}\vee x_{eN,0} $,

    \item $y_{(e-1)N+1,0}\vee x_{0,e}= y_{1,1} \vee x_{N,1}^{e-1}$,

    \item for any $1\le l\le p\le e$, $u_{l;\; 0,N-1} \vee u_{p;\; 0,N-1}=0$ and $u_{l;\; 0,N-1} \vee x_{e N,0}= u_{p;\; 0,N-1} \vee x_{e N,0} $,

    \item for any $1\le l\le e$, $x_{0,e} \vee u_{l;\; 0,N-1}=0$, $y_{1,1} \vee u_{l;\; 0,N-1} =0$, and $x_{N,1}\vee u_{l;\; 0,N-1} =0$.
\end{itemize}
Denote by $A$ the graded commutative algebra introduced in the above statement. Then there exists an isomorphism $\phi$ of graded algebras from $A$ to $\rmHH^*(\Lambda)$ sending $x_0$ to $x_{0,e}$, $\widetilde{x_l}$ to $u_{l;\; 0,N-1}$, $l=1,\dots,e$, $y$ to $y_{1,1}$,  $z$ to $x_{N,1}$, $w$ to $y_{(e-1)N+1,0}$ and $z'$ to $x_{eN,0}$.
It is easy to see that the set of relations at the assumption constitutes a Gr\"{o}bner basis with respect to the lexicographic order $w>y > z>x_{0 }>\widetilde{x_1}>\cdots >\widetilde{x_e}> z'$. 

\end{Proof}

\bigskip

\section{Gerstenhaber algebra structure}\label{Sect: Gerstenhaber}\

In this section, following Xu and Zhang \cite{XZ11}, we investigate the Gerstenhaber algebra structure on the Hochschild cohomology groups of a self-injective Nakayama algebra. We begin by establishing general formulae for truncated quiver algebras and then specialise to truncated basic cycle algebras.

\subsection{Gerstenhaber bracket for truncated quiver algebras}\

Let $\Lambda={\K Q}/J^N$ be a truncated quiver algebra. In Section~\ref{Section: preliminaries}, we recall that there exists a Gerstenhaber Lie bracket over the Hochschild cohomology complex $C^*(\Lambda, \Lambda)$. S\'{a}nchez-Flores defined the reduced bracket using reduced bar resolution $\bbB^*$ in \cite[Section~2]{San08}. Precisely, the bilinear map
\[
\bar{\circ}_i:\;\bbB^n\times \bbB^m\lra \bbB^{n+m-1}
\]
is defined by
\[
(f\bar\circ_i g)(1[{a_1}\mid  \cdots\mid  {a_{n+m-1}}]1):=f(1[{a_1}\mid  \cdots\mid a_{i-1} \mid 
g(1[{a_i}\mid  \cdots \mid  {a_{i+m-1}}]1)\mid  {a_{i+m}}\mid  \cdots \mid  {a_{n+m-1}}]1),
\]
for $f\in \bbB^n=\Hom_{\Lambda^e}(\Lambda\ot {\Lambda_+}^{\ot n}\ot \Lambda, \Lambda)$ and $g\in \bbB^m$, and define
\[
f\bar\circ g=\sum\limits_{i=1}^n(-1)^{(i-1)(m-1)}f\bar\circ_i g,\;\text{ and }\; [f , g ]_{\bbB}   =f \bar{\circ} g -(-1)^{(m-1)(n-1)} g \bar{\circ} f.
\]
S\'{a}nchez-Flores also showed that the \textit{reduced bracket} $[- , - ]_{\bbB} $ endows $s{\bbB}^{*}=\bigoplus_{n=1}^{\infty}{\bbB}^n$ with the structure of graded Lie algebra and the Gerstenhaber bracket and the reduced bracket provide the same graded Lie algebra structure
on $s\rmHH^{*}(\Lambda,\Lambda)$.

Through the comparison morphisms $\mu^*$ and $\omega^*$, Xu and Zhang introduce a bracket over $\bbP^*$ \cite{XZ11}, called the minimal bracket, as follows. Firstly, define a bilinear map
\[
\widetilde{\circ}_i:\;\bbP^n\times \bbP^m\lra \bbP^{n+m-1}
\]
to be 
\[
f \widetilde{\circ}_i g =\mu^*_{n+m-1}(\omega^*_n f \bar{\circ}_i\omega^*_m g )=(f \omega_n \bar{\circ}_i g \omega_m)\mu_{n+m-1},
\]
for $f\in \bbP^n=\Hom_{\Lambda^e}(\Lambda\ot \K Q_{\chi(n)} \ot \Lambda, \Lambda)$ and $g\in \bbP^m$, and define
\[
f \widetilde{\circ} g =\sum_{i=1}^n (-1)^{(m-1)(i-1)}f \widetilde{\circ}_i g , \;\text{ and }\;
  [f , g ]_{\bbP}   =f \widetilde{\circ} g -(-1)^{(n-1)(m-1)} g \widetilde{\circ} f.
\]

By Lemma~\ref{Lem: Parallel paths}, we know that $\bbP^n$ has a $\K$-basis $\{f_{(\gamma,b)} \di(\gamma,b) \in (Q_{\chi(n)}\pl  \calB)\}$, where $f_{(\gamma,b)}(1\ot \gamma'\ot 1)=\delta_{\gamma,\gamma'}b$, where $\delta_{\gamma,\gamma'}$ is the Kronecker symbol. For any $(\gamma_n,b_1)\in(Q_{\chi(n)}\pl  \calB)$ and $(\gamma_m,b_2)\in(Q_{\chi(m)}\pl  \calB)$, for short we write $(\gamma_n,b_1)\widetilde{\circ}_i(\gamma_m,b_2)$ instead of $f_{(\gamma_n,b_1)} \widetilde{\circ}_i f_{(\gamma_m,b_2)}$.

\bigskip

Let $\gamma, a,\eta\in Q$, if $\gamma=\gamma'a\gamma''$, we define $\gamma\vee^a\eta:=\gamma'\eta\gamma''$. Now we recall the definition of an overlap subpath.

\begin{Def}[{\cite[Section 4]{XZ11}}]
Let $\gamma\in Q_n$ and $b\in Q_r$ with $n,r\ge 1$. If $n>1$, the set of \textit{overlaps} of $b$ and $\gamma$ is 
\[
\calL(b,\gamma):=\{ a\in \calB\setminus Q_0\di b=b_a a,\;\gamma=a\gamma_a \},
\]
and if  $n=1$, $\gamma\in Q_1$, define
\[
\calL(b,\gamma):=\{\gamma_{b'}=\gamma \di b=b'\gamma b''\}.
\]
\end{Def}

\begin{Def}[{\cite[Definition 4.2]{XZ11}}]

Let $(\gamma_n,b_1)\in(Q_{\chi(n)}\pl Q_l)$ and $(\gamma_m,b_2)\in(Q_{\chi(m)}\ Q_r)$. Assume that  $\gamma_n=\alpha_{1,{\chi(n)}}$ with $\alpha_i\in Q_1$. For $s\in\mathds{Z}$ and $1\leq i\leq n$, if
$\alpha_{{\chi(i-1)+s+1},{\chi(i-1)+s+r}}=b_2$, then we define
\[
\gamma_n\vee^{b_2}_{i,s}\gamma_m:=\alpha_{1,{\chi(i-1)+s}}\gamma_m \alpha_{{\chi(i-1)+s+r+1},{\chi(n)}},
\]
which means that we replace the subpath $\alpha_{{\chi(i-1)+s+1},{\chi(i-1)+s+r}}$ of $\gamma_n$ with $\gamma_m$;
if
$\alpha_{{\chi(i-1)+s+1},{\chi(i-1)+s+r}}\neq b_2$, then 
\[
\gamma_n\vee^{b_2}_{i,s}\gamma_m:=0.
\]
\end{Def}

\begin{Prop}[{Compare with \cite[Theorem 4.3]{XZ11}}]\label{prop: circ product-i}

Let $(\gamma_n,b_1) \in (Q_{\chi(n)} \pl  \calB)$ and $(\gamma_m,b_2) \in (Q_{\chi(m)} \pl  \calB)$, denote $l( b_2) =r$, then for $1\leq i\leq n$, we have the following results:
\begin{itemize}
\item[(1)] Let $n=2h$ and $m=2t$, then
\[
(\gamma_n,b_1) \widetilde{\circ}_i (\gamma_m,b_2) =  
\begin{cases}
(\gamma_n\vee^{b_2}_{i,0}\gamma_m,b_1), & \text{if } r=N-1, \\
0,     & \text{otherwise}.\\
\end{cases}
\]

\item[(2)] Let $n=2h$ and $m=2t+1$, if $r\ge 1$, then
\begin{equation*}
(\gamma_n,b_1)\widetilde{\circ}_i(\gamma_m,b_2) =  
\left\{
\begin{aligned}
& \sum\limits_{s=N-r}^{N-1}\sum\limits_{\substack{ q }}\big((\gamma_n\vee^{b_2}_{i,s-1}\gamma_m)q,b_1q\big),      &      &\text{if $i< 2h$ is even}, \\
&\sum_{a\in\calL(\gamma_n,b_2)} \sum_{q'} ((\gamma_n \vee^a \gamma_m)q',(b_2\vee^a b_1)q'), & & \text{if } i=2h,\\
& \sum\limits_{s=0}^{N-r-1}\sum\limits_{\substack{ q }}\big((\gamma_n\vee^{b_2}_{i,s}\gamma_m)q,b_1q\big),   &      &\text{if } i \text{ is odd},
\end{aligned}
\right.
\end{equation*}
where $q$ (resp. $q'$) takes over paths of length $r-1$ (resp. $l(a) -1$) such that the sums lie in $\K(Q_{\chi(m+n-1)}\pl  \calB)$. If $r=0$, then $(\gamma_n,b_1)\widetilde{\circ}_i(\gamma_m,b_2)=0$.

\item[(3)] Let $n=2h+1$ and $m=2t$. When $i=2j$,
\begin{equation*}
(\gamma_n,b_1)\widetilde{\circ}_i(\gamma_m,b_2) =  
   \begin{cases}
    \sum\limits_{s=0}^{N-2}\sum\limits_{\substack{a_s, b_s}}\big(a_s(\gamma_n\vee^{b_2}_{i,s}\gamma_m)b_s,a_sb_1b_s\big),        &\text{if } r=N-1,\\
    0,        &\text{otherwise}, \\
   \end{cases}
\end{equation*}
where $a_s$ and $b_s$ take over paths of length $0\leq l( a_s) \leq N-2-s$, $l( a_s) +l( b_s) =N-2$ satisfying the sum lies in $\K(Q_{\chi(m+n-1)}\pl  \calB)$;
when $1<i=2j-1<2h+1$,
\begin{equation*}
(\gamma_n,b_1)\widetilde{\circ}_i(\gamma_m,b_2) = 
\begin{cases}
\sum\limits_{s=2-N}^{0} \sum\limits_{\substack{\widetilde{a}_s,\widetilde{b}_s}} \big( \widetilde{a}_s(\gamma_n \vee^{b_2}_{i,s} \gamma_m) \widetilde{b}_s, \widetilde{a}_s b_1 \widetilde{b}_s \big),    &\text{if }  r=N-1,\\
0,   &\text{otherwise}, \\
\end{cases}
\end{equation*}
where $\widetilde{a}_s$ and $\widetilde{b}_s$ take over paths of length $0\leq l( \widetilde{a}_s) \leq -s$, $l( \widetilde{a}_s) +l( \widetilde{b}_s) =N-2$ such that the sum lies in $\K(Q_{\chi(m+n-1)}\pl  \calB)$; 
when $i=1$,
\[
(\gamma_n,b_1)\widetilde{\circ}_1(\gamma_m,b_2) =
\sum\limits_{\substack{a\in\mathcal L(b_2,\gamma_n)}}\sum\limits_{\substack{b}}
\big((\gamma_n\vee^{a}\gamma_m)b,(b_2\vee^{a}b_1)b\big),   
\]
where $b$ takes over paths of length $l(a) -1$ satisfying the sum lies in $\K(Q_{\chi(m+n-1)}\pl  \calB)$;
when $i=2h+1$,
\[
(\gamma_n,b_1)\widetilde{\circ}_{2h+1} (\gamma_m,b_2)=
\begin{cases}
\sum\limits_{a\in\calL(\gamma_n,b_2)} \sum\limits_{b_a',b_a''} \big( b_a'(\gamma_n \vee^a \gamma_m)b_a'', b_a'(b_2\vee ^a b_1)b_a'' \big),    & \text{if } r=N-1, \\
0 ,   & \text{otherwise}, 
\end{cases}
\]
where $b_a'$ and $b_a''$ take over paths of length $0\le l( b_a') \le l(a) -1$, $l( b_a') +l( b_a'') =l(a) -1$ such that the sum lies in $\K(Q_{\chi(m+n-1)}\pl  \calB)$.

\item[(4)] Let $n=2h+1$ and $m=2t+1$, if $r\ge 1$, then
\begin{equation*}
(\gamma_n,b_1)\widetilde{\circ}_i(\gamma_m,b_2) =  
\left\{
\begin{aligned}
&\sum\limits_{s=0}^{N-r-1}  \sum\limits_{\substack{q}} \big((\gamma_n \vee^{b_2}_{i,s} \gamma_m)q, b_1q \big), & & \text{if $i$ is even},\\
&\sum\limits_{s=1-r}\limits^{0}\sum\limits_{\substack{q}}\big((\gamma_n\vee^{b_2}_{i,s}\gamma_m)q ,b_1q \big), & &\text{if $1< i<2h+1$ is odd},\\
&\sum\limits_{\substack{a\in\mathcal L(b_2,\gamma_n)}}\sum_b \big((\gamma_n \vee^a \gamma_m)b,(b_2 \vee^a b_1)b\big), & & \text{if }  i=1,\\
&\sum\limits_{a\in\calL(\gamma_n,b_2)}\sum_b \big((\gamma_n\vee^a \gamma_m)b,(b_2\vee^a b_1)b\big), & & \text{if } i=2h+1,
\end{aligned}
\right.
\end{equation*}
where $q$ and $b$ take over paths with length $l( q) =r-1$ and $l(b) =l(a) -1$ satisfying the sums lie in $\K(Q_{\chi(m+n-1)}\pl  \mathcal B)$.
If $r=0$, $(\gamma_n,b_1)\widetilde{\circ}_i(\gamma_m,b_2)=0$.
\end{itemize}

\end{Prop}

\begin{Proof}

We restrict our attention to Case (2) with $i$ even, Case (3) with $i=2h+1$, the remaining cases have already been addressed in \cite{XZ11} or can be obtained similarly.

In Case (2) and Case (3), we have that $n+m-1=2h+2t$ is even, then
\begin{align}\label{equation:case(23)}
& f_{(\gamma_n,b_1)}\widetilde{\circ}_if_{(\gamma_m,b_2)}(1\ot \alpha_{1,(h+t)N}\ot 1)\nonumber\\
&= (f_{(\gamma_n,b_1)}\omega_n \bar{\circ}_i f_{(\gamma_m,b_2)} \omega_m) \mu_{n+m-1} (1\ot \alpha_{1,(h+t)N}\ot 1) \nonumber\\
&= \sum\limits_{\substack{ x_1,\dots,x_{h+t}}}(f_{(\gamma_n,b_1)}\omega_n\bar{\circ}_i f_{(\gamma_m,b_2)}\omega_m) (1[p_1\mid \cdots\mid p_{2(h+t)}]q) \\
&= \sum\limits_{\substack{x_1,\dots,x_{h+t}}} f_{(\gamma_n,b_1)}\omega_n (1[p_1\mid \cdots\mid p_{i-1}\mid f_{(\gamma_m,b_2)}\omega_m(1[p_i\mid \cdots\mid p_{i+m-1}]1) \mid p_{i+m}\mid \cdots\mid p_{2(h+t)}]q)\nonumber,
\end{align}	
where $p_1\cdots p_{2(h+t)}\beta=\alpha_{1,(h+t)N}$ and  for  $k=1,\dots,h+t$, $1\le x_k\le N-1$, and
\[
l( p_a)  =
\begin{cases}
	1,    &\text{if } a=2 k,\\
	x_k,      &\text{if } a=2k-1.
\end{cases}  
\]

\textbf{Case (2) with $i$ even}  

First, we revise the proof of Case (2) in \cite[Theorem 4.3]{XZ11} for the situation $i=2j$, where $j\ge 1$.

Let $n=2h$ and $m=2t+1$, by definition of $\omega_*$, the non-zero terms of Equation~\ref{equation:case(23)} satisfy the following conditions: $x_1=\cdots=x_{j-1}=x_{j+1}=\cdots=x_{h+t}=N-1$, $r\ge 1$ and $x_j+r\ge N$. Moreover, one should have  
\[
\gamma_m= \alpha_{{(j-1)N+x_j+1},{(j+t-1)N+x_j+1}}.
\]

If $i\ne 2h$, the Equation~\ref{equation:case(23)} can be written as
\[
\sum\limits_{x_j=N-r}^{N-1} f_{(\gamma_n,b_1)}(1\ot \underbrace{\alpha_{1,{(j-1)N+x_j}}b_2 \alpha_{{(j+t-1)N+x_j+2},{(h+t)N-r+1}}}_{hN} \ot q)=\sum\limits_{x_j=N-r}^{N-1} b_1 q,
\]
with the condition that $\gamma_n=\alpha_{1,{(j-1)N+x_j}}b_2 \alpha_{{(j+t-1)N+x_j+2},{(h+t)N-r+1}} $, i.e., $\alpha_{1,{(h+t)N}}=(\gamma_n\vee^{b_2}_{i,x_j-1}\gamma_m)q$, where $q=\alpha_{{(h+t)N-r+2},{(h+t)N}}$.
Therefore 
\[
{(\gamma_n,b_1)}\widetilde{\circ}_i{(\gamma_m,b_2)}=\sum\limits_{s=N-r}^{N-1}\sum\limits_{q}((\gamma_n\vee^{b_2}_{i,s-1}\gamma_m)q,b_1q),
\]
where $q$ takes over paths of length $r-1$ such that the sum  lies in $\K(Q_{\chi(m+n-1)}\pl  \calB)$.

If $i=2h$, let $b_2=\beta_1\cdots\beta_r$, the Equation~\ref{equation:case(23)} can be written as
\[
\begin{aligned}
& \sum\limits_{x_h=N-r}^{N-1} f_{(\gamma_n,b_1)} \omega_{2h}( 1[p_1\mid  \cdots\mid  p_{2h-1}\mid  b_2 ]q)\\
& = \sum\limits_{x_h=N-r}^{N-1} f_{(\gamma_n,b_1)}(1\ot \underbrace{\alpha_{1,{(h-1)N+x_h}}\beta_{1,N-x_h}}_{hN} \ot \beta_{N-x_h+1,r}\,q)\\
& = \sum\limits_{x_h=N-r}^{N-1} b_1  \beta_{N-x_h+1,r}\,q,
\end{aligned}
\]
if $\gamma_m=\alpha_{{(h-1)N+x_h+1}, {(h+t-1)N+x_h+1}}$ and $\gamma_n=\alpha_{1,{(h-1)N+x_h}}\beta_{1,N-x_h}$, which means that for $a=\beta_{1,N-x_h}\in \calL(\gamma_n,b_2)$, $\alpha_{1,(h+t)N}=(\gamma_n \vee^a \gamma_m)q$ and $ b_1  \beta_{N-x_h,r}q=(b_2 \vee^a b_1)q $. 
Therefore
\[
(\gamma_n,b_1)\widetilde{\circ}_{2h} (\gamma_m,b_2)= \sum_{a\in\calL(\gamma_n,b_2)} \sum_{q} ((\gamma_n \vee^a \gamma_m)q,(b_2\vee^a b_1)q),
\]
where $q$ takes over paths of length $l(a) -1$ such that the sum lies in $\K(Q_{\chi(m+n-1)}\pl  \calB)$.

\bigskip

\textbf{Case (3) with $i=2h+1$}

For $n=2h+1$, $m=2t$, and $i=2h+1$, the definition of $\omega_*$ implies the non-zero terms in Equation~\eqref{equation:case(23)} satisfy 
\[
x_2=\cdots=x_{h+t}=N-1 \text{ and } r=N-1.
\]
Assume $b_2=\beta_1\cdots\beta_{N-1}$ If $\gamma_m=\alpha_{{(h-1)N+x_1+2},{(h+t-1)N+x_1+1}}$, the Equation~\eqref{equation:case(23)} can be written as 
\begin{equation*}
   \begin{aligned}
    &\sum_{x_1=1}^{N-1} f_{(\gamma_n,b_1)}\omega_n(1[p_1\mid \cdots\mid p_{2h}\mid b_2]q')\\
    =&\sum_{x_1=1}^{N-1}\sum_{p=1}^{x_1} f_{(\gamma_n,b_1)}(\alpha_{1,p-1}\ot \underbrace{\alpha_{p,{(h-1)N+x_1+1}}\beta_{1,N-x_1+p-1}}_{hN+1}\ot \beta_{N-x_1+p,N-1} \alpha_{(h+t-1)N+x_1+2,(h+t)N}).
    \end{aligned}
\end{equation*}
If $\gamma_n=\alpha_{p,{(h-1)N+x_1+1}}\beta_{1,N-x_1+p-1}$, the result is 
\[
\sum_{x_1=1}^{N-1}\sum_{p=1}^{x_1} \alpha_{1,p-1}b_1 \beta_{N-x_1+p,N-1} \alpha_{(h+t-1)N+x_1+2,(h+t)N}.
\]
Let $a=\beta_{1,N-x_1+p-1}\in \calL(\gamma_n,b_2)$, $b_a'=\alpha_{1,p-1}$ and $b_a''= \alpha_{(h+t-1)N+x_1+2,(h+t)N}$, then we can write 
\[
\alpha_{1,(h+t)N}=b_a'(\gamma_n \vee^a \gamma_m)b_a'' \text{ and }\alpha_{1,p-1}b_1 \beta_{N-x_1+p,N-1} \alpha_{(h+t-1)N+x_1+2,(h+t)N}=b_a'(b_2\vee ^a b_1)b_a''.
\]
Therefore, 
\[
(\gamma_n,b_1)\widetilde{\circ}_{2h+1} (\gamma_m,b_2)=\sum_{a\in\calL(\gamma_n,b_2)} \sum_{b_a',b_a''} ( b_a'(\gamma_n \vee^a \gamma_m)b_a'', b_a'(b_2\vee ^a b_1)b_a'' ),
\]
where for fixed $a\in\calL(\gamma_n,b_2)$, $b_a'$, $b_a''$ satisfy $0\le l( b_a') \le l(a) -1$ and $l( b_a') +l( b_a'') =l(a)-1$.

\end{Proof}

\begin{Rem}

The above result differs from \cite[Theorem 4.3]{XZ11} in the following cases: Case (2) where $n$ is even, $m$ is odd and $i$ even; Case (3) where $n$ is odd, $m$ is even and $i=n$; and Case (4) where $n, m$ are odd and $i>1$ is also odd.

\end{Rem}

\begin{Rem}

Let $f\in \bbP^m$ and $g\in \bbP^n$, then we have
\begin{equation*}
	\begin{aligned}
	d^*[f,g]_{\bbP}
    &=d^* \mu^*_{m+n-1}[\omega^*_mf,\omega^*_ng]_{\bbB}=\mu^*_{m+n}{b^*} [\omega^*_mf, \omega^*_ng] _{\bbB}\\
	&=\mu^*_{m+n}[{b^*}\omega^*_mf,\omega^*_ng]_{\bbB}+(-1)^{m-1}\mu^*_{m+n}[\omega^*_mf,{b^*} \omega^*_ng]_{\bbB}\\
	&=\mu^*_{m+n}[\omega^*_{m+1} d^*f,\omega^*_ng]_{\bbB}+(-1)^{m-1}\mu^*_{m+n}[\omega^*_mf,\omega^*_{n+1} d^*g]_{\bbB}\\
	&=[d^*f,g]_{\bbP}+(-1)^{m-1}[f,d^*g]_{\bbP}.
	\end{aligned}
\end{equation*}
This shows the minimal bracket induce a well-defined bracket in Hochschild cohomology
\[
[-,-]_{\bbP}:\; \rmHH^m(\Lambda)\times \rmHH^n(\Lambda) \to \rmHH^{m+n-1}(\Lambda).
\]
This endows $s\rmHH^*(A)$ with a Lie algebra structure, which is isomorphic to the Lie algebra structure induced by the reduced bracket (or the Gerstenhaber bracket).

 However, on the chain level, $s\bbP^*$ endowed with the differential $d^*$ and the minimal bracket $[-,-]_\bbP$ is NOT a differential graded Lie algebra as is shown in Remark~\ref{Remark: example for P not be a dg lie algebra}, since it is not even true for truncated basic cycle algebras.
It would be interesting to determine the $L_\infty$-structure on $s\bbP^*$.
\end{Rem}

\bigskip

\subsection{Gerstenhaber brackets for truncated basic cycle algebras}\label{Subsection: Gerstenhaber bracket for truncated basic cycle algebras}\

Recall that a truncated basic cycle algebra $\Lambda=\K Z_e/J^N$ with $N,\, e\ge 2$ is given by the quiver $Z_e$ with vertex set $\{1,2,\dots,e\}$ and arrow set $\{\alpha_1,\dots,\alpha_e\}$ satisfying $o(\alpha_i)\equiv t(\alpha_{\ul{i-1}})$, $i=1,\dots,e$. 

\medskip

Since a truncated basic cycle algebra is a special instance of truncated quiver algebra, we may invoke Proposition~\ref{prop: circ product-i} to determine the operation $\widetilde{\circ}$ between two elements. We begin by considering the case $\gcd(N,e)=1$.

\begin{Prop}[{Compare with \cite[Proposition 5.1]{XZ11}}]\label{prop: circ product}

Assume that $\gcd(N,e)=1$. Let $n\ge 1$, $m\ge 0$, $1\le l,\, p\le e$, and $0\le j,\, r\le N-1$,  with  $\chi(n)\equiv j\mode$ and $\chi(m)\equiv r\mode$. For any integer $s$, denote by $x_s\in[1,e]$ and $x_s'\in[2,e+1]$ the unique solution of the congruence $(x-1)N+l+s\equiv p\mode$.

Then, for $(\alpha_l^{\chi(n)},\alpha_l^j) \in (Q_{\chi(n)}\pl  \calB)$ and $(\alpha_p^{\chi(m)},\alpha_p^r) \in (Q_{\chi(m)}\pl  \calB)$, we have the following statements:

\begin{itemize}

\item[(1)] Let $n=2h$, $m=2t$. Then we have
\[
(\alpha_l^{hN},\alpha_l^j)\widetilde{\circ}(\alpha_p^{tN},\alpha_p^r) = \delta_{r,N-1} \Big(\Big[\frac{h-x_0}{e}\Big]-\Big[\frac{h-x_1}{e}\Big]\Big) (\alpha_l^{(h+t-1)N+1},\alpha_l^j).
\]

\item[(2)] Let  $n=2h$ and  $m=2t+1$. Then we have
\begin{equation*}
(\alpha_l^{hN},\alpha_l^j)\widetilde{\circ}(\alpha_p^{tN+1},\alpha_p^r)=
\begin{cases}  
\Big(\sum\limits_{s=0}^{N-1} \Big[\frac{h+e-x_s}{e}\Big]\Big)(\alpha_l^{(h+t)N},\alpha_l^{j+r-1}),&\text{if }r\ge 1,\\
0,&\text{if }r=0.
\end{cases}
\end{equation*}

\item[(3)] Let $n=2h+1$ and $m=2t$. Then we have
\begin{equation*}
\begin{aligned}
(\alpha_l^{hN+1},\alpha_l^j)\widetilde{\circ}(\alpha_p^{tN},\alpha_p^r)
&=\Big[\frac{r-1+e-\overline{l-p}}{e}\Big](\alpha_p^{(h+t)N},\alpha_p^{j+r-1})\\
&\quad+\delta_{r,N-1}\sum_{s=0}^{N-2}\sum_{i=0}^{N-2-s}\Big(\Big[\frac{h-x_{s+2}}{e}\Big]-\Big[\frac{h-x_{s+1}}{e}\Big]\Big)(\alpha_{\ul{l-i}}^{(h+t)N},\alpha_{\ul{l-i}}^{j+N-2}).
\end{aligned}
\end{equation*}

\item[(4)] Let $n=2h+1$ and $m=2t+1$. If $r\ge 1$, then we have
\begin{equation*}
\begin{aligned}
(\alpha_l^{hN+1},\alpha_l^j)\widetilde{\circ}(\alpha_p^{tN+1},\alpha_p^r)
&=\Big[\frac{r-1+e-\overline{l-p}}{e}\Big](\alpha_p^{(h+t)N+1},\alpha_p^{j+r-1})\\
&\quad+\Big(\sum\limits_{s=-r}^{-1}\Big[\frac{h+1-x_{s+1}'}{e}\Big]+\sum_{s=0}^{N-r-1}\Big[\frac{h-x_{s+1}}{e}\Big]+N\Big)(\alpha_l^{(h+t)N+1},\alpha_l^{j+r-1}).
\end{aligned}
\end{equation*}
If $r=0$, $(\alpha_l^{hN+1},\alpha_l^j)\widetilde{\circ}(\alpha_p^{tN+1},\alpha_p^r)=0$.
\end{itemize}
\end{Prop}

\begin{Proof}

Here we only provide the proofs of the cases (4), since the proofs of (1) and (2) are entirely analogous.

Let $n=2h+1$ and $m=2t+1$ be both odd.
By Proposition~\ref{prop: circ product-i} (4), if $r=0$, 
\[
(\alpha_l^{hN+1},\alpha_l^j)\widetilde{\circ}(\alpha_p^{tN+1},\alpha_p^r)=0.
\]
If $r\ge 1$, for $1\le i\le 2h+1$, we have
\begin{equation*}
(\alpha_l^{hN+1},\alpha_l^j)\widetilde{\circ}_i(\alpha_p^{tN+1},\alpha_p^r)=
\begin{cases}
\sum\limits_{s=0}^{N-r-1}\big(\alpha_l^{hN+r}\vee^{\alpha_p^r}_{i,s}\alpha_p^{tN+1},\alpha_l^{j+r-1}\big),   & \text{if } i=2k,\\
\sum\limits_{s=1-r}^{0}\big(\alpha_l^{hN+r}\vee^{\alpha_p^r}_{i,s}\alpha_p^{tN+1},\alpha_l^{j+r-1}\big), & \text{if } 1<i=2k-1<2h+1, \\
\sum\limits_{a\in\calL(\alpha_p^r,\alpha_l^{hN+1})} \big(\alpha_p^{(h+t)N+1},\alpha_p^{j+r-1} \big), &\text{if } i=1,\\
\sum\limits_{a\in\calL(\alpha_l^{hN+1},\alpha_p^r)} \big( \alpha_l^{(h+t)N+1},\alpha_l^{j+r-1} \big), & \text{if } i=2h+1.
\end{cases}
\end{equation*}
Moreover,
\begin{equation*}
\alpha_l^{hN+r}\vee^{\alpha_p^r}_{i,s}\alpha_p^{tN+1}=  
\left\{
\begin{array}{ll}
\alpha_l^{(h+t)N+1},&\text{if $i=2k$ and $(k-1)N+l+s+1\equiv p \mode$},  \\
\alpha_l^{(h+t)N+1},&\text{if $1<i=2k-1<2h+1$ and $(k-1)N+l+s\equiv p \mode$}, \\
0, &\text{otherwise}.
\end{array}
\right.
\end{equation*}
For $a\in\mathcal L(\alpha_p^r,\alpha_l^{hN+1})$, it should satisfy $p+r-l(a) \equiv l\mode$. Then the cardinality of the set $\calL(\alpha_p^r,\alpha_l^{hN+1})$ is $\big[\frac{r-1+e-\ol{l-p}}{e}\big]$, which implies
\[
(\alpha_l^{hN+1},\alpha_l^j)\widetilde{\circ}_1(\alpha_p^{tN+1},\alpha_p^r)=\big[\frac{r-1+e-\ol{l-p}}{e}\big] (\alpha_p^{(h+t)N+1},\alpha_p^{j+r-1}).
\]
If $i=2h+1>1$, $a\in\calL(\alpha_l^{hN+1},\alpha_p^r)$ satisfies $h N+l-l(a) +1\equiv p\mode$ and $1\le l(a) \le r$, let $s=-l(a) +1$, then we have 
\[
(\alpha_l^{hN+1},\alpha_l^j)\widetilde{\circ}_{2h+1}(\alpha_p^{tN+1},\alpha_p^r)=\sum\limits_{s=1-r}^{0}\big(\alpha_l^{hN+r}\vee^{\alpha_p^r}_{2h+1,s}\alpha_p^{tN+1},\alpha_l^{j+r-1}\big).
\]

Finally we obtain
\begin{equation*}
\begin{aligned}
&(\alpha_l^{hN+1},\alpha_l^j)\widetilde{\circ}(\alpha_p^{tN+1},\alpha_p^r)\\
& = \sum_{k=2}^{h+1}(\alpha_l^{hN+1},\alpha_l^j)\widetilde{\circ}_{2k-1}(\alpha_p^{tN+1},\alpha_p^r)
+\sum_{k=1}^h(\alpha_l^{hN+1},\alpha_l^j)\widetilde{\circ}_{2k}(\alpha_p^{tN+1},\alpha_p^r) +(\alpha_l^{hN+1},\alpha_l^j)\widetilde{\circ}_1(\alpha_p^{tN+1},\alpha_p^r)\\
&=\sum\limits_{s=1-r}^{0} c_1(\alpha_l^{(h+t)N+1},\alpha_l^{j+r-1})+
\sum\limits_{s=0}^{N-r-1}c_2(\alpha_l^{(h+t)N+1},\alpha_l^{j+r-1}) + \Big[\frac{r-1+e-\overline{l-p}}{e}\Big](\alpha_p^{(h+t)N+1},\alpha_p^{j+r-1})\\
&= \Big(\sum\limits_{s=1-r}^{0}\Big[\frac{h+1-x_{s}'}{e}\Big]+
    \sum_{s=0}^{N-r-1}\Big[\frac{h-x_{s+1}}{e}\Big]+N\Big)(\alpha_l^{(h+t)N+1},\alpha_l^{j+r-1}) \\
    &\qquad \qquad \quad \quad\qquad \qquad\qquad \qquad +\Big[\frac{r-1+e-\overline{l-p}}{e}\Big](\alpha_p^{(h+t)N+1},\alpha_p^{j+r-1}),
    \end{aligned}
\end{equation*}
where $c_1=\big[ \frac{h+1-x_{s}'}{e} \big]+1$ is the number of solutions $x$ satisfying $(x-1)N+l+s\equiv p \mode$ and $2\le x\le h+1$, $c_2=\big[ \frac{h-x_{s+1}}{e}\big]+1$ is the number of solutions $x$  satisfying $(x-1)N+l+s+1\equiv p \mode$ and $1\leq x\leq h$.
\end{Proof}

Proposition ~\ref{prop: circ product} yields the following  results.

\begin{Prop}\label{Prop: gerstenhaber bracket1}
Assume that $\gcd(N,e)= 1$. Let $n,\, m\ge 0$, $1\le l,p\le e$, and $0\le j,\, r\le N-1$. For any integer $s$, denote by $x_s\in[1,e]$ and $x_s'\in[2,e+1]$ (resp. $y_s\in[1,e]$ and $y_s'\in[2,e+1]$) the unique solution of the congruence $(x-1)N+l+s\equiv p\mode$ (resp. $(y-1)N+p+s\equiv l\mode$).

Assume that $\chi(n)\equiv j\mode$, $\chi(m)\equiv r\mode$.
Then for $(\alpha_l^{\chi(n)},\alpha_l^j)\in(Q_{\chi(n)}\pl  \mathcal B)$, $(\alpha_p^{\chi(m)},\alpha_p^r)\in(Q_{\chi(m)}\pl  \mathcal B)$, we have the following results:
\begin{itemize}

\item[(1)] Let $n=2h$ and $m=2t$, then we have
	\begin{equation*}
		\begin{aligned}
		\big[(\alpha_l^{hN},\alpha_l^j),(\alpha_p^{tN},\alpha_p^r)\big]
		=&\delta_{r,N-1}\Big(\Big[\frac{h-x_0}{e}\Big]-\Big[\frac{h-x_1}{e}\Big]\Big)(\alpha_l^{(h+t-1)N+1},\alpha_l^j)\\
		&+\delta_{j,N-1}\Big(\Big[\frac{t-y_0}{e}\Big]-\Big[\frac{t-y_1}{e}\Big]\Big)(\alpha_p^{(h+t-1)N+1},\alpha_p^r).
		\end{aligned}
	\end{equation*}

\item[(2)] Let $n=2h$ and $m=2t+1$, then we have
	\begin{equation*}
		\begin{aligned}
		\big[(\alpha_l^{hN},\alpha_l^j),(\alpha_p^{tN+1},\alpha_p^r)\big]
		=&\Big(\sum_{s=0}^{N-1} (1-\delta_{r,0})\Big[\frac{h+e-x_s}{e}\Big]-\Big[\frac{j-1+e-\overline{p-l}}{e}\Big]\Big)(\alpha_l^{(h+t)N},\alpha_l^{j+r-1})\\
		&-\delta_{j,N-1}\sum_{s=0}^{N-2}\sum_{i=0}^{N-2-s}\Big(\Big[\frac{t-y_{s+2}}{e}\Big]-\Big[\frac{t-y_{s+1}}{e}\Big]\Big)(\alpha_{\ul{p-i}}^{(h+t)N},\alpha_{\ul{p-i}}^{r+N-2}).
		\end{aligned}
	\end{equation*}
    In particular, if $n=0$, 
    \begin{equation*}
	  \begin{aligned}
		\big[(e_l,\alpha_l^j),(\alpha_p^{tN+1},\alpha_p^r)\big]
		=&-\Big[\frac{j-1+e-\overline{p-l}}{e}\Big](\alpha_l^{tN},\alpha_l^{j+r-1})\\
		&-\delta_{j,N-1}\sum_{s=0}^{N-2}\sum_{i=0}^{N-2-s}\Big(\Big[\frac{t-y_{s+2}}{e}\Big]-\Big[\frac{t-y_{s+1}}{e}\Big]\Big)(\alpha_{\ul{p-i}}^{tN},\alpha_{\ul{p-i}}^{r+N-2}).
		\end{aligned}
	\end{equation*}

\item[(3)] Let $n=2h+1$ and $m=2t+1$, then $\big[(\alpha_l^{hN+1},\alpha_l^j),(\alpha_p^{tN+1},\alpha_p^r)\big]$ is given by
	\begin{equation*}
		\begin{aligned}
		&\left((1-\delta_{r,0}) \Big(\sum\limits_{s=1-r}^{0}\Big[\frac{h+1-x_{s}'}{e}\Big]+\sum_{s=0}^{N-r-1}\Big[\frac{h-x_{s+1}}{e}\Big]+N \Big) -\Big[\frac{j-1+e-\overline{p-l}}{e}\Big]\right)(\alpha_l^{(h+t)N+1},\alpha_l^{j+r-1})\\
		&-\left((1-\delta_{j,0}) \Big(\sum\limits_{s=1-j}^{0} \Big[\frac{t+1-y_{s}'}{e}\Big]+\sum_{s=0}^{N-j-1}\Big[\frac{t-y_{s+1}}{e}\Big]+N \Big)-\Big[\frac{r-1+e-\overline{l-p}}{e}\Big]\right)(\alpha_p^{(h+t)N+1},\alpha_p^{j+r-1}).
		\end{aligned}
	\end{equation*}
    
\end{itemize}

\end{Prop}

Now consider the case  $\gcd(N,e)=d>1$. Write $N=dN_1$ and $e=de_1$. Denote by 
$\kappa_{u,v}=
\begin{cases}
       1,   &\text{if } u\di v, \\
       0,   &\text{otherwise}.
\end{cases}$

\begin{Prop}\label{prop: Circle product2} 
Let $\Lambda=\K Z_e/J^N$ be the truncated  basic cycle algebra with $N,e\ge 2$ and $\gcd(N,e)= 1$.
Assume that $n\ge1$, $m\ge 0$, $1\le l,p\le e$, and $0\le j,\, r\le N-1$ with  $\chi(n)\equiv j\mode$, $\chi(m)\equiv r\mode$. For an integer $s$, let $x_s\in[1,e_1]$ and $x_s'\in [2,e_1+1]$ be the unique solution of the congruence $(x-1)N_1\equiv \frac{p-s-l}{d} \model$ if $d\di p-s-l$.

Then for $(\alpha_l^{\chi(n)},\alpha_l^j)\in(Q_{\chi(n)}\pl  \mathcal B)$, $(\alpha_p^{\chi(m)},\alpha_p^r)\in(Q_{\chi(m)}\pl  \mathcal B)$, the following statements hold:

\begin{itemize}

\item[(1)] Let $n=2h$, $m=2t$, then we have
\[
(\alpha_l^{hN},\alpha_l^j)\widetilde{\circ}(\alpha_p^{tN},\alpha_p^r)=0.
\]

\item[(2)] If $n=2h$, $m=2t+1$, then we have
 \[
 (\alpha_l^{hN},\alpha_l^j)\widetilde{\circ}(\alpha_p^{tN+1},\alpha_p^r)=\Big(\sum\limits_{s=0}^{N -1}\kappa_{d,p-s-l}\Big[\frac{h+e_1-x_s}{e_1}\Big]\Big)(\alpha_l^{(h+t)N},\alpha_l^{j+r-1}).
 \]

\item[(3)] If $n=2h+1$, $m=2t$, then we have
\begin{equation*}
     (\alpha_l^{hN+1},\alpha_l^j)\widetilde{\circ}(\alpha_p^{tN},\alpha_p^r)
    =  \Big[\frac{r-1+e-\overline{l-p}}{e}\Big](\alpha_p^{(h+t)N},\alpha_p^{j+r-1}).
  \end{equation*}

\item[(4)] If $n=2h+1$, $m=2t+1$, then we have
\begin{equation*}
\begin{aligned}
& (\alpha_l^{hN+1},\alpha_l^j)\widetilde{\circ}(\alpha_p^{tN+1},\alpha_p^r)
=\Big[\frac{r-1+e-\overline{l-p}}{e}\Big](\alpha_p^{(h+t)N+1},\alpha_p^{j+r-1})\\
&\quad\quad\quad\quad+\Big( \sum_{s=1-r}^{ 0}\kappa_{d,p-l-s}\Big[\frac{h+1+e_1-x_s'}{e_1}\Big] +\sum_{s=0}^{N-r-1}\kappa_{d,p-l-s-1}\Big[\frac{h+e_1-x_{s+1}}{e_1}\Big]\Big)(\alpha_l^{(h+t)N+1},\alpha_l^{j+r-1}).
\end{aligned}
\end{equation*}
  
\end{itemize}

\end{Prop}

\begin{Proof}

We only consider the cases (3) and (4).

For (3), since $\gcd(e,N)=d>1$, then for any integer $b$, the congruence $Nx\equiv b\mode$ has solution in $\mathbb{Z}$ if and only if $d\di b$. The fact that $(\alpha_p^{\chi(m)},\alpha_p^r)\in(Q_{\chi(m)}\pl  \calB)$ implies that $tN\equiv r\mode$. Thus, $r\neq N-1$, because otherwise $tN\equiv N-1\mode$ would force $d=1$, which contradicts our assumption. By Proposition~\ref{prop: circ product-i} (3), since $r\neq N-1$, we obtain
\begin{equation*}
\begin{aligned}
(\alpha_l^{hN+1},\alpha_l^j)\widetilde{\circ}(\alpha_p^{tN},\alpha_p^r)
&=(\alpha_l^{hN+1},\alpha_l^j)\widetilde{\circ}_1(\alpha_p^{tN},\alpha_p^r)\\
&=\sum\limits_{\substack{a\in\mathcal L(\alpha_p^r, \alpha_l^{hN+1})}}
\big(\alpha_l^{hN+l(a) }\vee^{a} \alpha_p^{tN}, \alpha_p^{r}\vee^{a} \alpha_l^{j+l(a) -1}\big)\\
&= \Big[\frac{r-1+e-\overline{l-p}}{e}\Big](\alpha_p^{(h+t)N},\alpha_p^{j+r-1}).
\end{aligned}
\end{equation*}

Next, we discuss Case (4).
By a proof similar to that of Proposition \ref{prop: circ product} (4), we can obtain
\begin{equation*}
\begin{aligned}
&(\alpha_l^{hN+1},\alpha_l^j)\widetilde{\circ}(\alpha_p^{tN+1},\alpha_p^r)\\
= &\sum_{k=2}^{h+1}(\alpha_l^{hN+1},\alpha_l^j)\widetilde{\circ}_{2k-1}(\alpha_p^{tN+1},\alpha_p^r)
+\sum_{k=1}^h(\alpha_l^{hN+1},\alpha_l^j)\widetilde{\circ}_{2k}(\alpha_p^{tN+1},\alpha_p^r)+(\alpha_l^{hN+1},\alpha_l^j)\widetilde{\circ}_1(\alpha_p^{tN+1},\alpha_p^r) \\
= &\sum_{k=2}^{h+1} \sum\limits_{s=1-r}^{ 0} \alpha_l^{hN+r}\vee^{\alpha_p^r}_{2k-1,s} \alpha_p^{tN+1}+
\sum_{k=1}^h   \sum\limits_{s=0}^{ N-r-1} \alpha_l^{hN+r}\vee^{\alpha_p^r}_{2k,s}\alpha_p^{tN+1}+\sum\limits_{\substack{a\in\mathcal L(\alpha_p^r,\alpha_l^{hN+1})}}
(\alpha_p^{(h+t)N+1},\alpha_p^{j+r-1})   \\
= & \sum\limits_{s=1-r}^{0} c_{1,s}(\alpha_l^{(h+t)N+1},\alpha_l^{j+r-1})+\sum\limits_{s=0}^{N-r-1}c_{2,s}(\alpha_l^{(h+t)N+1},\alpha_l^{j+r-1})  +c_3(\alpha_p^{(h+t)N+1},\alpha_p^{j+r-1}),
    \end{aligned}
\end{equation*}
where $c_{1, s}$ is the number of solutions satisfying $(x-1)N+l +s\equiv p\mode$ and $2\leq x\leq h+1$, $c_{2,s}$ is the number of solutions satisfying $(x-1)N+l+s+1\equiv p\mode$ and $1\leq x\leq h$, $c_3$ is the number of solutions satisfying
$l\equiv p+x\mode$ and $x\in[0,r-1]$.
If $d\nmid p-l-s$, the congruence $(x-1)N+l+s+1\equiv p\mode$ has no solution, then $c_{1,s}=0$. If $d\di p-l-s$,  $c_{1,s}=\Big[\frac{h+1+e_1-x_s'}{e_1}\Big]$. If $d\di p-l-s-1$, $c_{2,s}=\Big[\frac{h+e_1-x_{s+1}}{e_1}\Big]$. And it is clear that we have $c_3=\Big[\frac{r-1+e-\overline{l-p}}{e}\Big]$.
Hence, we prove the equation
\begin{equation*}
\begin{aligned}
& (\alpha_l^{hN+1},\alpha_l^j)\widetilde{\circ}(\alpha_p^{tN+1},\alpha_p^r)
=\Big[\frac{r-1+e-\overline{l-p}}{e}\Big](\alpha_p^{(h+t)N+1},\alpha_p^{j+r-1})\\
&\quad\quad\quad\quad+\Big( \sum_{s=1-r}^{ 0}\kappa_{d,p-l-s}\Big[\frac{h+1+e_1-x_s'}{e_1}\Big] +\sum_{s=0}^{N-r-1}\kappa_{d,p-l-s-1}\Big[\frac{h+e_1-x_{s+1}}{e_1}\Big]\Big)(\alpha_l^{(h+t)N+1},\alpha_l^{j+r-1}).
\end{aligned}
\end{equation*}

\end{Proof}

Following Proposition ~\ref{prop: Circle product2}, we obtain the following result.

\begin{Prop}\label{Prop:gerstenhaber bracket2}

Let $\Lambda=\K Z_e/J^N$ be the truncated  basic cycle algebra with $N,\,e\ge 2$ and $\gcd(N,e)=d> 1$ with $N=dN_1$ and $e=de_1$. Assume that $n,\, m\ge 0$, $1\le l,\, p\le e$, and $0\le j,\,r\le N-1$ with $\chi(n)\equiv j\mode$, $\chi(m)\equiv r\mode$. For an integer $s$, let $x_s\in[1,e_1]$ and $x_s'\in [2,e_1+1]$ (resp. $y_s\in[1,e_1]$ and $y_s'\in [2,e_1+1]$) be the unique solution of the congruence $(x-1)N_1\equiv \frac{p-s-l}{d} \model$ (resp. $(y-1)N_1\equiv \frac{l-s-p}{d} \model$) if $d\di p-s-l$ (resp. $d\di l-s-p$).

Then for $(\alpha_l^{\chi(n)},\alpha_l^j)\in(Q_{\chi(n)}\pl  \mathcal B)$, $(\alpha_p^{\chi(m)},\alpha_p^r)\in(Q_{\chi(m)}\pl  \mathcal B)$, we have the following results:

\begin{itemize}

\item[(i)] If $n=2h$, $m=2t$, then we have
\[
\big[(\alpha_l^{hN},\alpha_l^j),(\alpha_p^{tN},\alpha_p^r)\big]=0.
\]

\item[(ii)] If $n=2h$ and $m=2t+1$, then we have
\[
\big[(\alpha_l^{hN},\alpha_l^j),(\alpha_p^{tN+1},\alpha_p^r)\big]
=\Big(\sum_{s=0}^{N -1}\kappa_{d,p-s-l}\Big[\frac{h+e_1-x_s}{e_1}\Big]-\Big[\frac{j-1+e-\overline{p-l}}{e}\Big]\Big)(\alpha_l^{(h+t)N},\alpha_l^{j+r-1}).
\]
In particular, if $n=0$, 
\[
\big[(e_l,\alpha_l^j),(\alpha_p^{tN+1},\alpha_p^r)\big]
= -\Big[\frac{j-1+e-\overline{p-l}}{e}\Big](\alpha_l^{tN},\alpha_l^{j+r-1}).
\]

\item[(iii)] If $n=2h+1$ and $m=2t+1$, then $\big[(\alpha_l^{hN+1},\alpha_l^j),(\alpha_p^{tN+1},\alpha_p^r)\big]$ is give by
\begin{equation*}
\begin{aligned}
&\Big(\sum_{s=1-r}^{0} \kappa_{d,p-l-s}\Big[\frac{h+1-x_s'}{e_1}\Big] +\sum_{s=0}^{N-r-1}\kappa_{d,p-l-s-1} \Big[\frac{h-x_{s+1}}{e_1}\Big]+N_1-\Big[\frac{j-1+e-\overline{p-l}}{e}\Big]\Big) (\alpha_l^{(h+t)N+1}, \alpha_l^{j+r-1})\\
&-\Big( \sum_{s=1-j}^{0} \kappa_{d,l-p-s}\Big[\frac{t+1-y_s'}{e_1}\Big]+\sum_{s=0}^{N-j-1}\kappa_{d,l-p-s-1}\Big[\frac{t-y_{s+1}}{e_1}\Big]+N_1-\Big[\frac{r-1+e-\overline{l-p}}{e}\Big]\Big)(\alpha_p^{(h+t)N+1},\alpha_p^{j+r-1}).
\end{aligned}
\end{equation*}

\end{itemize}

\end{Prop}

\medskip

From Propositions~\ref{Prop: gerstenhaber bracket1} and \ref{Prop:gerstenhaber bracket2}, we obtain the Gerstenhaber bracket of two chosen basis elements in $\rmHH^*(\Lambda)$.

\begin{Thm}\label{Thm:Gerstenhaber bracket of basic elements}

Let $\Lambda=\K Z_e/J^N$ be a truncated basic cycle algebra with $N,\, e\ge 2$. Let $m,n\ge 0$. Then the Gerstenhaber bracket of two basis elements in $\rmHH^m(\Lambda)$ and $\rmHH^n(\Lambda)$ $(m+n\ge 1)$ are given as follows:

\begin{enumerate}
	
\item ($[\rmHH^{0},\, \rmHH^{2k}]$) Assume that $m=0$ and $n=2k$ with $k\ge 1$.
	
\begin{itemize}

\item[(i)] Suppose that $N\not \equiv 1 \mode$ and that $\rmchar (\K)\nmid N$ or $kN\not \equiv N-1\mode$. Then for any $0\le i,j\le N-2$ with $i\equiv 0\mode$ and $j\equiv k N\mode$,
\[
\big[x_{0,i},\,x_{kN,j}\big]=0.
\]

\item[(ii)] Suppose that $N\not \equiv 1 \mode$ and that $\rmchar (\K)\di  N$ and $kN\equiv N-1\mode$. Then for any $0\le i\le N-2$ with $i\equiv 0\mode$ and  $0\le j\le N-1$ with   $j\equiv kN\mode$,
\[
\big[x_{0,i},\, x_{kN,j} \big]=0.
\]

\item[(iii)] Suppose that $N \equiv 1 \mode$ and that $\rmchar (\K)\nmid N$ or $kN\not \equiv N-1\mode$. Then

\begin{itemize}

\item[(a)] for any $1\le l\le e$ and $0\le j\le N-2$ with $j\equiv kN\mode$,
\[
\big[u_{l;\;0,N-1} ,\,  x_{k N,j}\big]=0.
\]
In particular, if $kN\not \equiv N-1\mode$, $\big[u_{l;\;0,N-1} ,\,  x_{k N,j}\big]$ is a coboundary;
            
\item[(b)] for any $0\le i,j\le N-2$ with $i\equiv 0\mode$ and $j\equiv kN\mode$,
\[
\big[x_{0,i},\, x_{kN,j}\big]=0.
\]

\end{itemize}

\item[(iv)] Suppose that $N \equiv 1 \mode$ and that $\rmchar (\K)\mid  N$ and $kN\equiv N-1\mode$. Then

\begin{itemize}

\item[(a)] for any $1\le l\le e$ and $0\le j\le N-1$ with $j\equiv kN \mode$,
\[
\big[u_{l;\;0,N-1} ,\,  x_{k N,j}\big]=0.
\]

\item[(b)] for any $0\le i\le N-2$ and  $0\le j\le N-1$ with  $i\equiv 0\mode$ and $j\equiv kN \mode$,
\[
\big[x_{0,i},\,x_{kN,j}\big]=0.
\]

\end{itemize}
		
\end{itemize}

\item $[\rmHH^{2h},\,\rmHH^{2k}]$ Assume that $m=2h$ and $n=2k$ with $h,k\ge 1$.

\begin{itemize}

\item[(i)] Suppse $\rmchar(\K)\nmid N$ or $\gcd(N,e)\not =1$. Then for any  $0\le i,j\le N-2$ with $i\equiv hN\mode$ and $j\equiv kN\mode$,
\[
\big[x_{hN,i},\, x_{kN,j}\big]=0.
\]			
		
\item[(ii)] Suppse $\rmchar(\K)\mid  N$ and $\gcd(N,e) =1$. Then for any  $0\le i,j\le N-1$ with $i\equiv hN\mode$ and $j\equiv kN\mode$,
\[
\big[x_{hN,i},\,x_{kN,j}\big]=0.
\]
		
\end{itemize}

\item $[\rmHH^{0},\,\rmHH^{2k+1}]$ Assume that $m=0$ and $n=2k+1$ with $k\ge 0$.

\begin{itemize}

\item[(i)] Suppose that $N\not \equiv 1 \mode$ and $\rmchar(\K)\di e$. Then for any $0\le i\le N-2$, $1\le j\le N-1$ with $i=i_1e\equiv 0\mode$ and $j\equiv kN+1\mode$,
\begin{equation*}
\big[x_{0,i},\, v_{1;\; kN+1,j}\big]=
\begin{cases}
-i_1\, x_{kN,i+j-1}, & \text{if }   i+j-1\le N-2 ,\\
0, &\text{otherwise}.
\end{cases}
\end{equation*}
In particular, if $i+j-1=N-1$, $\big[x_{0,i},\, v_{1;\; kN+1,j}\big]=-i_1\, x_{kN,N-1}$ is a coboundary.

\item[(ii)] Suppose that $N\not \equiv 1 \mode$, $\rmchar(\K)\nmid e$ and ($\rmchar(\K)\nmid N$ or $kN+1\not\equiv 0\mode$). Then for any $0\le i\le N-2$, $1\le j\le N-1$ with $i\equiv 0\mode$ and $j\equiv kN+1\mode$,
\begin{equation*}
\big[x_{0,i}, y_{kN+1,j}\big]=
\begin{cases}
-i\, x_{kN,i+j-1}, & \text{if }  i+j-1\le N-2 ,\\
0, &\text{otherwise}.
\end{cases}
\end{equation*}
In particular, if $i+j-1=N-1$, $\big[x_{0,i},\, v_{1;\; kN+1,j}\big]=-i\, x_{kN,N-1}$ is a coboundary.

\item[(iii)] Suppose that $N\not \equiv 1 \mode$, $\rmchar(\K)\di  N$ and $kN+1\equiv 0\mode$. Then for any $0\le i\le N-2$, $0\le j\le N-1$ with $i\equiv 0\mode$ and $j\equiv 0\mode$,
\begin{equation*}
\big[x_{0,i},\, y_{kN+1,j}\big]=
\begin{cases}
-i\, x_{kN,i+j-1}, & \text{if }  0\le i+j-1\le N-1\\
0, &\text{otherwise}.
\end{cases}
\end{equation*}

\item[(iv)] Suppose $N \equiv 1 \mode$ ($\ie,\; N=n_1 e+1$) and $\rmchar(\K)\di e$. Then

\begin{itemize}

\item[(a)] for any $1\le l\le e$ and $1\le j\le N-1$ with  $j\equiv kN+1\mode$,
\begin{equation*}
\big[u_{l;\;0,N-1},\,v_{1;\;k N+1,j}\big]=
\begin{cases}
-{n_1\, u_{l;\; 0,N-1}},  & \text{if } k=0  \text{ and } j=1\\
0, &\text{otherwise}.
\end{cases}	
\end{equation*}
In particular, if $k>0$ and $j=1$, 
\[
\big[u_{l;\;0,N-1},\,v_{1;\;k N+1,j}\big]=-n_1 (\alpha_l^{kN},\alpha_l^{N-1})=d_0^{2k-1}\big( n_1^2 \sum\limits_{p=1}^e (\alpha_p^{(k-1)N+1},e_p) -n_1 (\alpha_l^{(k-1)N+1},e_l) \big)
\]
is a coboundary;

\item[(b)] for any $0\le i\le N-2$, $1\le j\le N-1$ with $i=i_1 e\equiv 0\mode$ and $j\equiv kN+1\mode$,
\begin{equation*}
\big[x_{0,i},\, v_{1;kN+1,j}\big]=
\begin{cases}
-i_1\, x_{kN,i+j-1}, & \text{if } i+j-1\le N-2 ,\\
-i_1\sum\limits_{l=1}^{e}u_{l;\;0,N-1}, & \text{if } i+j-1= N-1 \text{ and } k=0,\\
0, &\text{otherwise}.
\end{cases}
\end{equation*}
In particular, if $k>0$ and $i+j-1=N-1$, 
\[
\big[x_{0,i},\, v_{1;kN+1,j}\big]=-i_1\, x_{kN,N-1} =-i_1 d_0^{2k-1}\big(  \sum\limits_{p=1}^e (\alpha_p^{(k-1)N+1},e_p) \big)
\]
is a coboundary.

\end{itemize}

\item[(v)] Suppose that $N \equiv 1 \mode$, $\rmchar(\K)\nmid e$ and ($\rmchar(\K)\nmid N$ or $kN+1\not\equiv 0\mode$). Then

\begin{itemize}

\item[(a)] for any $1\le l\le e$ and $1\le j\le N-1$ with  $j\equiv kN+1\mode$,
\begin{equation*}
\big[u_{l;\; 0,N-1},\, y_{kN+1,j}\big]=
\begin{cases}
-(N-1)\, u_{l;\; 0,N-1}, & \text{if } k=0 \text{ and } j=1 ,\\
0, &\text{otherwise}.
\end{cases}
\end{equation*}
In particular, if $k>0$ and $j=1$, $\big[u_{l;\; 0,N-1},\, y_{kN+1,j}\big]=-(N-1)(\alpha_l^{kN},\alpha_l^{N-1})$ is a coboundary;

\item[(b)] for any $0\le i\le N-2$, $1\le j\le N-1$ with $i\equiv 0\mode$ and $j\equiv kN+1\mode$,
\begin{equation*}
\big[x_{0,i}, y_{kN+1,j}\big]=
\begin{cases}
-i\, x_{kN,i+j-1}, & \text{if } i+j-1\le N-2 ,\\
-i\sum\limits_{l=1}^{e} u_{l;\;0,N-1}, & \text{if } i+j-1= N-1 \text{ and } k=0 ,\\
0, &\text{otherwise}.
\end{cases}
\end{equation*}
In particular, if $i+j-1=N-1$ and $k>0$, then $\big[x_{0,i}, y_{kN+1,j}\big]=-i \, x_{kN,N-1}$ is a coboundary.

\end{itemize}

\item[(vi)] Suppose that $N \equiv 1 \mode$ ($\ie,\; N=n_1 e+1$), $\rmchar(\K)\di N$ and $kN+1\equiv 0\mode$. Then

\begin{itemize}

\item[(a)] for any $1\le l\le e$ and $0\le j\le N-1$ with  $j\equiv 0\mode$,
\begin{equation*}
\big[u_{l;\; 0,N-1},\,y_{k N+1,j} \big]=
\begin{cases}
-n_1\, x_{kN,N-2}, & \text{if $j=0$},\\
0, &\text{otherwise}.
\end{cases}
\end{equation*}

\item[(b)] for any $0\le i\le N-2$, $0\le j\le N-1$ with $i\equiv 0\mode$ and $j\equiv 0\mode$,
\begin{equation*}
\big[x_{0,i},\, y_{kN+1,j}\big]=
\begin{cases}
-i\,x_{kN,i+j-1}, & \text{if $0\le i+j-1\le N-2$},\\
0, &\text{otherwise}.
\end{cases}
\end{equation*}
\end{itemize}
\end{itemize}

\item $[\rmHH^{2h},\rmHH^{2k+1}]$ Assume that $m=2h$ and $n=2k+1$, $h\ge 1$ and $k\ge 0$.

\begin{itemize}

\item[(i)] Suppose that $\rmchar(\K)\di e$. Then for any $0\le i\le N-2$, $1\le j \le N-1$, with $i\equiv hN\mode$ ($\ie \ hN-i=h_1 e$) and $j\equiv kN+1\mode$,
\begin{equation*}
\big[x_{hN,i},\, v_{1;\;kN+1,j}\big]=
\begin{cases}
h_1\,x_{(h+k)N,i+j-1}, & \text{if }  i+j-1\le N-2 ,\\
0, &\text{otherwise}.
\end{cases}
\end{equation*}
In particular, if $i+j-1=N-1$, $\big[x_{hN,i},\, v_{1;\;kN+1,j}\big]= h_1\,x_{(h+k)N,N-1}$ is a coboundary.

\item[(ii)] Suppose that $\rmchar(\K)\nmid e$ and ($\rmchar(\K)\nmid N$ or $\gcd(N,e)\not=1$). Then for any $0\le i\le N-2$, $1\le j \le N-1$, with $i\equiv hN\mode$ and $j\equiv kN+1\mode$,
\begin{equation*}
\big[x_{hN,i},\,y_{kN+1,j}\big]=
\begin{cases}
(hN-i)\,x_{(h+k)N,i+j-1}, & \text{if } i+j-1\le N-2 ,\\
0, &\text{otherwise}.
\end{cases}
\end{equation*}
In particular, if $i+j-1=N-1$, $\big[x_{hN,i},\, y_{kN+1,j}\big]= (hN-i)\,x_{(h+k)N,N-1}$ is a coboundary.

\item[(iii)] Suppose that $\rmchar(\K)\mid  N$ and $\gcd(N,e)=1$. Then
for any $0\le i,\,j\le N-1$ with $i\equiv hN\mode$ and $j\equiv kN+1\mode$,
\begin{equation*}
\big[x_{hN,i},\, y_{kN+1,j}\big]=
\begin{cases}
-i\, x_{(h+k)N,i+j-1}, & \text{if } 0\le i+j-1\le N-1 ,\\
0, &\text{otherwise}.
\end{cases}
\end{equation*}

\end{itemize}

\item $[\rmHH^{2h+1},\rmHH^{2k+1}]$  Assume that $m=2h+1$ and $n=2k+1$ with  $h, k\ge 0$.

\begin{itemize}

\item[(i)]  Suppose that $\rmchar(\K)\di e$. Then for any $1\le i,j\le N-1$ with $i\equiv hN+1\mode$ and $j\equiv kN+1\mode$ ($\ie,\; hN-i+1=h_1 e$ and $kN-j+1=k_1 e$ for some $h_1,k_1\ge 0$),
\begin{equation*}
\big[v_{1;\;hN+1,i},\, v_{1;\;kN+1,j}\big]=
\begin{cases}
(h_1-k_1)\, v_{1;\;(h+k)N+1,i+j-1}, & \text{if } i+j-1\le N-1 ,\\
0, &\text{otherwise}.
\end{cases}
\end{equation*}

\item[(ii)] Suppose that $\rmchar(\K)\nmid e$ and ($\rmchar(\K)\nmid N$ or $\gcd(N,e)\not=1$). Then for any $1\le i,j\le N-1$ with $i\equiv hN+1\mode$ and $j\equiv kN+1\mode$,
\begin{equation*}
\big[y_{hN+1,i},\, y_{kN+1,j}\big]=
\begin{cases}
(hN-kN-i+j)\, y_{(h+k)N+1,i+j-1}, & \text{if } i+j-1\le N-1 ,\\
0, &\text{otherwise}.
\end{cases}
\end{equation*}

\item[(iii)] Suppose that $\rmchar(\K)\di N$ and $\gcd(N,e)=1$. Then for any $0\le i,j\le N-1$ with $i\equiv hN+1\mode$ and $j\equiv kN+1\mode$,
\begin{equation*}
\big[y_{hN+1,i},\, y_{kN+1,j}\big]=
\begin{cases}
(-i+j)\, y_{(h+k)N+1,i+j-1}, & \text{if }  0\le i+j-1\le N-1 ,\\
0, &\text{otherwise}.
\end{cases}
\end{equation*}

\end{itemize}
	
\end{enumerate}	

\end{Thm}

To prove the above  proposition, we need a technical  lemma, whose proof is elementary and will be omitted.

\begin{Lem}\label{Lemma:two claims}

Let $N \ge 2$ be the integer such that $\gcd(N,e)=d\geq 1$. Write $N=dN_1$ and $e=de_1$. For fixed integer $s$, if $d\di p-s$, assume that $x_p\in [1,e_1]$ and $x_p'\in[2,e_1+1]$ and be the unique solution of the congruence $(x-1)N_1\equiv \frac{p-s}{d}\model$, respectively. Then

\begin{itemize}

\item[(1)]   $\{x_p\di 1\le p\le e, \text{ and } d\di p-s\}=\{1,\ldots,e_1\}$,
		
\item[(1')]  $\{x_p'\di 1\le p\le e, \text{ and } d\di p-s\}=\{2,\ldots,e_1+1\}$,
		
\item[(2)] for any $k\ge 0$, $\sum\limits_{x=1}^{e_1}[\frac{k+e_1-x}{e_1}]=k=\sum\limits_{x=1}^{e}[\frac{k+e-x}{e}]$.

\end{itemize}

\end{Lem}






\begin{SProof}[Sketch of proof of Theorem~\ref{Thm:Gerstenhaber bracket of basic elements}]

We only prove the cases (1)(iv)(a) and (4)(i), while the others are similar or follow directly from Lemma~\ref{Lemma:two claims} by straightward computation.

\medskip

\textbf{Case (1)(iv)(a):} By assumption, $N\equiv 1\mode$, which implies $\gcd(N,e)=1$. For $1\le l\le e$ let $y_p\in [1,e]$ (resp. $y_p'\in [1,e]$) be the unique solution of the congruence $(x-1)N\equiv l-p-1\mode$ (resp. $(x-1)N\equiv l-p\mode$).  Following from Proposition ~\ref{Prop: gerstenhaber bracket1} (1), we have
\begin{equation*}
\begin{aligned}
\big[u_{l;\;0,N-1},\,x_{kN,j}\big]
&=\sum\limits_{p=1}\limits^e\Big(\Big[\frac{k+e-y_p'}{e}\Big]-\Big[\frac{k+e-y_p}{e}\Big]\Big) (\alpha_p^{(k-1)N+1},\alpha_p^j).
\end{aligned}
\end{equation*}
If $j=0$, then $e\di k$, we have 
\[
\Big[\frac{k+e-y_p'}{e}\Big]=\Big[\frac{k+e-y_p}{e}\Big]=\frac{k}{e}, \text{ for all } p\in [1,e],
\]
then $\big[u_{l;\;0,N-1},\,x_{kN,j}\big]=0$. If $1\le j\le N-1$, it is obvious that $y_p=y_{\ul{p+1}}'$ for $1\le p\le e$, then
\begin{equation*}
\begin{aligned}
\big[u_{l;\;0,N-1},\, x_{kN,j}\big]
&=\sum\limits_{p=1}\limits^e\Big[\frac{k+e-y_p'}{e}\Big]\big[(\alpha_p^{(k-1)N+1},\alpha_p^j)-(\alpha_{\ul{p-1}}^{(k-1)N+1},\alpha_{\ul{p-1}}^j)\big]\\
&= d_{j-1}^{2k-2} \big( \sum\limits_{p=1}\limits^e \Big[\frac{k+e-y_p'}{e}\Big] (\alpha_p^{(k-1)N},\alpha_p^{j-1}) \big) \\
 &\equiv 0.
	\end{aligned}
\end{equation*}

\medskip

\textbf{Case (4)(i):} Following from Proposition~\ref{Prop: gerstenhaber bracket1} (2), if $\gcd(N,e)=1$, then
\begin{equation*}
	\begin{aligned}
	\big[x_{hN,i},\, v_{1;\,kN+1,j}\big]
    =& \sum\limits_{l=1}^e \Big(\sum\limits_{s=0}^{N-1} \Big[\frac{h+e-x_{s^l}}{e}\Big] -\Big[\frac{i-1+e-\overline{1-l}}{e}\Big]\Big)(\alpha_l^{(h+k)N},\alpha_l^{i+j-1}),
	\end{aligned}
\end{equation*}	
where $x_{s^l}\in [1, e]$ is the unique solution of the congruence $(x-1)N+l+s\equiv 1\mode$ for each $1\le l\le e$. Set 
\[
\lambda_l=\sum\limits_{s=0}^{N-1} \Big[\frac{h+e-x_{s^l}}{e}\Big] -\Big[\frac{i-1+e-\overline{1-l}}{e}\Big].
\]
We now prove that $\lambda_l=\lambda_{l+1}$ for all $l\in [1,e-1]$. For $l=1$, we have
\[
\lambda_1-\lambda_2= \Big[\frac{h+e-1}{e} \Big]- \Big[\frac{h}{e} \Big]+\Big[\frac{i}{e} \Big]-\Big[\frac{i+e-1}{e} \Big].
\]
We distinguish two cases:
\begin{itemize}
    \item if $i\equiv 0\mode$, then also $h\equiv 0\mode$. In this case, 
    \[
    \Big[\frac{h+e-1}{e} \Big]= \Big[\frac{h}{e} \Big]\quad \text{and}\quad \Big[\frac{i}{e} \Big]=\Big[\frac{i+e-1}{e} \Big],
    \]
    and therefore $\lambda_1=\lambda_2$.

    \item if $i\equiv 0\mode$, then $h\equiv 0\mode$. In this case, 
    \[
    \Big[\frac{h+e-1}{e} \Big]= \Big[\frac{h}{e} \Big]+1 \quad \text{and}\quad \Big[\frac{i}{e} \Big]+1=\Big[\frac{i+e-1}{e} \Big],
    \]
    so again $\lambda_1=\lambda_2$.
\end{itemize}

In the same way, we obtain $\lambda_l=\lambda_{l+1}$ for $2\le l\le e-1$. By Lemma~\ref{Lemma:two claims} we have $\sum\limits_{l=1}^e\lambda_l=hN-i$, and hence $\lambda_l=h_1$ for each $l\in [1,e]$, where $hN-i=h_1 e$.

If $\gcd(N,e)=d>1$, then by Proposition ~\ref{Prop:gerstenhaber bracket2}~(2), the proof proceeds in the same way.

\end{SProof}

By Theorem~\ref{Thm:Gerstenhaber bracket of basic elements}, we can also determine the Gerstenhaber algebra structures on the Hochschild cohomology ring $\rmHH^*(\Lambda)$ of a truncated basic cycle algebra $\Lambda=\K Z_e/J^N$, where $\gcd(N,e)=d\ge 1$, $N=dN_1$ and $e=de_1$. We then distinguish $8$ different cases for further analysis.

\begin{Prop}\label{Prop:Gerstenhaber algebraic structure 1}

Assume that $N\le e$. If $\rmchar(\K) \nmid N$ or $\gcd(N,e)\not= 1$, then the graded commutative algebra $\rmHH^*(\Lambda)$ in Proposition~\ref{Prop:ring structure1} is a Gerstenhaber algebra, where the Gerstenhaber brackets between generators are given as follows: for each $i$, let $k'_i$ denotes the integer satisfying $k_i N-di=k'_i e$. Then
\[
[y,z_i]=
\begin{cases}
-(k_iN-di)z_i,     & \text{if }\rmchar(\K)\nmid e, \\
- k'_i z_i     & \text{if }\rmchar(\K)\di e,
\end{cases} \quad \text{for all $0\le i\le [\frac{N-2}{d}]$},
\]
with all other Gerstenhaber brackets between basis elements equal to zero. 

\end{Prop}

\begin{Proof}

The generators and relations of $\rmHH^*(\Lambda)$ are described in Proposition~\ref{Prop:ring structure1}. In what follows, we compute the Gerstenhaber brackets among the generators.

By the isomorphism constructed in Proposition~\ref{Prop:ring structure1}, if $\rmchar(\K)\nmid e$, $y$ corresponds to $y_{1,1}$, and Theorem~\ref{Thm:Gerstenhaber bracket of basic elements} (5)(ii) gives $\big[y_{1,1},\, y_{1,1}\big]=0$. For any $0\le i,j\le [\frac{N-2}{d}]$, $[z_i,z_j]$ corresponds to $\big[ x_{k_i N, di},\, x_{k_j N, dj}\big]$, which vanishes by Theorem~\ref{Thm:Gerstenhaber bracket of basic elements} (2)(i). For $0\le i\le [\frac{N-2}{d}]$, $[y,z_i]$ corresponds to $\big[y_{1, 1},\, x_{k_iN, di}\big]=-(k_iN-di) x_{k_iN, di}$, which follows from Theorem~\ref{Thm:Gerstenhaber bracket of basic elements}  (4)(ii).

If $\rmchar(\K)\di e$, the result can be obtained in the same manner.
\end{Proof}

\begin{Prop}\label{Prop:Gerstenhaber algebraic structure 2}

Assume that $N\le e$. If $\rmchar(\K) \di N$ and $\gcd(N,e)=1$. Then the graded commutative algebra $\rmHH^*(\Lambda)$ in Proposition~\ref{Prop:ring structure2} is a Gerstenhaber algebras, where the Gerstenhaber bracket between generators are given as follows: 
\[
[y,z_i]=iz_i,\; [y,w]=-w,\; \text{and}\; [z_i,w]=
		\begin{cases}
			-i z_{i-1},&\text{if }  k_i+I-1\le e\\
			-i z_0z_{i-1},&\text{if } k_i+I-1> e
		\end{cases}	
\; \text{for all } 0\le i\le N-1,
\]
with all other Gerstenhaber brackets between basis elements vanish. 

\end{Prop}

\begin{Proof} 

The generators and relations are described in Proposition~\ref{Prop:ring structure2}, in which the isomorphism sends $y$ to $y_{1,1}$, $z_j$ to $x_{k_jN,j}$ and $w$ to $y_{(I-1)N+1,0}$.

Similar to the proof of Proposition~\ref{Prop:Gerstenhaber algebraic structure 1}, $[y,y]=0$ and $[y,z_i]=iz_i$ for $0\le i\le N-1$. It is obvious that $[w,w]=0$ and 
\[
[y,w]=[y_{1,1},\, y_{(I-1)N+1,0} ]= -w
\]
by Theorem~\ref{Thm:Gerstenhaber bracket of basic elements} (5)(iii).

For $0\le i,j\le N-1$, 
\[
[z_i,z_j]=[x_{k_iN,i},\, x_{k_jN,j}]=0
\]
by Theorem~\ref{Thm:Gerstenhaber bracket of basic elements} (2)(ii).

For $0\le i\le N-1$, by  Theorem~\ref{Thm:Gerstenhaber bracket of basic elements} (4)(iii). 
\[
[z_i,w]=[x_{k_iN,i}, y_{(I-1)N+1,0}]=
-i\, x_{(k_i+I-1)N, i-1}.
\]
Assume that $1\le i\le N-1$. If $k_i+I-1\le e$, by definition of $k_{i-1}$ we have $k_{i-1}=k_i+I-1$, then   $[z_i,w]=-iz_{i-1}$; if $k_i+I-1> e$, then $k_i+I-1=k_{i-1}+k_0$, by Proposition~\ref{Prop: cup product of basis elements} (3)(ii), $[z_i,w]=-i\, x_{k_0N,0}\vee x_{k_{i-1}N,i-1}=-iz_0z_{i-1}$.

\end{Proof}

The proofs of the following Propositions follow directly from Theorem~\ref{Thm:Gerstenhaber bracket of basic elements}, as in the cases of Proposition~\ref{Prop: gerstenhaber bracket1} and ~\ref{Prop:Gerstenhaber algebraic structure 2}. Hence, we omit the details.

\begin{Prop}\label{Prop:Gerstenhaber algebraic structure 3}

Assume that $N> e$ and $N\not\equiv 1\mode$. If $\rmchar(\K) \nmid N$ or  $\gcd(N,e)\neq 1$, then the graded commutative algebra $\rmHH^*(\Lambda)$ in Proposition~\ref{Prop:ring structure3} is a Gerstenhaber algebra, where the Gerstenhaber brackets between generators are given as follows: for each $i$, let $k'_i$ denotes the integer satisfying $k_i N-di=k'_i e$. Then for all $0\le i\le e_1-1$
\[
[x_0,y]=
\begin{cases}
-ex_0,     & \text{if } \rmchar(\K)\nmid e \\
- x_0 ,    & \text{if } \rmchar(\K)\di e
\end{cases}
\quad \text{and}\quad 
[y,z_i]=
\begin{cases}
-(k_iN-di)\, z_i,     & \text{if } \rmchar(\K)\nmid e \\
-k_i' \, z_i,     & \text{if } \rmchar(\K)\di e
\end{cases}
\]
with all other Gerstenhaber brackets between basis elements equal to zero. 

\end{Prop}

\begin{Ex}\label{Example: Gerstenhaber brackets}

Let $\Lambda$ be the algebra of Example~\ref{Example: generators and relations}. Then
\[
\rmHH^*(\Lambda)\cong \K[x,y,z]/\langle x^2,y^2\rangle
\]
with $|x| =0$, $|y| =1$ and $|z| =2$. By Proposition~\ref{Prop:Gerstenhaber algebraic structure 3}, the Gerstenhaber brackets between generators are trivial except that $[x,y]=-x$.

\end{Ex}

\begin{Prop}\label{Prop:Gerstenhaber algebraic structure 4}

Assume that $N>e$, $N\not\equiv 1\mode$, $\rmchar(\K)\di N$ and $\gcd(N,e)=1$, then the graded commutative algebra $\rmHH^*(\Lambda)$ described in Proposition~\ref{Prop:ring structure4} is a Gerstenhaber algebra with Gerstenhaber brackets between generators are given as follows:
\[
 [x_0,y]=-ex_0,\, [x_0,w]=-ez_{e-1},\,  [y,w]=-w,
\]
and for $0\le i\le e-1$, we have $[y,z_i]=iz_i$ and
\begin{equation*}
	[z_i,w]=
	\begin{cases}
		-i z_{i-1},&\text{if } 1\le i\le e-1 \text{ and } k_i+I-1\le e,\\
		-i z_0z_{i-1},&\text{if } 1\le i\le e-1 \text{ and } k_i+I-1>e,\\
		0,&\text{if } i=0.
	\end{cases}
\end{equation*}
All other Gerstenhaber brackets between basis elements equal to zero. 

\end{Prop}

\begin{Prop}\label{Prop:Gerstenhaber algebraic structure 5}

Assume that $N=e+1$ and $\rmchar(\K)\nmid N$, then the graded commutative algebra $\rmHH^*(\Lambda)$ in Proposition~\ref{Prop: ring structure5} is a Gerstenhaber algebra with Gerstenhaber brackets between generators are given as follows:
\[
[\widetilde{x_l},y]=
\begin{cases}
-e\widetilde{x_l},     & \text{if } \rmchar(\K)\nmid e \\
- \widetilde{x_l},    & \text{if } \rmchar(\K)\di e
\end{cases} \quad \text{for all }1\le l\le e,
\]
\[
[y,z]=
\begin{cases}
-ez ,    &  \text{if } \rmchar(\K)\nmid e \\
-z,     &  \text{if } \rmchar(\K)\di e
\end{cases}
\quad \text{and}\quad [y,z']=
\begin{cases}
-eNz' ,    &  \text{if } \rmchar(\K)\nmid e \\
-z',     &  \text{if } \rmchar(\K)\di e
\end{cases}.
\]
    
\end{Prop}

\begin{Rem}\label{Remark: example for P not be a dg lie algebra}

As stated at the end of the previous  Subsection, on the chain level, $s\bbP^*$ endowed with the differential $d^*$ and the minimal bracket $[-,-]_\bbP$ is NOT a differential graded Lie algebra. Assume that $e=2$, $N=3$ and $\rmchar(\K)=2$. Let $a=(e_1,\alpha_1^2)\in \bbP^0$, $b=(\alpha_1^3,\alpha_1)\in \bbP^2$ and $c=(\alpha_1^4,\alpha_1^2)\in \bbP^{3}$. Then by Proposition~\ref{Prop: gerstenhaber bracket1}, we have 
	\begin{align*}
		\big[[a,b],c\big]&=[(\alpha_1,\alpha_1),(\alpha_1^4,\alpha_1^2)]=-(\alpha_1^4,\alpha_1^2)=-c,\\
		\big[[b,c],a\big]&=[(\alpha_1^6,\alpha_1^2),(e_1,\alpha_1^2)]=0,\\
		{[c,a]}&=0.
	\end{align*}
So $a$, $b$ and $c$ do not satisfy the graded Jacobi identity.

\end{Rem}

\begin{Prop}\label{Prop:Gerstenhaber algebraic structure 6}

Assume that $N=e+1$ and $\rmchar(\K)\mid N$. Then the graded commutative algebra $\rmHH^*(\Lambda)$ in Proposition~\ref{Prop: ring structure 6} is a Gerstenhaber algebra with all the nontrivial Gerstenhaber brackets between generators are given as follows:
\[
[\widetilde{x_l},y]=\widetilde{x_l},\; [\widetilde{x_l},w]=-z^{e-1},\,
[y,z]=z,\; [y,w]=-w,\; [z,w]=-z',
\]
for $1\le l \le e$.
    
\end{Prop}

\begin{Prop}\label{Prop:Gerstenhaber algebraic structure 7}

Assume that $N=me+1$, $m>1$ and $\rmchar(\K)\nmid N$, then the graded commutative algebra $\rmHH^*(\Lambda)$ in Propositions~\ref{Prop: ring structure 7} is a Gerstenhaber algebra, and all the nontrivial Gerstenhaber brackets between generators are given as follows: for $1\le l \le e$,
\begin{itemize}
    \item if $\rmchar(\K)\nmid e$, 
    \[
    [x_0,y]= -ex_0,\, [\widetilde{x_l},y]=-me\widetilde{x_l},\,[y,z]=-mez, \,[y,z']=-eNz',
    \]

    \item if $\rmchar(\K)\di e$, 
    \[
    [x_0,y]= -x_0,\, [\widetilde{x_l},y]=-m\widetilde{x_l},\,[y,z]=-mz, \,[y,z']=-z'.
    \]
\end{itemize}

\end{Prop}

\begin{Prop}\label{Prop:Gerstenhaber algebraic structure 8}

Assume that $N=me+1$, $m>1$ and $\rmchar(\K)\mid  N$, then the graded commutative algebra $\rmHH^*(\Lambda)$ in Propositions~\ref{Prop: ring structure 8} is a Gerstenhaber algebra, and all the nonzero Gerstenhaber brackets between generators are given as follows: for $1\le l\le e$,
\[
[x_0,y]=-e x_0,\; [x_0,w]=-ez^{e-1},\;[\widetilde{x_l},y]= \widetilde{x_l},\;[\widetilde{x_j},w]=-m z^{e-1} x_0^{m-1},
\]
\[
\text{and}\quad [y,w]=-w,\; [y,z]=z,\;[z,w]=-z'.
\]
\end{Prop}

\bigskip

\section{BV structure}\label{section:BV-structure}\

In this section, we determine the BV algebra structure on the Hochschild cohomology of a selfinjective Nakayama algebra.

\medskip

Let $\Lambda=\K Z_e/J^N$ be the truncated basic cycle algebra with $N,\,e\ge 2$. Then $\Lambda$ is a Frobenius algebra, endowed with the associative nondegenerate $\K$-bilinear form $\langle\ ,\ \rangle\colon \Lambda\times \Lambda \rightarrow \K$: 
\begin{equation*}
\langle\alpha_l^i,\alpha_{l'}^j \rangle=
	\begin{cases}
		1,&\text{if } \, i+j=N-1 \text{ and } \, l'\equiv l+i \mode ;\\
		0, &\text{otherwise},
	\end{cases}
\end{equation*}
for $\alpha_l^i,\,\alpha_{l'}^j\in \Lambda$ with $1\le l,\, l'\le e,\; 0\le i,\, j\le N-1$. The Nakayama automorphism $\sigma\colon\Lambda\rightarrow\Lambda$ of $\Lambda$ associated to $\langle\ ,\ \rangle$ is defined by, for each $\alpha_l^i \in \Lambda$ with $1\le l\le e,\; 0\le i\le N-1$, one has
\[
\sigma(\alpha_l^i)=\alpha_{\ul{l+N-1}}^i.
\]


We begin by analyzing the conditions under which $\sigma$ is semisimple. This consideration yields the following criterion.

\begin{Cri}\label{Criterion: semisimple Nakayama automorphism}

Let $\Lambda=\K Z_e/J^N$ be the truncated basic cycle algebra. $\sigma$ is the Nakayama automorphism of $\Lambda$ described as above. Then $\sigma$ is semisimple if and only if it satisfies one of the following conditions:

\begin{enumerate}

    \item   $\rmchar(\K)\nmid e$;
		
	\item $2\le \rmchar(\K)\le e$,  $\rmchar(\K)\di e$, $\gcd(N-1,e)=c>1$ and $\rmchar(\K)\nmid \frac{e}{c}$.
    
\end{enumerate}

\end{Cri}

\begin{Proof}

Assume that $\rmchar(\K)\nmid e$, if $\rmchar(\K)=0$ or $\rmchar(\K)>e$, from \cite[Criterion 5.1]{LZZ16} we know that $\sigma$ is semisimple.

Take a $\K$-basis $ \{e_1,\ldots,e_e,\alpha_1,\ldots,\alpha_e,\ldots,\alpha_1^{N-1}, \ldots, \alpha_e^{N-1}\}$ of $\Lambda$. By the definition of $\sigma$ we can write $\sigma=\sigma_1 \times \sigma_2\times\cdots\times\sigma_N$, where $\sigma_i=\sigma\mid _{\K\{\alpha_1^{i-1},\ldots,\alpha_e^{i-1}\}}$ for $1\le i\le N$ (with convention $\alpha_l^0=e_l$). For all $1\le i\le N$, the representative matrix of the automorphism $\sigma_i$ has the form
	$$ \begin{bmatrix}
		0&\cdots & 0& 1&0&\cdots &  0\\
		0&\cdots & 0& 0&1& \cdots &0\\
		\vdots&\vdots &\vdots&\vdots &\vdots &\ddots&\vdots\\
		0&\cdots & 0& 0&0&\cdots &  1\\
		1&\cdots&0 &0& 0&\cdots & 0\\
		\vdots&\vdots &\vdots&\vdots &\vdots &\ddots&\vdots\\
		0&\cdots&1 &0& 0&\cdots & 0
	\end{bmatrix}_{e}$$
In the matrix above, the entry in the $j$-th row and $\ul{j+N-1}$-th column is 1 for $1\le j\le e$, and all other entries are 0. It is straightforward to verify that this matrix can be written as the product of $N-1$ matrices $M$ of the following from
\[
M=
\begin{bmatrix}
	0 & 1& 0&\cdots &  0\\
	0 & 0& 1& \cdots &0\\
	\vdots &\vdots &\vdots &\ddots&\vdots\\
	0 & 0& 0&\cdots &  1\\
	1 &0 &0& \cdots & 0
\end{bmatrix}_{e}
\]
$\ie$, the matrix of $\sigma_i$ is $M^{N-1}$.

Assume that $2\le \rmchar(\K)\le e$ and $\rmchar(\K)\nmid e$. To prove $\sigma$ is semisimple, it suffice to prove that $M$ is semisimple over the field $\K$. By simple calculation, we know that the characteristic polynomial of matrix $M$ is $\lambda^e-1$, which has $e$ different roots in $\K$ because $\rmchar(\K)\nmid e$. This means $M$ is semisimple over $\K$, then we have proved that if condition (1) is satisfied, $\sigma$ is semisimple.

Assume that $2\le \rmchar(\K)\le e$ and $\rmchar(\K)\mid e$.

If $\gcd(N-1,e)=1$, there exist $1\le a\le e-1$ such that $M^{(N-1)a}=M$ because $M^e=E$. So in order to prove $M^{N-1}$ is not semisimple, we just need to prove that $M$ is not semisimple. Since $\rmchar(\K)\mid e$, we have $(\lambda-1)^2\mid \lambda^e-1$, then the algebraic multiplicity (AM for short) of the eigenvalue $1$ is greater than or equal to $2$, $\ie\ AM(1)\ge 2$. And it is easy to see the geometric multiplicity (GM for short) of $1$ is $1$, then $1=GM(1)<AM(1)$, which means that $M$ is not semisimple over $\K$.

If $\gcd(N-1,e)=c>1$, there exist $1\le a\le e-1$ such that $M^{(N-1)a}=M^c$. The characteristic polynomial of matrix $M^c$ is $(\lambda^{\frac{e}{c}}-1)^c$. If $\rmchar(\K)\nmid \frac{e}{c}$, let $w$ be one of the roots of the polynomial $\lambda^{\frac{e}{c}-1}+\cdots+\lambda+1$, it's not difficult to verify that $GM(w^l)=c=AM(w^l)$ for $0\le l\le \frac{e}{c}-1$. Then $M^c$ is semisimple over $\K$. If $\rmchar(\K)\mid  \frac{e}{c}$, $AM(1)\ge 2c>GM(1)=c$. Then $M^c$ is not semisimple over $\K$.
This completes the proof of the criterion.

\end{Proof}

\begin{Rem}

The condition in the Criterion~\ref{Criterion: semisimple Nakayama automorphism} is equivalent to the condition that $\rmchar(\K)\nmid \operatorname{ord}\,(\sigma)$.
    
\end{Rem}

\medskip

Now we assume that $\sigma$ is semisimple, i.e., $\rmchar(\K)$ satisfies condition (1) or (2) in Criterion~\ref{Criterion: semisimple Nakayama automorphism}. From Section~\ref{Section: preliminaries}, $\rmHH^*(\Lambda)$ is a BV-algebra.
For $f'\in H_1^n(\Lambda,\Lambda)$, the corresponding element be written as $f\in \bbB^n=\Hom_{A^e}(\Lambda\ot {\Lambda_+}^{\ot n}\ot \Lambda,\Lambda)$, with $n\ge 1$,
\[
 \Delta_{\bbB}(f)(a\ot a_1\ot\cdots\ot a_{n-1}\ot b) =\sum_{l=1}^e \sum_{j=0}^{N-1} \sum_{i=1}^n  (-1)^{i(n-1)} \langle f'(a_{i}\ot \cdots\ot a_{n-1}\ot \alpha_l^j \ot \sigma(a_1)\ot \cdots\ot \sigma(a_{i-1})), \; 1\rangle a\alpha_{\ul{l+j-N+1}}^{N-1-j} b.
\]
Then the differential $\Delta$ on the corresponding element in $\Hom_{\Lambda^e}(\Lambda\ot \K Q_{\chi(n)}\ot \Lambda,\Lambda)$ is defined by the following commutative diagram:
\begin{eqnarray*}
\xymatrix{
\Hom_{\Lambda^e}(\Lambda\ot\K Q_{\chi(n)}\ot\Lambda,\Lambda) \ar[d]_{\omega_n^*} \ar[r]^{\Delta_n} 
&\Hom_{\Lambda^e}(\Lambda\ot\K Q_{\chi(n-1)}\ot \Lambda,\Lambda)\\
\Hom_{{\Lambda^e}}(\Lambda\ot\Lambda_+^{\ot n}\ot \Lambda,\Lambda) \ar[r]^{\Delta_{\bbB}}
&\Hom_{{\Lambda^e}}(\Lambda\ot\Lambda_+^{\ot n-1}\ot \Lambda,\Lambda) \ar[u]_{\mu_{n-1}^*},}
\end{eqnarray*}
that is, $\Delta_n=\mu_{n-1}^*\circ \Delta_{\bbB} \circ \omega_n^*$.

Before calculating, we still need to analyze the basis elements of $H_1^*(\Lambda,\Lambda)$.

Let $\lambda\in I_{\Lambda}$ and $\Lambda_{\lambda}$ be the eigenspace $\Ker(\sigma-\lambda \mathrm{id})$ associated with $\lambda$.

If $\rmchar(\K)\nmid e$ and $\gcd(N-1,e)=1$, by the proof of Criterion~\ref{Criterion: semisimple Nakayama automorphism}, the representative matrix of the Nakayama automorphism $\sigma$ restricted on $\K\{\alpha_1^{i},\ldots,\alpha_e^i\}$ is $M^{N-1}$, 
and it is similar to a diagonal matrix
$$\begin{bmatrix}
	1&0&0&\cdots&0\\
	0&w&0&\cdots&0\\
	0&0&w^2&\cdots&0\\
	\vdots&\vdots &\vdots&\ddots&\vdots\\
	0&0&0&\cdots&w^{e-1}
\end{bmatrix}_e$$
where the element $w\in\K$ is one of roots of the polynomial $1+x+\cdots+x^{e-1}$, then we have a decomposition $\Lambda=\Lambda_1\oplus \Lambda_{w}\oplus\cdots\oplus \Lambda_{w^{e-1}}$ of $\Lambda$ by $e$ $\K$-vector spaces
\begin{equation*}
	\begin{array}{ll}
		\Lambda_1=\K 1_{\Lambda}\oplus \K(\sum\limits_{l=1}^e\alpha_l)\oplus\cdots\oplus\K(\sum\limits_{l=1}^e\alpha_l^{N-1}),\\
		\Lambda_w=\bigoplus_{i=0}^{N-1}\K(\sum\limits_{j=0}^{e-1}w^j\alpha_{\ul{1+j(N-1)}}^i),\\
		\Lambda_{w^2}=\bigoplus_{i=0}^{N-1}\K(\sum\limits_{j=0}^{e-1}w^{2j}\alpha_{\ul{1+j(N-1)}}^i),\\
		\vdots\\
		\Lambda_{w^{e-1}}=\bigoplus_{i=0}^{N-1}\K(\sum\limits_{j=0}^{e-1}w^{(e-1)j}\alpha_{\ul{1+j(N-1)}}^i).
	\end{array}
\end{equation*}
Moreover, for $1\le k\le e-1$, $\Lambda_{w^k+}=\Lambda_{w^k}$ and $\Lambda_{1+}= \K(\sum\limits_{l=1}^e\alpha_l)\oplus\cdots\oplus\K(\sum\limits_{l=1}^e\alpha_l^{N-1})$.

If ($\rmchar(\K)\nmid e$ and $\gcd(N-1,e)=c>1$) or ($\rmchar(\K)\mid e$, $\gcd(N-1,e)=c>1$ and $\rmchar(\K)\nmid \frac{e}{c}$), assume that $N-1=cN_2$ and $e=ce_2$, then by similar analysis, we have a decomposition  $\Lambda=\Lambda_1\oplus\Lambda_{v}\oplus\cdots\Lambda_{v^{e_2-1}}$ of $\Lambda$ by $e_2$ $\K$-vector spaces
\begin{equation*}
	\begin{array}{ll}
		\Lambda_1=\bigoplus_{i=0}^{N-1}\bigoplus_{p=1}^{c}\K (\sum\limits_{j=0}^{e_2-1}\alpha_{\ul{p+j(N-1)}}^i),\\
		\Lambda_v=\bigoplus_{i=0}^{N-1}\bigoplus_{p=1}^{c}\K (\sum\limits_{j=0}^{e_2-1}v^j\alpha_{\ul{p+j(N-1)}}^i),\\
		\Lambda_{v^2}=\bigoplus_{i=0}^{N-1}\bigoplus_{p=1}^{c}\K (\sum\limits_{j=0}^{e_2-1}v^{2j}\alpha_{\ul{p+j(N-1)}}^i),\\
		\vdots\\
		\Lambda_{v^{e_2-1}}=\bigoplus_{i=0}^{N-1}\bigoplus_{p=1}^{c}\K (\sum\limits_{j=0}^{e_2-1}v^{(e_2-1)j}\alpha_{\ul{p+j(N-1)}}^i),
	\end{array}
\end{equation*}
where $v$ is one of roots of the polynomial $1+x+\cdots+x^{e_2-1}$.

For $\lambda\in I_{\Lambda}$, write $\Lambda_{\lambda+}=\Lambda_{\lambda}$ for $\lambda\neq 1$ and $\Lambda_{1+}=\Lambda_1/(\K Q_0)$, let $C_{(1)}^n(\Lambda,\Lambda)$ be those Hochschild cochians $\phi\in C^n(\Lambda,\Lambda)$ such that we have $\phi(\Lambda_{\lambda_1+}\ot\cdots \ot \Lambda_{\lambda_n+})\subset \Lambda_{\lambda_1\cdots\lambda_n}$ for all eigenvalues $\lambda_1,\ldots,\lambda_n\in I_{\Lambda}$. The restriction $b_{(1)}^n$ of the Hochschild differential $b^n: C^n(\Lambda,\Lambda)\to C^{n+1} (\Lambda,\Lambda)$ to $C_{(1)}^n(\Lambda,\Lambda)$ has values in $C_{(1)}^{n+1}(\Lambda,\Lambda)$. Put
\[
\rmH_{(1)}^n(\Lambda,\Lambda):= \rmH^n(C_{(1)}^*(\Lambda,\Lambda),b_{(1)}^*).
\]
The subcomplex $(C_{(1)}^*(\Lambda,\Lambda),b_1^*)$ of $(C^*(\Lambda,\Lambda),b^*)$ defines a morphism of graded vector spaces
\[
\rmH_{(1)}^*(\Lambda,\Lambda)\to \rmHH^*(\Lambda).
\]
Moreover, if $\sigma$ is semisimple, this is an isomorphism. We set 
\[
\Hom_{\Lambda^e}^{(1)}(P_n,\Lambda)=\{ f\in \Hom_{\Lambda^e}(P_n,\Lambda)\di \omega^n(f)=f\circ \omega_n \in C_{(1)}^n(\Lambda,\Lambda) \}.
\]
Then we obtain the following Lemma.

\begin{Lem}
Let $\Lambda={\K Z_e}/J^N$ with $N ,\, e\ge 2$, and the Nakayama automorphism $\sigma$ is semisimple. Then 
\begin{itemize}
    \item [(1)] if $\rmchar(\K)\nmid e$ and $\gcd(N-1,e)=1$, then
    \[
    \Hom_{\Lambda^e}^{(1)}(P_n,\Lambda)\cong \K \{ x_{\chi(n),i}\di 0\le i\le N-1 \text{ and } \chi(n)\equiv i\mode \},
    \]

    \item [(2)] if ($\rmchar(\K)\nmid e$ and $\gcd(N-1,e)=c>1$) or ($\rmchar(\K)\mid e$, $\gcd(N-1,e)=c>1$ and $\rmchar(\K)\nmid \frac{e}{c}$), assume that $N-1=cN_2$ and $e=ce_2$, then
    \[
    \Hom_{\Lambda^e}^{(1)}(P_n,\Lambda)\cong \K \left\{ \sum\limits_{j=0}^{e_2-1}\big( \alpha_{\ul{l+j(N-1)}}^{\chi(n)}, \alpha_{\ul{l+j(N-1)}}^i\big)\, \middle|\, 1\le l\le c, \; 0\le i\le N-1, \text{ and } \chi(n)\equiv i\mode  \right\}.
    \]
\end{itemize}  

\end{Lem}

The restriction $d_{(1)}^n$ of the Hochschild differential $d^n: \Hom_{\Lambda^e}(P_n,\Lambda)\to \Hom_{\Lambda^e}(P_{n+1},\Lambda)$ to $\Hom_{\Lambda^e}^{(1)}(P_n,\Lambda)$ has values in $\Hom_{\Lambda^e}^{(1)}(P_{n+1},\Lambda)$.
Then we can obtain the basis elements of vector space $\rmH_{(1)}^*(\Lambda,\Lambda)$ by using Corollary~\ref{Cor:new description of P for basic cycle}.

\begin{Cor}\label{Thm: basis of HH1} 
Let $\Lambda={\K Z_e}/J^N$ with $N ,\, e\ge 2$, and the Nakayama automorphism $\sigma$ is semisimple. The following statements hold:
\begin{itemize}

\item [(1)] for $k\ge 1$,
      \begin{equation*}
	 \rmHH^{2k}_{(1)}(\Lambda)=
	 \begin{cases}
	 \K  \left\{x_{kN, i}\di 0\le i\le N-1,\; i \equiv N-1 \mode\right\} , & \text{if }    \rmchar(\K) \di N  \text{ and }  k N  \equiv N-1 \mode,\\
	 \K\left\{x_{kN, i}\di 0\le i\le N-2,\; i \equiv kN \mode\right\} , &\text{otherwise};
	 \end{cases}
      \end{equation*}

\item [(2)] for $k\ge 0$, if $\rmchar(\K)\nmid e$,
\begin{equation*}
	\rmHH^{2 k+1}_{(1)}(\Lambda)=
	 \begin{cases}
	 \K \left\{y_{kN+1, i }\di 0\le i\le N-1,\; i \equiv 0  \mode\right\}, &\text{if }  \rmchar (\K) \di N  \text{ and } kN+1  \equiv 0 \mode,\\
	 \K\left\{y_{kN+1, i}\di 1\le i\le N-1, i \equiv kN+1  \mode\right\}. &\text{otherwise};
	 \end{cases}
    \end{equation*}
    if $\rmchar(\K)\di e$, $\gcd(N-1,e)=c> 1$ and $\rmchar(\K)\nmid e_2:=\frac{e}{c}$, then 
    \[
    \rmHH^{2 k+1}_{(1)}(\Lambda)= \K \left\{ \sum_{j=0}^{e_2-1} (\alpha_{\ul{1+j(N-1)}}^{kN+1},\alpha_{\ul{1+j(N-1)}}^i) \,\middle|\, 1\le i\le N-1\text{ and }\;i \equiv kN+1  \mode \right\}.
    \]

\end{itemize}

\end{Cor}

Consequently, we explicitly construct a bijection between $\rmHH^n(\Lambda)$ and $\rmHH^n_{(1)}(\Lambda)$ for all $n\ge 1$. More precisely, 
\[
\rmHH^{2 k}(\Lambda)=\rmHH^{2 k}_{(1)}(\Lambda)\quad \text{for all }\; k\ge 1.
\]
If $\rmchar(\K)\nmid e$, then $\rmHH^{2 k+1}(\Lambda)=\rmHH^{2 k+1}_{(1)}(\Lambda)$ for all $k\ge 0$. If $\rmchar(\K)\di e$, $\gcd(N-1,e)=c> 1$ and $\rmchar(\K)\nmid e_2:=\frac{e}{c}$, then the element $v_{1;\; kN+1,j}\in \rmHH^{2k+1}(\Lambda)$ corresponds to $\frac{1}{e_2} \sum\limits_{j=0}^{e_2-1} (\alpha_{\ul{1+j(N-1)}}^{kN+1},\alpha_{\ul{1+j(N-1)}}^j)$.

\medskip

The following result provides the values of   $\Delta$ on the basis elements in $\rmHH^*(\Lambda)$.

\begin{Thm}\label{Thm: BV operation}

Let $\Lambda=\K Z_e/J^N$ be the truncated basic cycle algebra with $N,\, e\ge 2$, and the Nakayama automorphism $\sigma$ is semisimple. 

\begin{itemize}

    \item [(1)] Let $n=2h\ge 2$. For any $x \in \rmHH^n(\Lambda)$, $\Delta_n(x )=0$.

    \item[(2)] Let $n=2h+1\ge 1$ and $\rmchar(\K)$ satisfies condition (1) in Criterion~\ref{Criterion: semisimple Nakayama automorphism}, then for $y_{hN+1,j}$ with $0\le j\le N-1$ and $j\equiv hN+1\mode$ (Note that $j$ can take the value $0$ if and only if $\rmchar(\K)\di N$ and $hN+1\equiv 0\mode$), we have
		\begin{equation*}
		\Delta_n\big(y_{hN+1,j}\big)=
		\begin{cases}
		(hN+N-j)\,x_{hN,j-1}, &\text{if } 1\le j\le N-1,\\
		0,& \text{if } j=0.
		\end{cases}
		\end{equation*}

     \item[(3)] Let $n=2h+1\ge 1$ and $\rmchar(\K)$ satisfies condition (2) in Criterion~\ref{Criterion: semisimple Nakayama automorphism}. Assume that $\gcd(N-1,e)=c>1$, $N-1=c N_2$ and $e=c e_2$. Then for $v_{1;\;h N+1,j}$ with $1\le j\le N-1$ and $j\equiv h N+1\mode$, assume that $hN-j+1=h_1 e$, we have
		\begin{equation*}
		\Delta_n\big(v_{1;\; hN+1,j}\big)=(h_1+\frac{N_2}{e_2}) \, x_{hN,j-1}.
		\end{equation*}

\end{itemize}

\end{Thm}

\begin{Proof}

\begin{itemize}

\item [(1)] If $n=2h\ge 2$, let $x=x_{hN,j}$ with $0\le j\le N-1$ and $j\equiv h N\mode$ (Note that $j$ can take the value $N-1$ if and only if $\rmchar(\K)\di N$ and $hN\equiv N-1\mode$), then for $1\le q\le e$, we have
\[
\begin{aligned}
&\Delta_n(f_{x_{hN,j}})(1\ot \alpha_q^{(h-1)N+1}\ot 1)\\
=&(\mu_{2h-1}^*\circ \Delta_{\bbB} \circ\omega_{2h}^*)(f_{x_{hN,j}})(1\ot \alpha_q^{(h-1)N+1}\ot 1)\\
=&\Delta_{\bbB}(f_{x_{hN,j}}  \omega_{2h})\big(\mu_{2h-1}(1\ot \alpha_q^{(h-1)N+1}\ot 1)\big)\\
=&\sum\limits_{l=1}^e\sum\limits_{\substack{x_1,\cdots,x_{h-1}}}\Delta_{\bbB}(f_{(\alpha_l^{hN},\alpha_l^{j})}\omega_{2h})(1[p_1\mid \cdots\mid p_{2h-1}]\gamma)\\
=&\sum\limits_{l,l'=1}^e\sum\limits_{\substack{x_1,\cdots,x_{h-1}}} \sum\limits_{i=0}^{N-1} \sum\limits_{r=1}^{2h} (-1)^{r} \langle  f_{(\alpha_l^{hN},\alpha_l^j)}\omega_{2h}(1[p_r\mid \cdots\mid p_{2h-1}\mid \alpha_{l'}^i \mid \sigma(p_1)\mid \cdots \mid \sigma(p_{r-1})]1),1\rangle\alpha_{l''}^{N-1-i}\gamma,
\end{aligned}
\]
with $l''=\ul{l'+i-N+1}$. By definition of $\mu_{2h-1}$, we have $l( p_1) =l( p_3) =\cdots=l( p_{2h-1}) $, and $l( p_{2k}) =x_k\in [1,N-1]$ for $k=1,\ldots,h-1$.
By definitions of $\omega_{2h}$ and the operation $\langle- ,- \rangle$, the nonzero items in the above equation should satisfy $x_1=\cdots=x_{h-1}=i=j=N-1$ and $l'=q$, in this situation, $l''=q$, then the above equation can be written as
\begin{equation*}
	\begin{aligned}
	&\delta_{j,N-1}\sum\limits_{l=1}^e\sum\limits_{r=1}^{2h}(-1)^r\langle  f_{(\alpha_l^{hN},\alpha_l^{N-1})}\omega_{2h}(1[p_r\mid \cdots\mid p_{2h-1}\mid p_0\mid \sigma(p_1)\mid \cdots \mid \sigma(p_{r-1})]1),1\rangle e_q\\
	&=\delta_{j,N-1}\sum\limits_{l=1}^e\langle  f_{(\alpha_l^{hN},\alpha_l^{N-1})}(1\ot \alpha_q^{hN}\ot 1),1\rangle e_q &\textcolor{blue}{(r=2h)}\\
	& -\delta_{j,N-1}\sum\limits_{l=1}^e\langle  f_{(\alpha_l^{hN},\alpha_l^{N-1})}(1\ot \alpha_{q}^{hN}\ot 1),1\rangle e_q &\textcolor{blue}{(r=1)}\\
	& +\delta_{j,N-1}\sum\limits_{l=1}^e\sum\limits_{t=1}\limits^{ h-1}\langle  f_{(\alpha_l^{hN},\alpha_l^{N-1})}(1\ot \alpha_{\ul{q+(t-1)N+1}}^{hN}\ot 1),1\rangle e_q&\textcolor{blue}{(r=2t)}\\
	& -\delta_{j,N-1}\sum\limits_{l=1}^e\sum\limits_{t=1}\limits^{ h-1}\langle  f_{(\alpha_l^{hN},\alpha_l^{N-1})}(1\ot \alpha_{\ul{q+tN}}^{hN}\ot 1),1\rangle e_q&\textcolor{blue}{(r=2t+1)}\\
	&= \delta_{j,N-1}(h-1)e_q-\delta_{j,N-1}(h-1)e_q=0,
	\end{aligned}
\end{equation*}
so we have $\Delta_n\big(x_{hN,N-1}\big)=0$, (i) has been proved.

\item[(2)]  Let $n=2h+1$, $h\ge 0$ and $y_{hN+1,j}=\sum\limits_{l=1}^e (\alpha_l^{hN+1},\alpha_l^j)$ with $0\le j\le N-1$ ($j$ can be $0$ if and only if $\rmchar(\K)\di N$ and $hN+1\equiv 0\mode$) be a basis element of $\rmHH^n(\Lambda)$, for any $1\le q\le e$, we have
\begin{equation}\label{equation: delta1}
\begin{aligned}
&\Delta_n(f_{y_{hN+1,j}})(1\ot \alpha_q^{hN}\ot 1)\\
=&\Delta_{\bbB}(f_{y_{hN+1,j}}\omega_{2h+1})\big(\mu_{2h}(1\ot \alpha_q^{hN}\ot 1)\big)\\
=&\sum\limits_{l=1}^e\sum\limits_{\substack{x_1,\cdots,x_{h}}}\Delta_{\bbB}
(f_{(\alpha_l^{hN+1},\alpha_l^{j})}\omega_{2h+1})(1[p_1\mid \cdots\mid p_{2h}]\gamma)\\
=&\sum\limits_{l=1}^e\sum\limits_{\substack{x_1,\cdots,x_{h}}}\sum\limits_{i=0}^{N-1}\sum\limits_{l'=1}^{e}\sum\limits_{r=1}^{2h+1}\langle  f_{(\alpha_l^{hN+1},\alpha_l^{j})}\omega_{2h+1}(1[p_r\mid \cdots\mid p_{2h}\mid \alpha_{l'}^i\mid \sigma(p_1)\mid \cdots \mid \sigma(p_{r-1})]1),1\rangle\alpha_{l''}^{N-1-i}\gamma,
\end{aligned}
\end{equation}
with $l''=\ul{l'+i-N+1}$. By definition of $\mu_{2h}$, for $k=1,\ldots,h$, we have
\begin{equation*}
l( p_a) =
\begin{cases}
x_k,&\text{if } a=2k-1,\\
1,&\text{if } a=2k.
\end{cases}
\end{equation*}

The nonzero items in the equation~\ref{equation: delta1} with $r=2h+1$ should satisfy $x_1=\cdots=x_h=N-1$, $i=N-j$, $j\ge 1$ and $l'\equiv \ul{q+j-1} $, then $l''=q$, the item with $r=2h+1$ can be written as
\[
\begin{array}{ll}
&\sum\limits_{l=1}^e\langle  f_{(\alpha_l^{hN+1},\alpha_l^j)}\omega_{2h+1}(1[\alpha_{l'}^{N-j}\mid \sigma(p_1)\mid \cdots \mid \sigma(p_{2h})]1) ,1\rangle\alpha_q^{j-1}     \\
=& \sum\limits_{l=1}^e\sum\limits_{j'=0}^{N-j-1} \langle f_{(\alpha_l^{hN+1},\alpha_l^j)}(\alpha_{l'}^{j'}\ot \alpha_{\ul{l'+j'}}^{hN+1}\ot \alpha_{\ul{l'+j'+j}}^{N-j-j'-1}),1\rangle \alpha_q^{j-1} \\
=& \sum\limits_{l=1}^e a_{l,q}\; \alpha_q^{j-1},
\end{array}
\]
where $a_{l,q}$ is the cardinality of the set $\{ j'\di 0\le j'\le N-j-1 \text{ and } q+j+j'-1\equiv l\mode \}$.

The nonzero items in the equation~\ref{equation: delta1} with $r=2k-1$, $1\le k\le h$ should satisfy $x_1=\cdots=x_{k-1}=x_{k+1}=\cdots=x_h=i=N-1$, $x_k=N-j$, $j\ge 1$, $\gamma=\alpha_q^{j-1}$ and $l'= q$, then the item with $r=2k-1$ can be written as
\[
\begin{array}{ll}
     &\sum\limits_{l=1}^e\langle  f_{(\alpha_l^{hN+1},\alpha_l^j)}\omega_{2h+1}(1[p_{2k-1}\mid \cdots\mid p_{2h}\mid \alpha_q^{N-1} \mid \sigma(p_1)\mid \cdots \mid \sigma(p_{2k-2})]1) ,1\rangle\alpha_q^{j-1}  \\
     =&\sum\limits_{l=1}^e \sum\limits_{j'=0}^{N-j-1} \langle  f_{(\alpha_l^{hN+1},\alpha_l^j)}(\alpha_{\ul{q+(k-1)N}}^{j'} \ot \alpha_{\ul{q+(k-1)N+j'}}^{hN+1}\ot \alpha_{\ul{q+j+j'+(k-1)N}}^{N-j-j'-1}),1\rangle \alpha_q^{j-1} \\
     =&\sum\limits_{l=1}^e b_{k,l,q}\; \alpha_q^{j-1},
\end{array}
\]
where $b_{k,l,q}$ denotes the cardinality of the set $\{ j'\di 0\le j'\le N-j-1 \text{ and } q+j'+(k-1)N\equiv l\mode \}$.

The nonzero items in the equation~\ref{equation: delta1} with $r=2k$, $1\le k\le h$ should satisfy $x_2=\cdots=x_h=N-1$, $x_1+i=2N-1-j$, $j\ge 1$, $\alpha_{l''}^{N-1-i}\gamma=\alpha_q^{j-1}$ and $l'=\ul {q+x_1+j-N}$, then the item with $r=2k$ can be written as
\[
\begin{array}{ll}
     &\sum\limits_{l=1}^e\sum\limits_{x_1=N-j}\limits^{N-1}\langle  f_{(\alpha_l^{hN+1},\alpha_l^j)}\omega_{2h+1}(1[p_{2k }\mid \cdots\mid p_{2h}\mid \alpha_{l'}^{i}\mid \sigma(p_1)\mid \cdots \mid \sigma(p _{2k-1})]1) ,1\rangle\alpha_q^{j-1}  \\
     =&\sum\limits_{l=1}^e \sum\limits_{x_1=N-j}\limits^{N-1}\langle  f_{(\alpha_l^{hN+1},\alpha_l^j)}(1\ot \alpha_{\ul{q+(k-1)N+x_1}}^{hN+1}\ot \alpha_{\ul{q+(k-1)N+x_1+j}}^{N-j-1}),1 \rangle \alpha_q^{j-1} \\
     =&\sum\limits_{l=1}^e c_{k,l,q}\; \alpha_q^{j-1},
\end{array}
\]
where $c_{k,l,q}$ denotes the cardinality of the set $\{ j'\di N-j \le j'\le N-1 \text{ and } q+j'+(k-1)N\equiv l\mode \}$.

Then for any $1\le q\le e$, $(\mu_{2h}^*\circ \Delta_{\bbB} \circ\omega_{2h+1}^*)(f_{y_{ hN+1,j}})(1\ot \alpha_q^{hN}\ot 1)=\sum\limits_{l=1}^e d_{l,q}\, \alpha_q^{j-1}=(hN+N-j)\, \alpha_q^{j-1}$, with $d_{l,q}$ the cardinality of the set $\{ j'\di 0 \le j'\le hN+N-j-1 \text{ and } q+j'\equiv l\mode \}$.

Then we obtain that $\Delta_{2h+1}(y_{hN+1,j})=(hN+N-j)\, x_{hN,j-1}$ if $j\ge 1$; otherwise $\Delta_{2h+1}(y_{hN+1,j})=0$.

\item[(3)] If $n=2h+1\ge 1$ and $\rmchar(\K)$ satisfies condition (2) in Criterion~\ref{Criterion: semisimple Nakayama automorphism}, then the basis element of $\rmHH^{2h+1}(\Lambda)$ is of the form $v_{1;\;h N+1,j}$ with $1\le j\le N-1$ and $j\equiv hN+1\mode$. We have 
\[
\Delta_n(v_{1;\; hN+1,j})=\frac{1}{e_2}\Delta_n\big( \sum\limits_{r=0}^{e_2-1} (\alpha_{\ul{1+r(N-1)}}^{hN+1},\alpha_{\ul{1+r(N-1)}}^j) \big) .
\]

By the proof of (ii), we have 
\[
\Delta_n\big( \sum\limits_{r=0}^{e_2-1} (\alpha_{\ul{1+r(N-1)}}^{hN+1},\alpha_{\ul{1+r(N-1)}}^j) \big)= \sum_{q=1}^e \sum\limits_{r=0}^{e_2-1} d_{r,q} (\alpha_q^{hN},\alpha_q^{j-1}),
\]
where $d_{r,q}$ is the cardinality of the set $\{ j'\di 0 \le j'\le hN+N-j-1,\; q+j'\equiv 1+r(N-1)\mode \}$.
 Then $\sum\limits_{r=0}^{e_2-1} d_{r,q}=h_1e_2+N_2 $. Therefore, $\Delta_n(v_{1;\; hN+1,j})=(h_1+\frac{N_2}{e_2}) \, x_{hN,j-1} $.
\end{itemize}
	
\end{Proof}

\begin{Rem}In case of semisimple Nakayama automorphism, our formula for the $\Delta$-operator coincides with that in \cite[Theorem 4.6]{Ita23}, except that in the statement of \cite[Theorem 4.6]{Ita23}, there is a typos that ``$Ni+N-j-1$" should read ``$\frac{Ni+N-j-1}{g_0}$", as the correct formula is given in the proof. 
\end{Rem}

Now we turn to the case where the Nakayama automorphism $\sigma$ is not semisimple, that is, when $\rmchar(\K)\di \mathrm{ord}(\sigma)$. In this case, the Hochschild cohomology $\rmHH^*(\Lambda)$ remains a BV algebra.

Before presenting the explicit construction of the BV differential, we first state a more general criterion.
\begin{Cri}\label{Criterion: nonsemisimple case}
	Let $(\rmH^{\bullet},\vee,[-,-])$ be a Gerstenhaber algebra with $(\rmH^{\bullet},\vee)$ finitely generated and  $\mathcal{G}$ a  minimal set of  homogeneous generators of $(\rmH^{\bullet},\vee)$.  Assume that the following assumptions hold: 
	\begin{itemize}
		\item [(1)] $\mathcal{G}\cap \rmH^{2\bullet-1}\subset  \rmH^1$, and for any $y, y'\in \mathcal{G}\cap \rmH^{2\bullet-1}$,  $y\vee y'=0$  and $[y,\,y']=0$.
		\item [(2)]Given  $\lambda_i\in\K$, $y_i\in \mathcal{G}\cap \rmH^{2\bullet-1}$,  and $f_i^{\even}$   being  products of some elements in $(\mathcal{G}\cap \rmH^{2\bullet}) \cup \{1\}$,  whenever $\sum_{i=1}^n \lambda_i\, f_i^{\even}\vee y_i=0$ holds in $\rmH^\bullet$, we have
			\[
			\sum_{i=1}^n k_i \,[f_i^{\even},\, y_i]=0.
			\]
		  
		\item [(3)] For any $f,\, g\in \rmH^{2\bullet}$, 
		\[
		[f,\,g]=0.
		\]
	\end{itemize} 
	Then $(\rmH^{\bullet},\vee,[-,-])$ becomes a BV algebra endowed with the operator $\Delta: \rmH^{\bullet}\to \rmH^{\bullet-1}$: for $n\ge0 $,
	\[\begin{aligned}
		&\Delta_{2n}=0,\\
		&\Delta_{2n+1}(f^{\even} \vee y) =[f^{\even}, \, y],
	\end{aligned}
	\]
	 where $y\in \mathcal{G}\cap \rmH^{2\bullet-1}$, $f^{\even}$ is a product of some elements in $(\mathcal{G}\cap \rmH^{2\bullet}) \cup \{1\}$.
\end{Cri}

\begin{Proof}
	 By Assumption (2), $\Delta$ is well-defined. 

It is clear that $\Delta\circ \Delta=0$ by definition.

	For $f,\,g\in \rmH^{2\bullet}$, by the definition of $\Delta$ and Assumption (3), we have
	\[
	[f,\, g]=0=(-1)^{|f| }(\Delta(f\vee g)- \Delta(f)\vee g-(-1)^{|f| }f\vee \Delta(g)).
	\] 

	Let us consider  $f=\sum_{i=1}^{n}\lambda_i\, f_i^{\even} \in \rmH^{2n}$ and $g=\sum_{j=1}^{m} \mu_j\, g_j^{\even}\vee z_j\in\rmH^{2m+1}$,  
	where $\lambda_i, \mu_j\in \K$, $f_i^{\even}$, $g_j^{\even}$ are products of some elements from $(\mathcal{G}\cap \rmH^{2\bullet}) \cup \{1\}$ and $z_j\in \mathcal{G}\cap \rmH^{2\bullet-1}$. By the definition of $\Delta$ and Assumption (1), we have
	\[
	\begin{aligned}
		&(-1)^{|f|}\big(\Delta(f\vee g)-\Delta(f)\vee g-(-1)^{|f|} f\vee \Delta(g)\big)\\
		&= \Delta(f\vee g)-f\vee \Delta(g)\\
		&=\sum_{i=1}^{n} \sum_{j=1}^{m} \lambda_i \mu_j\, \big( \Delta(f_i^{\even}\vee g_j^{\even}\vee z_j) -f_i^{\even}\vee \Delta(g_j^{\even}\vee z_j) \big)\\
		&= \sum_{i=1}^{n} \sum_{j=1}^{m} \lambda_i \mu_j\, \big( [f_i^{\even}\vee g_j^{\even}, z_j] -f_i^{\even}\vee [g_j^{\even}, z_j] \big)\\
		&= \sum_{j=1}^{m} \mu_j\, [f,z_j]\vee g_j^{\even} \quad (\text{by Poisson rule for $\vee$ and $[-,-]$})\\
		&=- \sum_{j=1}^{m} \mu_j\,  g_j^{\even}\vee [z_j,f] \\
		&=[f,g] \quad (\text{by Poisson rule and Assumption (3)}).
	\end{aligned}
	\]

It is obvious by Assumption (1) that for any $f,g\in\rmH^{2\bullet-1}$, we have 
	\[
	f\vee g=0.
	\]
Let us consider  $f=\sum_{i=1}^{n}\lambda_i\, f_i^{\even}\vee y_i \in \rmH^{2n+1}$ and $g=\sum_{j=1}^{m} \mu_j\, g_j^{\even}\vee z_j\in\rmH^{2m+1}$,  
	where $\lambda_i,\mu_j\in \K$, $  f_i^{\even}$, $g_j^{\even}$ are products of  elements from $(\mathcal{G}\cap \rmH^{2\bullet}) \cup \{1\}$ and $y_i, z_j\in \mathcal{G}\cap \rmH^{2\bullet-1}$.   By the definition of $\Delta$ and Assumption (1), we have
	\[
	\begin{aligned}
		&(-1)^{|f|}\big(\Delta(f\vee g)-\Delta(f)\vee g-(-1)^{|f|} f\vee \Delta(g)\big)\\
		&= \Delta(f)\vee g-f\vee \Delta(g)\\
		&=-\sum_{i=1}^{n} \lambda_i \, [y_i, f_i^{\even}]\vee g + \sum_{j=1}^{m} \mu_j\, f\vee [z_j, g_j^{\even}]\\
		&=-\sum_{i=1}^{n} \lambda_i \, y_i\vee [g ,f_i^{\even}] + \sum_{j=1}^{m} \mu_j\, [f, g_j^{\even}]\vee z_j \quad (\text{by Poisson rule and Assumption (1)})\\
		&=-\sum_{i=1}^{n} \sum_{j=1}^{m} \lambda_i \mu_j\,  y_i\vee [g_j^{\even}\vee z_j ,f_i^{\even}] + \sum_{j=1}^{m} \mu_j\, [f, g_j^{\even}]\vee z_j\\
		&=-\sum_{i=1}^{n} \sum_{j=1}^{m} \lambda_i \mu_j\,  y_i\vee g_j^{\even}\vee  [z_j ,f_i^{\even}] + \sum_{j=1}^{m} \mu_j\, [f, g_j^{\even}]\vee z_j \quad (\text{by Poisson rule and Assumption (3)}) \\
		&=\sum_{i=1}^{n} \sum_{j=1}^{m} \lambda_i \mu_j\,   g_j^{\even}\vee y_i \vee  [f_i^{\even}, z_j] + \sum_{j=1}^{m} \mu_j\, [f, g_j^{\even}]\vee z_j\\
		&=\sum_{i=1}^{n} \sum_{j=1}^{m} \lambda_i \mu_j\, g_j^{\even}\vee  [y_i \vee f_i^{\even}, z_j] + \sum_{j=1}^{m} \mu_j\, [f, g_j^{\even}]\vee z_j \quad (\text{by Poisson rule and Assumption (1)})\\
		&= \sum_{j=1}^{m} \mu_j\, g_j^{\even}\vee  [f, z_j] + \sum_{j=1}^{m} \mu_j\, [f, g_j^{\even}]\vee z_j \\
		&=[f,g] \quad (\text{by Poisson rule}).
	\end{aligned}
	\]
	Therefore, $\Delta$ is a BV differential on $(\rmH^{\bullet}, \vee, [-,-])$.
\end{Proof}

 By Theorem~\ref{Thm: basis}, Proposition~\ref{Prop: cup product of basis elements}, and Proposition~\ref{Thm:Gerstenhaber bracket of basic elements}, when the Nakayama automorphism $\sigma$ is not semisimple (i.e., $\rmchar(\K)\mid \operatorname{ord}\,(\sigma)$), the Gerstenhaber algebra $(\rmHH^*(\Lambda),\vee,[-,-])$ satisfies the assumptions of Criterion~\ref{Criterion: nonsemisimple case}. We therefore obtain the following result. 

\begin{Thm}\label{thm:delta-operation-nonsemisimple}
Let $\Lambda=\K Z_e/J^N$ be the truncated basic cycle algebra, and the Nakayama automorphism $\sigma$ is not semisimple. Then, $(\rmHH^*(\Lambda),\vee,[-,-])$ is a BV algebra. More precisely, the $\Delta$-operator is given as follows: 
\begin{itemize}

    \item [(1)] Let $n=2h\ge 2$. For any $x\in \rmHH^n(\Lambda)$, $\Delta_n(x)=0$.

     \item[(2)] Let $n=2h+1\ge 1$. For $v_{1;\;h N+1,j}$ with $1\le j\le N-1$ and $j\equiv h N+1\mode$, assume that $hN-j+1=h_1 e$, we have
		\begin{equation*}
		\Delta_n\big(v_{1;\; hN+1,j}\big)= h_1 \, x_{hN,j-1}.
		\end{equation*}
    
\end{itemize}

\end{Thm}

\begin{Rem}
Note that the BV operator $\Delta$ in Theorem~\ref{thm:delta-operation-nonsemisimple} is not unique. For instance, in case~(2), for any $0\neq a\in \K$, if we define
\[
\Delta_n\big(v_{1;\; hN+1,j}\big)= (h_1 +a)\, x_{hN,j-1}
\]
and keep values of $\Delta$ on even-degree elements unchanged, then this defines another BV algebra structure on $(\rmHH^*(\Lambda),\vee,[-,-])$.
\end{Rem}

The authors of \cite{LZZ16} suggested the following example of  a self-injective Nakayama algebra with non diagonalisable Nakayama automorphism,  which may be a possible example of a Frobenius algebra whose Hochschild cohomology ring is NOT  a BV algebra. However, this is not the case.

\begin{Ex}

Let $\K$ be a field of characteristic two. Consider the algebra $\Lambda$ defined by the quiver with relations
\[
\xymatrix{1\ar@<1ex>[r]^\alpha & 2\ar@<1ex>[l]^{\beta},  & \alpha\beta\alpha\beta=0=\beta\alpha\beta\alpha}.
\]
As described in Example~\ref{Example: generators and relations}, $\rmHH^*(\Lambda)\cong \K[x,y,z]/\langle x^2,y^2\rangle$ with $|x| =0$, $|y| =1$ and $|z| =2$. And by Example~\ref{Example: Gerstenhaber brackets}, $\rmHH^*(\Lambda)$ is a Gerstenhaber algebra with the Gerstenhaber brackets between generators are given by:
\[
[x,y]=-x,\; [x,x]=[x,z]=[y,y]=[y,z]=[z,z]=0.
\]
The BV differential $\Delta$ is given by
\[
\Delta(xy)=-x,\; \Delta(x)=\Delta(y)=\Delta(z)=\Delta(xz)=\Delta(yz)=\Delta(z^2)=0.
\]
\end{Ex}

 Combining Theorems~\ref{Thm: BV operation} and \ref{thm:delta-operation-nonsemisimple} yields the   main result of this paper.
\begin{Thm}\label{thm:main-result}
	Let $\Lambda=\K Z_e/ J^N$, $N\ge 2$, be a truncated basic cycle algebra. The the Hochschild cohomology ring $\rmHH^*(\Lambda)$ of $\Lambda$ is a BV algebra.
\end{Thm}

\bigskip

\textbf{Acknowledgements:} The first, third and fifth authors were supported by  the National Key R$\&$D Program of China (No. 2024YFA1013803), by the China Scholarship Council (CSC),  and by Shanghai Key Laboratory of PMMP (No. 22DZ2229014). The fourth author was supported by the National Natural Science Foundation of China (No. 12401039).

It is our great pleasure to express our sincere gratitude to  Satoshi Usui for his very useful comments.

\end{document}